
\def\LaTeX{%
  \let\Begin\begin
  \let\End\end
  \let\salta\relax
  \let\finqui\relax
  \let\futuro\relax}

\def\UK{\def\our{our}\let\sz s}
\def\USA{\def\our{or}\let\sz z}

\UK 



\LaTeX

\USA


\salta

\documentclass[twoside,12pt]{article}
\setlength{\textheight}{24cm}
\setlength{\textwidth}{16cm}
\setlength{\oddsidemargin}{2mm}
\setlength{\evensidemargin}{2mm}
\setlength{\topmargin}{-15mm}
\parskip2mm


\usepackage[usenames,dvipsnames]{color}
\usepackage{amsmath}
\usepackage{amsthm}
\usepackage{amssymb,bbm}
\usepackage[mathcal]{euscript}

\usepackage{cite}
\usepackage{enumitem}

\usepackage[ulem=normalem,draft]{changes}
%
%

%
 
\definecolor{ciclamino}{rgb}{0.5,0,0.5}
\definecolor{rosso}{rgb}{0.85,0,0}

\def\abramo #1{{\color{black}#1}}
\def\abramob #1{{\color{violet}#1}}
\def\an #1{{\color{rosso}#1}}
\def\last #1{{\color{black}#1}}
\def\andre #1{{\color{black}#1}}

\def\an #1{{#1}}
\def\abramo #1{{#1}}
\def\abramob #1{{#1}}




\bibliographystyle{plain}


%
\newtheorem{theorem}{Theorem}[section]
\newtheorem{remark}[theorem]{Remark}

\newtheorem{lemma}[theorem]{Lemma}

\finqui

\def\Bcenter{\Begin{center}}
\def\Ecenter{\End{center}}
\let\non\nonumber




\def\step #1 \par{\medskip\noindent{\bf #1.}\quad}
\def\jstep #1: \par {\vspace{2mm}\noindent\underline{\sc #1 :}\par\nobreak\vspace{1mm}\noindent}

\def\Lip{Lip\-schitz}
\def\Holder{H\"older}

\def\Poincare{Poincar\'e}
\def\lhs{left-hand side}
\def\rhs{right-hand side}




\def\multibold #1{\def\arg{#1}%
  \ifx\arg\pto \let\next\relax
  \else
  \def\next{\expandafter
    \def\csname #1#1\endcsname{{\boldsymbol #1}}%
    \multibold}%
  \fi \next}

\def\pto{.}

\def\multical #1{\def\arg{#1}%
  \ifx\arg\pto \let\next\relax
  \else
  \def\next{\expandafter
    \def\csname cal#1\endcsname{{\cal #1}}%
    \multical}%
  \fi \next}

\def\multigrass #1{\def\arg{#1}%
  \ifx\arg\pto \let\next\relax
  \else
  \def\next{\expandafter
    \def\csname grass#1\endcsname{{\mathbb #1}}%
    \multigrass}%
  \fi \next}


\def\multimathop #1 {\def\arg{#1}%
  \ifx\arg\pto \let\next\relax
  \else
  \def\next{\expandafter
    \def\csname #1\endcsname{\mathop{\rm #1}\nolimits}%
    \multimathop}%
  \fi \next}

\multibold
qweryuiopasdfghjklzxcvbnmQWERTYUIOPASDFGHJKLZXCVBNM.  

\multical
QWERTYUIOPASDFGHJKLZXCVBNM.

\multigrass
QWERTYUIOPASDFGHJKLZXCVBNM.

\multimathop
diag dist div dom mean meas sign supp .


\def\Accorpa #1#2 #3 {\gdef #1{\eqref{#2}--\eqref{#3}}%
  \wlog{}\wlog{\string #1 -> #2 - #3}\wlog{}}


\def\graffe #1{\mathopen\{#1\mathclose\}}

\def\<#1>{\mathopen\langle #1\mathclose\rangle}
\def\norma #1{\mathopen \| #1\mathclose \|}

\let\ov\overline

\def\cpto{\,\cdot\,}

\def\Beta{\widehat\beta}
\def\Pi{\widehat\pi}

\def\Betaeps{\Beta_\eps}
\def\betaeps{\beta_\eps}

\def\iot {\int_0^t}

\def\intQt{\int_{Q_t}}
\def\intQ{\int_Q}
\def\iO{\int_\Omega}

\def\dt{\partial_t}
\def\dn{\partial_{\nn}}

\def\checkmmode #1{\relax\ifmmode\hbox{#1}\else{#1}\fi}


\let\erre\grassR
\let\enne\grassN
\def\errebar{(-\infty,+\infty]}




\def\genspazio #1#2#3#4#5{#1^{#2}(#5,#4;#3)}
\def\spazio #1#2#3{\genspazio {#1}{#2}{#3}T0}
\def\spaziot #1#2#3{\genspazio {#1}{#2}{#3}t0}

\def\L {\spazio L}
\def\H {\spazio H}
\def\W {\spazio W}
\def\Lt {\spaziot L}
\def\Ht {\spaziot H}

\def\C #1#2{C^{#1}([0,T];#2)}


\def\Lx #1{L^{#1}(\Omega)}
\def\Hx #1{H^{#1}(\Omega)}
\def\Wx #1{W^{#1}(\Omega)}



\let\badtheta\theta
\let\theta\vartheta
\let\badeps\epsilon
\let\eps\varepsilon

\let\phi\varphi

\let\TeXchi\chi                         
\newbox\chibox
\setbox0 \hbox{\mathsurround0pt $\TeXchi$}
\setbox\chibox \hbox{\raise\dp0 \box 0 }
\def\chi{\copy\chibox}



\let\hat\widehat

\def\0{{\boldsymbol {0} }}

\let\I\grassI

\def\NN{{\cal N}}

\let\emb\hookrightarrow

\def\Vp{{V^*}}


\def\CP0{(${\mathcal{CP}}_0$)}

\def\d {{\delta}}

\def\b {{\beta}}
\def\ph {{\varphi}}
\def\phv {{\ph}}
\def\phd {{\ph}}
\def\pha {{\ph_a}}
\def\pham {{\ph_{a,\gal}}}
\def\phT {{\ph}}

\def\Sv {{\cal S}_v}

\def\Sn {{\cal S}_n}
\def\Sv {{\cal S}}
\def\SS {{\cal S}}
\def\Sc {{\cal S}_c}
\def\Sa {{\cal S}_a}

\def\mobm{{\mathbbm{m}}}
\def\mobn{{\mathbbm{n}}}

\def\Wn{{\cal W}_d}
\def\vv{{\boldsymbol{v}}}

\def\rate{{{R_d}}}
\def\esp{{{\sigma}}}
\def\gal{{\an{k}}}
\let\epsilon\eps
\def\gtest{{\an{\psi}}}

%

\Begin{document}


%
\title{
{Analysis of a} multi-species Cahn--Hilliard--Keller--Segel tumor growth model with chemotaxis and angiogenesis}

\author{}
\date{}
\maketitle
\Bcenter
\vskip-1.5cm
{\large\sc Abramo Agosti$^{(1)}$}\\
{\normalsize e-mail: {\tt abramo.agosti@unipv.it}}
\\[0.25cm]
{\large\sc Andrea Signori$^{(2)}$}\\
{\normalsize e-mail: {\tt andrea.signori@polimi.it}}
\\[0.25cm]
$^{(1)}$
{\small Dipartimento di Matematica ``F. Casorati''}\\
{\small Universit\`a di Pavia}\\
{\small via Ferrata 5, I-27100 Pavia, Italy}\\[.3cm] 
$^{(2)}$
{\small Dipartimento di Matematica, Politecnico di Milano}\\
{\small via E. Bonardi 9, I-20133 Milano, Italy}


\Ecenter

\Begin{abstract}
\noindent 
{We introduce a multi-species diffuse interface model for tumor growth, characterized by its incorporation of essential features related to chemotaxis, angiogenesis {and proliferation} mechanisms. 
We establish the weak well-posedness of the system within an appropriate variational framework, accommodating various choices for the nonlinear potentials.}
One of the primary novelties of the work lies in the rigorous establishment of the existence of a weak solution through the introduction of delicate approximation schemes. {To our knowledge,  this represents a novel advancement for both the intricate Cahn--Hilliard--Keller--Segel system and the Keller--Segel subsystem with source terms.}
Moreover,
{when specific conditions are met, such as having more regular initial data, \abramob{a smallness condition on the chemotactic constant with respect to the magnitude of initial conditions} and potentially focusing solely on the two-dimensional case,}
we provide regularity results for the weak solutions.
{Finally, we derive a continuous dependence estimate, which, in turn, leads to the uniqueness of the smoothed solution as a natural consequence.}

\End{abstract}
\salta
\pagestyle{myheadings}
\newcommand\testopari{\sc Agosti -- Signori }
\newcommand\testodispari{\sc Multiphase Cahn--Hilliard--Keller--Segel tumor model}
\markboth{\testopari}{\testodispari}
\finqui
%

\section{Introduction}
\label{SEC:INTRO}
\setcounter{equation}{0}
Suppose $\Omega \subset \erre^d$, $d \in \{2,3\}$, be a smooth and bounded domain,
and let $T>0$ be a given final time. 
Then, we introduce and analyze a  tumor growth model {that characterizes the dynamics of tumor development in the presence of chemotaxis and angiogenesis.}
The system assumes the form of a {\it multiphase Cahn--Hilliard--Keller--Segel} (MCHKS) model and it reads as follows:
\begin{alignat}{2}
	\label{system:1}
	& \dt \phv - \div\big(\mobm(\phT,\pha,n)(\nabla \mu -  \chi_\ph \nabla n)\big) = \Sv 
	\quad && \text{in ${Q}$,}
	\\
	\label{system:3}
	& \mu = - \Delta \phT + F'(\phT) 
	\quad && \text{in $Q$,}
	\\
	\label{system:4}
	& \dt \pha - \div\big(\pha \mobn(\pha,c)\nabla ({\log (\pha)} - \chi_a c)\big) = \Sa
	\quad && \text{in $Q$,}
	\\
	\label{system:5}
	& \dt n - \Delta n - \chi_\ph \phv = \Sn 
	\quad && \text{in $Q$,}
	\\
	\label{system:6}
	& \dt c - \Delta c  - \chi_a \pha = \Sc 
	\quad && \text{in $Q$,}	
\end{alignat}
where $Q:= \Omega \times (0,T)$, {$\mobm,\mobn$ are positive mobility functions, $F'$ is a local interaction potential} and $\Sv,\Sa,\Sn,$ and $ \Sc$ indicate some source terms accounting for the interplay between the different variables
{$\ph,\mu,\pha,n,$ and $c$ whose meaning will be discussed in the section to follow.}
Meaningful biological {for the source terms} examples will be {depicted} in the next sections.
{Besides, we anticipate that the magnitude of the chemotaxis sensitivities $\chi_\ph$ and $\chi_a$ will play a role in the forthcoming mathematical analysis, and will be chosen accordingly.}
As for the initial and boundary conditions,
{after setting $\Sigma:= \partial \Omega \times (0,T),$}
we require that
\begin{alignat}{2}
	\non
	\dn \phT & = (\mobm(\phT,\pha,n)\nabla \mu) \cdot \nn= 	(\pha \mobn(\ph_a,c)\nabla ({\log (\pha)} - \chi_a c)) \cdot \nn
	\\  \quad 	
	\label{system:7}
	& = \dn n =  \dn c =0 	\qquad && \text{on ${\Sigma}$,}
	\\
	\label{system:8}
	\phT(0) & = \ph^0,
	\quad 
	\pha(0) = \ph_a^0,
	\quad
	n(0) = n^0,
	\quad 
	c(0) = c^0
	\qquad && \text{in $\Omega$\an{,}}
\end{alignat}
\Accorpa\Sys {system:1} {system:8}
{with $\nn$ indicating the outer unit normal vector to $\partial\Omega$.}

Notice that as $\mobn \equiv 1$, at least for sufficiently regular solutions, equation \eqref{system:4} can be rewritten in the following equivalent forms:
\begin{align}
	\non
	&\dt \pha - \div(\pha \nabla ({\log (\pha)} - \chi_a c)) 
	= \dt \pha - \Delta \pha + \chi_a \div( \pha \nabla c ) 
	\\ & \quad 
	\label{system:4bis}
	= \dt \pha - \div(\nabla \pha - \chi_a \nabla c) =  \Sa.
\end{align}
{In the following, after deriving the model \Sys\ from basic principles of Mixture Theory and variational principles of Thermodynamics, we will {suggest} relevant biological constitutive assumptions for the {source terms in the system}. Then, {combining} a Faedo--Galerkin approximation {scheme along with further regularizations for the occurring nonlinearities}, we will prove, under proper assumptions on the regularity of initial data and on the growth laws {of} the source terms, the  existence of {global} weak solutions.
{Ultimately}, we will provide regularity results for the weak solutions, together with their continuous dependence from data, under {stronger} assumptions on initial data and in the two-dimensional {setting}, which lead to the uniqueness of the solution.}

\subsection{Modeling Considerations and Relevant  Biological  Choices}
\label{SUB:BIO}
In this section, we derive the {multiphase} model for tumor growth with {angiogenesis and chemotaxis} in \eqref{system:1}--\eqref{system:6}{.
The {above} model describes} a diffuse-interface mixture {composed} by a tumor component, a liquid component and an angiogenetic component, coupled with two massless chemicals representing a perfectly diluted nutrient and an angiogenetic factor. {The system} is a three-phase reduction of the multiphase model introduced and calibrated on patient-specific data in \cite{ALL}, which considered the tumor phase as composed by a viable and a necrotic components, typically observed in neuroimaging data, which exchange mass in hypoxic conditions. Since here we are interested in obtaining analytical results for the proposed model, we reduce its complexity by considering a single tumor component, hence neglecting the dynamics of necrosis formation. This does not modify the mathematical structure of the underlying PDEs system, which, as in the model proposed in \cite{ALL}, takes the form of a Cahn--Hilliard{--}Keller–Segel system for the mixture components with reaction-diffusion equations for the chemicals; yet, in the present work we will state and consider {more} general constitutive assumptions {than the ones considered in \cite{ALL}, identifying the {general conditions}} which will let us obtain analytical results regarding existence and regularity of solutions. 

The model \eqref{system:1}--\eqref{system:6} is derived from variational principles complying with the second law of thermodynamics in isothermal situations. In the following, we will only present the main steps {of the derivation}, referring to \cite{ALL} for {more} {details}. Let us consider a saturated, closed and incompressible mixture in {$\Omega$}, composed by a tumor phase with volume fraction $\phi$, a liquid phase composed by liquid, healthy cells and normal vasculature, with volume fraction $\phi_l$, and an angiogenetic phase composed by tumor-induced new vasculature with volume fraction $\phi_a$. We assume that all the phases have a constant density $\gamma$, equal to the water density (since the cells are mostly composed by water). The mixture dynamics is coupled with the evolution of massless chemicals, comprising a nutrient species, with concentration (number of moles) $n$, and an angiogenetic factor, with concentration (number of moles) $c$. Each mixture component satisfies a mass continuity equation, while the massless nutrient and chemical species satisfy generic transport equations:
\begin{align}
\label{eqn:mix1a}
    {\dt \phi}
    +\div(\phi \vv)+\div(\mathbf{J_{\phi}})=\frac{\Gamma_{\phi}(\phi,\phi_l,{\pha,}n,c)}{\gamma},\\
    \label{eqn:mix1b}
    {\dt \pha}
    +\div(\phi_a \vv)+\div({\mathbf{J}_a})
    =\frac{\Gamma_a(\phi,\phi_l,{\pha,}n,c)}{\gamma},
    \\
    \label{eqn:mix1c}
    {\dt \phi_l}
    +\div(\phi_l \vv)+\div({\mathbf{J}_l})=\frac{\Gamma_l(\phi,\phi_l,{\pha,}n,c)}{\gamma},\\
    \label{eqn:mix1d}
    {\dt n}
    +\div(n \vv)+F_n={\Sn}(\phi,\phi_l,{\pha,}n,c),\\
    \label{eqn:mix1e}
    {\dt c}
    +\div(c \vv)+F_c={\Sc}(\phi,\phi_l,{\pha,}n,c),
\end{align}
{subject to the constraints} $\phi+\phi_a+\phi_l=1$ {and} $\Gamma_{\phi}+\Gamma_a+\Gamma_l=0$, {with} {fluxes} $\mathbf{J_{\phi}}=\phi(\vv_{\phi}-\vv)$, ${\mathbf{J}_a}=\phi_a(\vv_a-\vv)$, {and} ${\mathbf{J}_l}=\phi_l(\vv_l-\vv)$. Here, {$\Gamma_\ph,\Gamma_a,$ and $\Gamma_l$ stand for source terms,} $\vv=\phi \vv_{\phi}+\phi_a \vv_a+\phi_l \vv_l$ is the volume-averaged mixture velocity, which satisfies the incompressibility condition 
\begin{equation}
\label{eqn:mix2}
\div\vv=0,
\end{equation}
as a consequence of the saturation and the closedness properties of the mixture. 
The terms $-F_n$ and $-F_c$ are generic transport terms to be determined in relation with the specific free energy of the system, while the source terms ${\Sn}$ and ${\Sc}$ represent source and consumption terms for the chemicals and must be constitutively assigned. 
{On the other hand, the source}  terms $\Gamma_{\phi}$ {and} $\Gamma_a$ represent cells proliferation and death, while we take $\Gamma_l{:}=-\left(\Gamma_{\phi}+\Gamma_a\right)$. We make the following modeling assumptions:
\begin{itemize}
\item The endothelial cells of the tumor-induced vasculature constitute a self-interacting phase in the mixture which can migrate to regions with higher angiogenetic factor concentration, being coupled to the angiogenetic factor by a chemotactic term;
\item The tumor cells can migrate to regions with higher nutrient concentration, being coupled to the nutrient by a chemotactic term.
\end{itemize}
\begin{remark}
In \cite{GLSS}{,} a prototype tumor growth model with nutrient diffusion and chemotaxis, further developed and studied in multiple subsequent works (see, e.g., \cite{EG,EGN,GKT}), was introduced{. There,} the nutrient dynamics was constrained as a modeling assumption to satisfy a mass continuity equation. This lead to a cross-diffusion term in the nutrient equation, representing a nutrient flux towards regions with higher cells concentration, which, as noted in \cite{RSchS} {(see also \cite{GSigSpr})}, may have an {nonphysical} interpretation. Moreover, with this term the nutrient equation does not satisfy the minimum and maximum principles, and the nutrient concentration may assume {nonphysical} values. In \cite{RSchS}, a different modeling approach was employed by assuming that the nutrient dynamics satisfy a mass continuity equation in the form of a Keller--Segel equation, coupled with a Cahn--Hilliard equation for the tumor concentration. The latter approach makes the nutrient flux towards regions with higher cells concentration proportional to the nutrient concentration, representing chemotactic aggregation of nutrients following the cells gradient, and also enforces a minimum principle for the nutrient. We observe that this picture of the nutrient dynamics is non-standard, since typically the dynamics of massless and passive chemicals are driven by random motion, with no self-aggregation, plus advection and source terms. Hence, in our modeling approach we constrain the massles chemicals to satisfy generic transport equations, in the form of reaction-advection-diffusion equations which satisfy both minimum and maximum principles, while we enforce the mass continuity equation for the mixture components which contribute to the mixture mass.
\end{remark}
With the given modeling assumptions, the free energy of the system, expressed as its internal energy minus its entropy, takes the following general {structure}
\begin{align}
    \non
    & 
    E(\phi,{\pha},n,c)=\int_{\Omega}e(\phi,{\pha},n,c)
    {= \int_{\Omega}{\Lambda} \big((F(\phi)+{\kappa_a}\phi_a(\log(\phi_a)-1) \big)}
    \\ 
    & \quad \label{eqn:mix3}
    + {\iO {\Lambda} \biggl(} \frac{{{\badeps}}^2}{2}|\nabla \phi|^2
    +\frac{{D_n}}{2}|\nabla n|^2-\chi_{\phi}n\phi+\frac{{D_c}}{2}|\nabla c|^2-\chi_a\phi_ac\biggr),
\end{align}
where $e(\phi,{\pha},n,c)$ is the free energy {density} per unit volume. Here, ${\Lambda}$ is the Young modulus of the tissue, in units of $[Pa]$. The term $-{\kappa_a}\phi_a(\log(\phi_a)-1)$ is the entropy associated to the self-interacting endothelial cells, with ${\kappa_a}$ a positive adimensional coefficient. The term $\frac{{{\badeps}}^2}{2}|\nabla \phi|^2$ represents the diffuse-interface internal energy between the tumor cells and the host tissue, with ${{\badeps}}$, in units of $[m]$ {denoting} the interfacial thickness. 
The terms $\frac{{D_n}}{2}|\nabla n|^2$ and $\frac{{D_c}}{2}|\nabla c|^2$ represent the contribution to the internal energy from the random motion of the chemical species resulting in diffusive behaviors along concentration gradients, where ${D_n}$ and ${D_c}$ are the isotropic mean deviations of the displacement of the particles, with units of $[mm^2/\text{Mol}^2]$. The terms $-\chi_{\phi}n\phi$ and $-\chi_a\phi_ac$ are interaction terms associated to chemotaxis, with positive chemotactic coefficients $\chi_{\phi}$ and $\chi_a$ in units of $[\text{Mol}^{-1}]$. Finally, the term $F(\phi)$ represents the entropy minus the internal energy associated to the binary interaction between the tumor cells and their surrounding. {A typical choice is represented by the \textit{Flory--Huggins} potential which} can be written as
\begin{equation}
    F_{\rm dw}({r}):=
  \begin{cases}
  		\frac{c_1}{2}\left({r}\log {r}+(1-{r})\log(1-{r})\right)\last{+\frac{c_2}{2}r(1-r)},
  		 \quad & \text{if } {r}\in[0,1],
  		\\		
  		+\infty,\quad  &\text{otherwise},
  \end{cases}
  \label{Flog}
\end{equation}
 with $0<c_1<{c_2}$ adimensional parameters {and, due to its singularities,
it} enforces the tumor concentration to take values in the physical range $[0,1]$. In applications, \eqref{Flog} is {often} substituted by its double-well smooth {polynomial} approximation
\begin{equation}
   F_{\rm reg}({r}):=\frac{c_3}{4} {r}^2 ({r}-1)^2, \quad  {c_3 >0}, \,{r}\in\erre\an{.}
  \label{Freg} 
\end{equation}
{Another possibility is to consider singular and nonregular potentials like}
the {{\it double obstacle potential}}
\begin{equation}
  F_{\rm dob}({r})
  :=\begin{cases}
  c_3 {r} (1-{r}), \quad&\text{if } {r}\in[0,1],\\
  +\infty, \quad&\text{otherwise},
  \end{cases}
  \label{F2ob}
\end{equation}
with $c_3>0$ an adimensional parameter.
We observe that the smooth potential \eqref{Freg} does not enforce the tumor concentration to take values in the physical range $[0,1]$. In biological applications, when cells interaction is predominant with respect to the adhesion between the cells and surrounding tissues, a single-well cellular potential of {\it Lennard--Jones} type is used, which expresses infinite repulsion when the cells are infinitely tight together ({i.e.,} in the situation $\phi\equiv 1$) and attraction when they are far apart ({i.e.,} for small values of $\phi$). The attraction must go to zero with no cells, with the potential having an unstable critical point at $\phi\equiv 0$. In \cite{ALL}, the following phenomenological form for the single-well potential is used
\begin{align}
	\label{LJ}
	 F_{\rm sw}({r}):=
  \begin{cases}
  		- (1-{r}^*)\log(1-{r}) - \frac {{r}^3}3 -(1-{r}^*)\frac {{r}^2}2 
        - (1-{r}^*){r} + \kappa,
  		\, & \text{if } {r}\in[0,1),
  		\\		
  		+\infty ,
  		\,
  		&\text{otherwise},
  \end{cases}
\end{align}
where $\kappa\geq 0$ and ${r}^* \in (0,1)$ corresponds to the volume fraction at which the cells would be at mutual equilibrium. Note that \eqref{LJ} enforces the tumor concentration to take values in the physical range $[0,1)$.\\
We now give general constitutive assumptions such that the equations \eqref{eqn:mix1a}--\eqref{eqn:mix2} satisfy the second law of thermodynamics in isothermal situations and with source terms, which takes the form of the following dissipation inequality {(see \cite{GLSS,G})}
\begin{equation}
\label{eqn:mix4}
\frac{d}{dt}\int_{R(t)}e\leq -\int_{\partial R(t)}{\mathbf{J}_E}\cdot {\nn} 
+\int_{R(t)}\biggl(\frac{\Gamma_{\phi}}{\gamma}m_{\phi}+\frac{\Gamma_a}{\gamma}m_a+{\Sn}m_n+{\Sc}m_c{\biggr)},
\end{equation}
for each material volume $R(t)\subset \Omega$, where ${\nn}$ is the outer normal to $\partial R(t)$, with the energy flux ${\mathbf{J}_E}$ and the multipliers $m_{\phi}, m_a, m_n, $ {and} $m_c$ to be determined. Following similar calculations as those reported in \cite[Section 2.1]{ALL}, we obtain
\begin{align}
\label{eqn:mix5}
\bar{p}&=p+\mu\phi+\mu_a\phi_a+\eta n+{{\badtheta}} c,\\
\label{eqn:mix6}
\vv&=-k\left(\nabla \bar{p}-\left(\mu \nabla \phi+\mu_a \nabla \phi_a+\eta \nabla n + {{\badtheta}} \nabla c\right)\right),\\
\label{eqn:mix7}
\mathbf{J_{\phi}}& =-b(\phi,\phi_a,n,c)\nabla \mu,\\
\label{eqn:mix8}
{\mathbf{J}_a}& =-b_a(\phi,\phi_a,n,c)\nabla \mu_a,\\
\label{eqn:mix9}
F_n& =\alpha_n\eta,\\
\label{eqn:10g}
F_c& =\alpha_c{{\badtheta}},\\
\non
{\mathbf{J}_E}& =\mu\mathbf{J}_{\phi}+\mu_a{{\mathbf{J}_a}}+{{\badeps}}^2
{(\dt \phi)}
\nabla \phi+{(\dt n)}
{D_n} \nabla n+{(\dt c)}
{D_c} \nabla c\\
& \label{eqn:mix10}
 \quad 
+\left(\mu\phi+\mu_a\phi_a+\eta n+{{\badtheta}} c+p-e\right)\vv,\\
m_{\phi}& =\mu, \quad  m_a=\mu_a, \quad m_n=\eta, \quad m_c={{\badtheta}},
\end{align}
where $p$ is the scalar Lagrange multiplier of the constraint \eqref{eqn:mix2}, $k$ is a positive friction parameter, with units of $[mm^2/(Pa\,s)]$, $\alpha_n,\alpha_c$ are positive coefficients related to the time scales of the dynamics of the chemical species, in units of $[\text{Mol}^2/ Pa \, s]$, $b, b_a$ are positive mobilities, and
\begin{alignat}{2}
\label{eqn:mix11}
    & \mu:=\frac{\delta E}{\delta \phi}={\Lambda}\left(F'(\phi)-{{\badeps}}^2 \Delta \phi-\chi_{\phi}n\right), 
    \quad  
    && 
    \mu_a:=\frac{\delta E}{\delta \phi_a}={\Lambda}\left({\kappa_a}\log(\phi_a)-\chi_{a}c\right),
\\
\label{eqn:mix12}
    & \eta:=\frac{\delta E}{\delta n}={\Lambda}\left(-D_n\Delta n-\chi_{\phi}\phi\right), \quad  && {{\badtheta}} :=\frac{\delta E}{\delta c}={\Lambda}\left(-D_c \Delta c-\chi_a\phi_a\right).
\end{alignat}
Inserting \eqref{eqn:mix5}--\eqref{eqn:mix10} in \eqref{eqn:mix1a}--\eqref{eqn:mix2}, we get the following system of equations
\begin{alignat}{2}
\label{eqn:mix13:1}
& \vv=-k\big(\nabla \bar{p}-\left(\mu \nabla \phi+\mu_a \nabla \phi_a+\eta \nabla n + {{\badtheta}} \nabla c\right)\big)
&& \qquad \text{in $\Omega$},
\\  \label{eqn:mix13:2}
& \div\vv=0,
&& \qquad \text{in $\Omega$},
\\ \label{eqn:mix13:3}
 & {\dt \phi}
+\vv\cdot \nabla \phi - \div\big(b(\phi,\phi_a,n,c)\nabla \mu\big)=\frac{\Gamma_{\phi}}{\gamma} ,
&& \qquad \text{in $\Omega$},
\\ \label{eqn:mix13:4}
& {\dt \pha}
+\vv\cdot \nabla \phi_a - \div\big(b_a(\phi,\phi_a,n,c)\nabla \mu_a\big)=\frac{\Gamma_a}{\gamma}
&& \qquad \text{in $\Omega$},
\\ \label{eqn:mix13:5}
& \mu={\Lambda}\left(F'(\phi)-{{\badeps}}^2 \Delta \phi-\chi_{\phi} n\right)
&& \qquad \text{in $\Omega$},
\\ \label{eqn:mix13:6}
& \mu_a={\Lambda}\left({\kappa_a}\log(\phi_a)-\chi_{a}c\right),
\\ \label{eqn:mix13:7}
& {\dt n} +\vv\cdot \nabla n - \alpha_n{D}_n\Delta n-\alpha_n\chi_{\phi} \phi= {\Sn}
&& \qquad \text{in $\Omega$},
\\ \label{eqn:mix13:8}
& {\dt c}+\vv\cdot \nabla c - \alpha_c{D}_c\Delta c-\alpha_c\chi_a \phi_a= {\Sc}
&& \qquad \text{in $\Omega$},
\end{alignat} 
\Accorpa\sysmix {eqn:mix13:1} {eqn:mix13:8}
{which} we endow with the homogeneous boundary conditions
{
\begin{align}
\label{eqn:mix14}
 {b \,\dn \mu=b_a \,\dn \mu_a= \dn \phi=\dn n
 = \dn c =\vv=0
 \qquad \text{on {$\partial \Omega$}}},
\end{align}}%
which imply that
\[
    \mathbf{J}_E\cdot {\nn}=0
    \qquad \text{on {$\partial \Omega$}},
\]%
{and with initial conditions
\begin{equation}
    \label{eqn:mix14ic}
    \phT(0) = \ph^0,\,\,
    \pha(0) = \ph_a^0, \,\,
    n(0) = n^0, \,	\,
    c(0) = c^0 \quad
    \text{in}\; \Omega.        
 \end{equation}
 }%
{Let us point out \cite{KS2}, where a related multiphase system with velocity field subject to Darcy and Brinkmann laws is analyzed.}
 
A solution of system \sysmix, supplemented with the boundary conditions \eqref{eqn:mix14}, formally satisfies the following energy equality
\begin{align}
\label{eqn:mix15}
&\frac{dE}{dt}+\frac{1}{k}\int_{\Omega}|{\vv}|^2  +\int_{\Omega}\left(b|\nabla \mu|^2+b_a|\nabla \mu_a|^2+\alpha_n\eta^2+\alpha_c{{\badtheta}}^2\right)\\
& \notag \quad  =\int_{\Omega}\left(\frac{\Gamma_{\phi}}{\gamma}\mu+\frac{\Gamma_a}{\gamma}\mu_a+{\Sn}\eta+{\Sc}{{\badtheta}}\right).
\end{align}
{We then need to complement the system \sysmix\ with particular} forms for the mobility functions $b,b_a$ and to assign biologically meaningful forms {to}  the source terms $\frac{\Gamma_{\phi}}{\gamma},\frac{\Gamma_{a}}{\gamma},{\Sn},$ {and} ${\Sc}$. The former task, following \cite{ALL}, is accomplished by applying the Onsager Variational Principle (OVP) \cite{O}, which defines the irreversible non-equilibrium dynamics for near-equilibrium systems in terms of linear fluxes-forces balance equations. In isothermal situations, the OVP takes the following form: given a set of slow state variables $\mathbf{x}_i, i=1, \dots, n$, the dynamics of the system is described by the thermodynamic fluxes which minimize the Onsager functional $\mathcal{O}(\dot{\mathbf{x}}_i)=\Phi(\dot{\mathbf{x}}_i)+\dot{E}(\mathbf{x}_i,\dot{\mathbf{x}}_i)$, where $\Phi$ is the dissipation functional, which is quadratic in $\dot{\mathbf{x}}_i$ as a near-equilibrium approximation, and $E$ {stands for} the free energy of the system. We thus minimize \eqref{eqn:mix15} with respect to the variables $\vv_{\phi}, \vv_a, \vv_l$, given the following quadratic approximations for the viscous dissipative terms in \eqref{eqn:mix15}:
\begin{align}
\label{eqn:mix16}
\int_{\Omega}b(\phi,\phi_a,n,c)|\nabla \mu|^2 & {=} \int_{\Omega}M_{\phi}k(\phi,\phi_l,n)|\vv_{\phi}-\vv_l|^2 ,\\
\non
\int_{\Omega}b_a(\phi,\phi_a,n,c)|\nabla \mu_a|^2 & {=} {\int_{\Omega}M_{al}k_{al}(\phi_a,\phi_l,c)|\vv_a-\vv_l|^2}
   \\ & \label{eqn:mix17} \quad  {+   \iO M_{av}k_{av}(\phi_a,\phi,c)|\vv_a-\vv_{\phi}|^2,}
\end{align}
where $M_{\phi}, M_{al},$ {and} $ M_{av}$ are positive friction parameters, related to the friction between the tumor cells and the liquid phase and between the endothelial cells and both the liquid and the tumor phases respectively, and $k(\phi,\phi_l,n)$, $k_{al}(\phi_a,\phi_l,c)$, {and} $k_{av}(\phi_a,\phi,c)$ are generic friction functions, whose form depend on the nature of the filtration processes driven by the drag between the mixture phases, to be empirically determined. The dependence of $k,k_{al},$ {and} $k_{av}$ on their arguments will be described later {on}. We note that here the drag laws \eqref{eqn:mix16} and \eqref{eqn:mix17} are more general than the ones introduced in \cite[Section 2.1]{ALL}, where {it was assumed} $k(\phi,\phi_l,n){=}\phi$ and $k_{al}(\phi_a,\phi_l,c)=k_{av}(\phi_a,\phi_l,c){=} \phi_a$. Substituting \eqref{eqn:mix16} and \eqref{eqn:mix17} in \eqref{eqn:mix15}, with similar calculations as those reported in \cite[Section 2.1]{ALL}, we find that the first order conditions with respect to variations in the variables $\vv_{\phi}, \vv_a,$ {and} $\vv_l$ take the form of the following Darcy {type} laws:
\begin{align}
\label{eqn:mix18}
&\vv_{\phi}-\vv_l=-\frac{(1-\phi)\phi}{M_{\phi}k(\phi,\phi_l,n)(1+\phi_a)}\nabla \mu,\\
\label{eqn:mix19}
&\vv_a-\vv_l=-\frac{\phi_a}{M_{al}k_{al}(\phi_a,\phi_l,c)}\nabla \mu_a, \quad \vv_a-\vv_{\phi}=-\frac{\phi_a}{M_{av}k_{av}(\phi_a,\phi,c)}\nabla \mu_a.
\end{align}
Inserting {those} in \eqref{eqn:mix16} and \eqref{eqn:mix17} {produces}
\begin{align}
\label{eqn:mix20}
&b(\phi,\phi_a,n)=\frac{\phi^2(1-\phi)^2}{M_{\phi}k(\phi,\phi_l,n)(1+\phi_a)^2},\\
\label{eqn:mix21}
&b_a(\phi_a,n,c)=\frac{\phi_a^2}{M_{al}k_{al}(\phi_a,\phi_l,c)+M_{av}k_{av}(\phi_a,\phi,c)}.
\end{align}
Considering \eqref{eqn:mix16} as the viscous dissipation due to a Darcy flow of the liquid phase through the porous-permeable solid matrix associated to the soft material of the tumor phase, a general expression for the friction function $k$ can be given as
\[
k(\phi,\phi_l,n)=\frac{\phi_l\nu_l}{\rho(\phi,\phi_l,n)},
\]
where $\nu_l$ is the viscosity of the liquid phase and $\rho(\phi,\phi_l,n)$ the intrinsic permeability of the tumor phase, assumed to depend on the tumor, the liquid and the nutrient concentrations. A possible expression for $\rho$ can be derived by assuming that the tumor tissue consists of homogeneous and isotropic parallel cylindrical pores, and {the} Poiseuille formula for a capillary tube \cite{HM} {yields} that
\[
\rho(\phi,\phi_l,n)=\xi(n)\frac{r(\phi,\phi_l)^2\phi_l}{8\delta^2},
\]
where $\xi(n)$ is an empirical positive and finite geometrical parameter, whose value may depend on the nutrient availability, $r$ is the effective radius of the pores, depending on $\phi$ and $\phi_l$, and $\delta$ is the tortuosity factor. Since the tumor is a soft tissue, its permeability should also depend on the strain level in the material \cite{HM}, which is neglected in the current modeling framework. A general expression for $r$ is of the form {(see, e.g., \cite{Civ})}
\[
r(\phi,\phi_l)=C\left(\frac{\phi_l}{\phi}\right)^\lambda,
\]
where $C$ is a positive parameter related to the specific internal surface area of the pores and $\lambda$ is a positive empirical parameter. With the latter relations, \eqref{eqn:mix20} becomes
\begin{equation}
\label{eqn:mix22}
b(\phi,\phi_a,n)=\frac{B_{\phi}\xi(n)\phi^{2-2\lambda}(1-\phi)^{2}(1-\phi-\phi_a)^{2\lambda}}{(1+\phi_a)^2},
\end{equation}
where $B_{\phi}$ is a positive parameter related to friction and geometrical coefficients. 
\begin{remark}
\label{rem:saturation}
 We observe that the degeneracy of \eqref{eqn:mix22} for $\phi+\phi_a=1$ enforces the condition $\phi+\phi_a\leq 1$ in the dynamics described by \sysmix. This, together with the conditions $0\leq \phi < 1$, enforced by the cellular potential \eqref{LJ}, and $\phi_a\geq 0$, enforced by the particular form of $\mu_a$ in \eqref{eqn:mix11}, {allow} us {to} interpret the solutions $\phi$ and $\phi_a$ of \sysmix\ as concentrations, implying also the validity of the saturation condition for the underlying mixture model. We also observe that, in the case $\lambda=1$, which corresponds to the well-known Kozeny--Carman law for the intrinsic permeability \cite{K,Car}, the mobility \eqref{eqn:mix22} does not degenerate at $\phi=0$.    
\end{remark}
{To} derive general expressions for the friction functions $k_{al}$ and $k_{av}$, we start by the relations
\[
k_{al}(\phi_a,\phi_l,c)=\frac{\phi_l\nu_l}{\rho_{al}(\phi_a,\phi_l,c)}, \quad k_{av}(\phi_a,\phi,c)=\frac{\phi\nu_{\phi}}{\rho_{av}(\phi_a,\phi,c)},
\]
where $\nu_l,\nu_{\phi}$ are the viscosity of the liquid and the tumor phase respectively and $\rho_{al},\rho_{av}$ are the intrinsic permeability of the endothelial cells phase with respect to the liquid and tumor cells filtration processes respectively. Assuming that the network of tumor induced vasculature made by endothelial cells is described by a random fractal of dimension two with no axis of symmetry embedded in the three dimensional space, the permeability of the endothelial cells network takes the form \cite{Cos}
\[
\rho_{al}(\phi_a,\phi_l,c)=\mathcal{A}(\phi_a,c)\frac{\phi_l}{\phi_a}, \quad \rho_{av}(\phi_a,\phi,c)=\mathcal{A}(\phi_a,c)\frac{\phi}{\phi_a},
\]
where $\mathcal{A}(\phi_a,c)$ is a positive and finite parameter related to geometrical quantities and to the pore cross-sectional area, which may generally depend both on $c$ and $\phi_a$. Introducing the function
\[
\mobn(\phi_a,c):=\frac{\mathcal{A}(\phi_a,c)}{M_{al}\nu_l+M_{av}\nu_{\phi}},
\]
we can write
\begin{equation}
\label{eqn:mix23}
b_a(\phi_a,c)=\phi_a\mobn(\phi_a,c),
\end{equation}
with
\begin{equation}
\label{eqn:mix24}
0<{m_0}\leq \mobn(\phi_a,c)<{M},
\end{equation}
for given positive real numbers ${m_0}$ {and} ${M}$. 
We finally assign biologically meaningful forms for the source terms $\frac{\Gamma_{\phi}}{\gamma},\frac{\Gamma_{a}}{\gamma},{\Sn},$ {and} ${\Sc}$. Following \cite{ALL,GLSS}, we assume that the tumor cells proliferate, with a rate $\nu$, proportionally to the nutrient concentration, as long as the nutrient concentration is above the hypoxia threshold $\delta_n$. Moreover, they die by apoptosis at a rate $\rate$. Hence, we write
\begin{equation}
\label{eqn:mix25}
\frac{\Gamma_{\phi}}{\gamma}=\nu (n-\delta_n)_+h(\phi)-\rate\phi,
\end{equation}
where $h:\mathbb{R}\to [0,1]$ is a continuous function which interpolates linearly between $h(0)=0$ and $h(1)=1$, and is extended as constant outside of the interval $[0,1]$. The source term for the endothelial cells is expressed as
\begin{equation}
\label{eqn:mix26}
\frac{\Gamma_{a}}{\gamma}=\big((c-\delta_a)_+(1-h(\phi))+\zeta\big)({\kappa_0}\phi_a-{\kappa_\infty}\phi_a^2),
\end{equation}
with $\zeta>0$. This means that new vessels form by accumulation of endothelial cells from the existing vasculature following a logistic growth, describing the growth of a population of self-interacting particles with saturation \cite{Ric}. This process is driven by random detachment of endothelial cells from their basement membrane, at a (small) rate $\zeta {\kappa_0}$; outside of the tumor mass, it is driven by the angiogenetic signal, when the concentration of the angiogenetic factor is greater than a proliferation threshold $\delta_a$. The nutrient supply is described by the law 
\begin{equation}
\label{eqn:mix27}
{\Sn} {= {\Sn} (\phi,\pha,n)}=R_n(\bar{n}-n)(1-\phi)+R_a(\bar{n}-n)\phi_a-C_n\phi n,
\end{equation}
where $\bar{n}$ is the typical nutrient concentration inside the capillaries. Nutrients are released from the normal vasculature at a rate $R_n$, as long as $n<\bar{n}$, with the normal capillaries being destroyed as the tumor cells proliferate, and consumed at a rate $C_n$. Moreover, nutrients are supplied by the tumor induced vasculature proportionally to $\phi_a$. Finally, for what concerns the source term of the angiogenetic factor, it is released by the tumor cells at a rate $R_c$ when the nutrient concentration
is below the hypoxia threshold and the angiogenetic factor concentration is below its saturation level $\bar{c}$, and it is consumed by endotelial cells at a rate $C_c$. Hence, {we  have}
\begin{equation}
\label{eqn:mix28}
{\Sc}= {{\Sc} (\phi,\pha,n,c)} = R_ch(\phi)(\delta_n-n)_+(\bar{c}-c)-C_c\phi_ac.
\end{equation}
We substitute now \eqref{eqn:mix22}, \eqref{eqn:mix23} and \eqref{eqn:mix25}--\eqref{eqn:mix28} in \sysmix. We consider the limit of high viscosity of the mixture, which corresponds to $k\to 0$ in \an{\eqref{eqn:mix13:1}}, (which is appropriate for the description of the tumor dynamics, see{, e.g.,} \cite{AZAMM}). {Furthermore}, we rearrange the terms by introducing a new chemical potential
\begin{equation}
\label{eqn:mix29}
\hat{\mu}={-{{\badeps}}^2 \Delta \phi + F^{\prime}(\phi)},
\end{equation}
inserting the chemotactic term as a chemotaxis flux in \eqref{eqn:mix13:3}. In order to derive an adimensionalized version of \sysmix, we also introduce the functions
\[
\hat{\mobm}:=\frac{b}{B_{\phi}}, 
\quad 
\hat{\mobn}:=(M_{al}\nu_l+M_{av}\nu_{\phi})\mobn, 
\quad 
\hat{\mu}_a:=\frac{\mu_a}{{\Lambda}}, 
\quad 
\hat{n}:=\frac{n}{\bar{n}}, 
\quad 
\hat{c}:=\frac{c}{\bar{c}},
\]
the nutrient penetration length
\[
l_n:=\sqrt{\frac{{\Lambda} \alpha_nD_n}{C_n}},
\]
the parameters $\hat{\delta}_n:=\frac{\delta_n}{\bar{n}}$, $\hat{\delta}_a:=\frac{\delta_a}{\bar{c}}$, ${\hat{m}}=\frac{\rate}{\nu}$, ${\hat{\kappa}}_0=\frac{{\kappa_0}}{\nu}$, {and} ${\hat{\kappa}}_{\infty}=\frac{{\kappa_\infty}}{\nu}$ and the change to adimensional space and time variables $\hat{{x}}=\frac{x}{l_n}$, {and} $\hat{t}=t\nu$. Then, system \sysmix\ becomes, without reporting the hat superscripts for ease of notation,
\begin{alignat*}{2}
& {\dt \phi}- \frac{{\Lambda} B_{\phi}}{\nu l_n^2}\div\bigl(\mobm(\phi,\phi_a,n)\nabla \mu\bigr)+\frac{{\Lambda} \chi_{\phi}B_{\phi}\bar{n}}{\nu l_n^2}\div\bigl(\mobm(\phi,\phi_a,n)\nabla n\bigr) && 
\\ &   \quad 
=(n-\delta_n)_+h(\phi)-m\phi
&& \quad \text{in $Q$},\\ 
& \mu={-\frac{{{\badeps}}^2}{l_n^2} \Delta \phi + F'(\phi)}
&& \quad \text{in $Q$},\\ 
& {\dt \pha} - \frac{{\Lambda}}{(M_{al}\nu_l+M_{av}\nu_{\phi})\nu l_n^2}\div \big(\phi_a{\mobn}(\phi_a,c)\nabla \left({\kappa_a}\log(\phi_a)-\chi_{a}\bar{c}c\right)\big) \\
    &  \quad =\left((c-\delta_a)_+(1-h(\phi))+\zeta\right)({\kappa_0}\phi_a-{\kappa_\infty}\phi_a^2)
    && \quad \text{in $Q$},\\ 
 &\frac{\nu}{C_n}{\dt n}-\Delta n-\frac{\alpha_n{\Lambda}\chi_{\phi}}{\bar{n}C_n} \phi= \frac{R_n}{C_n}(1-n)(1-\phi)+\frac{R_a}{C_n}(1-n)\phi_a-\phi n
 && \quad \text{in $Q$},\\ 
 & \frac{\nu}{C_c} {\dt c}- \frac{\alpha_c{D}_cC_n}{\alpha_n{D}_nC_c}\Delta c-\frac{\alpha_c{\Lambda}\chi_a}{\bar{c}C_c} \phi_a= \frac{R_c}{C_c}h(\phi)(\delta_n-n)_+(1-c)-\phi_a c
 && \quad \text{in $Q$}.
\end{alignat*}
Given the values of the optimized parameters reported in \cite[Section 4]{ALL}, we observe that the following adimensional combination of parameters, which will play a role in the analysis developed in the forthcoming sections, take values of the order of magnitude
\begin{equation}
\label{eqn:mix31}
\frac{\alpha_n{\Lambda}\chi_{\phi}}{\bar{n}C_n}\sim 0.01 <1, \quad \frac{\alpha_c{\Lambda}\chi_a}{\bar{c}C_c}\sim 0.001<1. 
\end{equation}
All other adimensional combination of parameters in \eqref{eqn:mix31} do not play a {significant} role in the analysis, so we can take them, without loss of generality and for ease of notation, as equal to one by choosing $\alpha_n=\alpha_c=1$, $D_c=D_n$, $C_n=C_c=\nu=R_n=R_a=R_c$, $\bar{n}=\bar{c}=1$, $l_n={{\badeps}}=1$, ${\Lambda}=\nu$, $B_{\phi}=\frac{1}{M_{al}\nu_l+M_{av}\nu_{\phi}}=1$, {and} ${\kappa_a}=1$. 
Therefore, we obtain
\begin{alignat}{2}
 & \non
 {\dt \phi}- \div\big(\mobm(\phi,\phi_a,n)\nabla  \mu\big)+\chi_{\phi}\div\big(\mobm(\phi,\phi_a,n)\nabla n\big)
 &&
 \\ & \quad  \label{eqn:mix32:1}
 =(n-\delta_n)_+h(\phi)-m\phi
  && \quad \text{in $Q$},
 \\ 
 &  \label{eqn:mix32:2}
 \mu={-\Delta \phi + F'(\phi)}
  && \quad \text{in $Q$},\\ 
& \non {\dt \pha} - \div\big(\phi_a{\mobn}(\phi_a,c)\nabla \left(\log(\phi_a)-\chi_{a}c\right)\big)
 && 
\\ & \quad \label{eqn:mix32:3}
=\left((c-\delta_a)_+(1-h(\phi))+\zeta\right)({\kappa_0}\phi_a-{\kappa_\infty}\phi_a^2)
 && \quad \text{in $Q$},\\ 
& \label{eqn:mix32:4}
{\dt n} -\Delta n-\chi_{\phi} \phi= (1-n)(1-\phi)+(1-n)\phi_a-\phi n
 && \quad \text{in $Q$},\\ 
\label{eqn:mix32:5}
& {\dt c} - \Delta c-\chi_a \phi_a= h(\phi)(\delta_n-n)_+(1-c)-\phi_a c
&& \quad \text{in $Q$},
\\ \non
	{\dn \phT} & {= (\mobm(\phT,\pha,n)\nabla \mu) \cdot \nn= 	(\pha \mobn(\ph_a,c)\nabla ({\log (\pha)} - \chi_a c)) \cdot \nn}
	\\  \quad 	
	\label{eqn:mix32:6}
	& {= \dn n =  \dn c =0} 	 && \quad {\text{on $\Sigma$,}}
	\\
	\label{eqn:mix32:7}
	{\phT(0)} & {= \ph^0,
	\quad 
	\pha(0) = \ph_a^0,
	\quad
	n(0) = n^0,
	\quad 
	c(0) = c^0}
	&& \quad  {\text{in $\Omega$,}}
\end{alignat}
\Accorpa\eqnmix {eqn:mix32:1} {eqn:mix32:7}
where we still maintain the properties \eqref{eqn:mix31}, {i.e., we {require that}} 
\begin{equation}
\label{eqn:mix33}
0<\chi_{\phi} <1, \quad 0<\chi_a<1. 
\end{equation}

\section{Notation, Assumptions and Main Results}
\label{SEC:NOT:ASS:MR}
\setcounter{equation}{0}

To begin with, we assume the set $\Omega$ to be a bounded, connected and smooth open subset of~$\erre^d$, $d \in \{2,3\}$,
with boundary $\Gamma:=\partial\Omega$. 
Given a final time $T>0$, we set, for every $t\in(0,T]$,
\begin{align*}
	Q_t:=\Omega\times(0,t),
	\quad 
	\Sigma_t:=\Gamma\times(0,t),
	\quad 
	Q:=Q_T,
	\quad 
	\Sigma:=\Sigma_T.
\end{align*}
Let $X$ denote a Banach space. {W}e indicate by 
$\norma\cpto_X$, $X^*$, and $\<{\cdot},{\cdot}>_X$
its norm, its dual space, and the associated duality pairing in the order.
As for the classical Lebesgue and Sobolev spaces on $\Omega$, for $1\leq p\leq\infty$ and $k\geq 0$ we use $L^p(\Omega)$ and $W^{k,p}(\Omega)$,
with the standard convention $\Hx k:=\Wx{k,2}$ and norms $\norma{\cdot}_{\Lx p}:= \norma{\cdot}_p$, $\norma{\cdot}_{W^{k,p}(\Omega)}$, and $\norma{\cdot}_{H^{k}(\Omega)}$. Similar symbols are employed to denote spaces and norms constructed on $Q$, $\Gamma$ and $\Sigma$.
For convenience, we set
\begin{align*}
  && H:=L^2(\Omega), \quad
  V:=H^1(\Omega), \quad
  W:=\graffe{ v\in H^2(\Omega): \partial_{\nn}v=0\;\,\text{a.e.~on } \Gamma},
\end{align*}
and endow them with their corresponding norms $\norma{\cdot}:= \norma{\cdot}_H, \norma{\cdot}_V,$ and $\norma{\cdot}_W$, respectively.
As usual, $H$ will be identified to its dual so that 
we have the following continuous, dense, and compact embeddings:
\begin{align*}
	W \emb V \emb H \emb V^*
\end{align*}
along with the identification 
\begin{align*}
	\<u,v >_V = \iO uv, \quad 
	u \in H, v\in V.
\end{align*}
Finally, for every $v  \in \Vp$, we employ $(v)_\Omega := \frac 1{|\Omega|} \<v, 1>_V$ to indicate the generalized mean value of $v$. Sometimes, when no confusion may arise, we simply use $v_\Omega$ instead of $(v)_\Omega$.
We {then use} $V_0$, $H_0$, and $V_0^*$ {to denote} the closed subspaces of functions with zero spatial
mean of $V$, $H$, and $V^*$, respectively. Then, the operator $-\Delta$
with homogeneous Neumann boundary conditions {may be} considered as
\begin{align*}
   (-\Delta) : V \to V^*, 
   \qquad 
   \<(-\Delta) v, z>_V:= \iO \nabla v\cdot \nabla z,
   \quad 
   \text{$v,z\in V$\an{.}}
\end{align*}
{It follows that it}
 is invertible when {restricted} to act on functions
with zero spatial {average}.
Namely, $-\Delta : V_0 \to V_0^*$ is invertible and we denote its inverse by $\calN: =(-\Delta )^{-1} : V_0^* \to V_0$. 
It is well-known that 
\begin{align*}
    \norma{v^*}_*:= \norma{\nabla ({\cal N} v^*)} = (\nabla ({\cal N} v^*),\nabla ({\cal N} v^*))^{1/2} = \< v^*,{\cal N} v^*>_V^{1/2},
     \quad v^* \in V_0^*,
\end{align*}
yields a Hilbert norm on $V_0^*$.
In addition, it holds that
\begin{align*}
\<{-\Delta v},{\calN v^*}>_V&= \<{v^*},{v}>_V, \quad 
\<{v^*},{\calN w^*}>_V= (v^*, w^* )_{*}\andre{,}
\quad 
\text{{$v \in V_0,$ \, $v^*,w^* \in V_0^*$,}}
\end{align*}
where the symbol $(\cdot, \cdot )_{*}$ denotes the
{standard} inner product of $V^*$.
Furthermore, if $v^* \in \H1 {V^*_0}$, we have{, for a.e~$t \in (0,T)$,} that
\begin{align*}
\<{\partial_t v^*(t)},{\calN v^*(t)}>_V= \frac 12 \frac d {dt} \norma{v^*(t)}^2_{*}\,.
\end{align*}
{
Finally, let us introduce the notation  $(\cdot)_\pm$  for the positive and negative part function, respectively. Namely, $(\cdot)_\pm:\erre \to [0,+ \infty)$ are defined as
\begin{align*}
    ({r})_+ := \max \{{r}, 0\},
    \quad
    \text{and}
    \quad
    ({r})_- := - \min \{{r}, 0\}{, \quad r \in \erre.}
\end{align*}
}%

{
Let us also mention here the following {standard result} that will be useful later on.
\begin{lemma}
    \label{LEM:dist}
    Let $f,g\in L^1(0,T)$, $g_0\in \mathbb{R}$, and, for any $\varrho \in C_c^{\infty}([0,T))$, let
    \[
    -\int_0^T\varrho'(\tau)g(\tau)\,d\tau+\int_0^T\varrho(\tau)f(\tau)\,d\tau-\varrho(0)g_0=0.
    \]
    Then, it holds that
    \[
    g(t)-g(s)+\int_s^tf(\tau)d\tau=0,
    \]
    for $a. e. \,\,t,s\in [0,T)$, including $s=0$, provided we replace $g(0)$ with $g_0$.
\end{lemma}
The proof of Lemma \ref{LEM:dist} is a consequence of the fundamental lemma of the calculus of variations \cite[Lemma 1.2.1]{JJ}, see also \cite[Lemma 3.1]{Lsz} for a similar result.
}%

To conclude, let set a useful convention for the appearing constants. From now on,  the capital $C$ 
will be used to denote a generic constant
whose actual values may change from line to line and even within the same line
and depends only on structural data of the system. When specific constants enter the computations, like $\delta$ for instance,  we will employ {self-explanatory}
subscript{s} like $C_\delta$ to indicate that the constant depend\an{s} on the parameter $\delta$, in addition.

Let us now 
 \abramo{make some preliminary remarks on the qualitative properties satisfied by System \eqref{eqn:mix32:1}--\eqref{eqn:mix32:5}, with the constraint \eqref{eqn:mix33}, which will be justified throughout the forthcoming calculations. In particular, we observe that:
\begin{itemize}
\item In the case with a smooth potential like \eqref{Freg}, the variable $\phv$ is not {guaranteed} to satisfy the {pointwise} property $0\leq \phv \leq 1$. As a consequence, a solution to \eqref{eqn:mix32:4} does not formally satisfy the maximum and minimum principles $0\leq n\leq 1$;
\item A solution of \eqref{eqn:mix32:3} formally satisfies the minimum principle $\pha\geq 0$ {due to the logarithmic term} {in the free energy}; 
\item A solution to \eqref{eqn:mix32:5} formally satisfies the {minimum and maximum} principles $0\leq c\leq 1$. As we will see in the following, this property is fundamental to ensure the coercivity of the chemotaxis term $-\chi_a\iO  \pha c $ arising in the free energy \an{$ E$}. 
\end{itemize} 
\begin{remark}
\label{rem:plog}
In \cite{RSchS}\an{,} a generalized logistic growth for \eqref{eqn:mix32:3} of the form
\begin{equation*}
\frac{\Gamma_{a}}{\gamma}\propto {\kappa_0}\phi_a-{\kappa_\infty}\phi_a^p
{\quad \text{ with $p \in (1,2]$}}
\end{equation*}
was considered for the Keller--Segel system. Considering the latter growth law in our model, a source term for \eqref{eqn:mix32:5} of the form
\begin{equation*}
{\Sc}= {{\Sc} (\phi,\pha,n,c)} = R_ch(\phi)(\delta_n-n)_+(\bar{c}-c)-C_c\phi_a^{p-1}c,
\end{equation*}
should be considered. In this situation, the property $c \leq 1$ would be valid also for $1<p<2$ only if {we can ensure that} $0\leq \pha\leq 1${. This property is not always trivial to obtain, but can be reached, for instance, assuming in System \eqnmix\ a degenerate mobility $\mobm(\phv,\pha,n)$ as in \eqref{eqn:mix22} (see Remark \ref{rem:saturation}).} With a non-degenerate mobility $\mobm(\phv,\pha,n)$ and $p<2$, in order to ensure the coercivity of the chemotaxis term $-\chi_a\iO  \pha c $, we would need to introduce the property that $c\leq 1$ as a constraint in equation \eqref{eqn:mix32:5}, e.g., by adding to the free energy \eqref{eqn:mix3} the indicator function of the set $c\leq 1$.
\end{remark}
In light of the above properties, we will make different structural assumptions corresponding to the cases with a smooth or a singular potential.}
{Before diving into listing the mathematical assumptions on the system, let us point out that 
the following structural assumptions on the source terms are motivated by the aforementioned discussion. 
From a mathematical perspective, instead of specifying a particular form for these terms, it suffices to {postulate} specific growth conditions.
However, it is worth noticing  that this approach is only applicable in certain cases (cf. \eqref{ass:Source:1}--\eqref{ass:Source:4}), and these conditions may vary when dealing with regular and singular potentials. To simplify the technical aspects as much as possible, we opt to adopt a specific structure {for the sources} that maintains a high degree of generality and facilitates the analysis.}
\begin{enumerate}[label={\bf A\arabic{*}}, ref={\bf A\arabic{*}}]
\item \label{ass:1:potential}
\abramo{
In the case of a {\it smooth potential}, we postulate that {$F$ is defined on the whole real line,} $F\in C^{2}(\mathbb{R})$, and there exists $c_1\geq 0$ such that
\begin{equation}
\label{ass:potsm:1}
|F'(r)|\leq c_1\left(F(r)+1\right), \; F(r)\geq 0, \;  r \in \mathbb{R}.
\end{equation}
{Besides, $F$ enjoys the decomposition} $F=\Beta+\Pi$, where $\Beta$ is convex and $\Pi$ is concave, and {there exists} $c_2\geq 0$ such that 
\begin{equation}
\label{ass:potsm:2}
|\Pi''(r)|\leq c_2\left(|r|^q+1\right), \;  r \in \mathbb{R},
{\quad \text{with $q\in [0,4)$.}}
\end{equation}
We observe that the smooth potential \eqref{Freg} satisfies Assumption \eqref{ass:potsm:1}.\\
\noindent
In the case of a {\it singular potential}, we postulate $F$ to be decomposed in a {singular}{, proper, and}  convex part $\Beta$ and a nonconvex, smooth, perturbation $\Pi$ with a quadratic {growth}. 
Namely, we require that 
\begin{align}
	&F:\erre\to\errebar \; \hbox{enjoys the splitting} \; F=\Beta +\Pi, \; \hbox{where}
	\qquad
	\label{ass:pot:1}
	\\
	& \Beta :\erre\to [0,+\infty] \; \text{is {proper,} convex and l.s.c.\
	  with subdifferential} \;
	  \beta := \partial\Beta,
	\non
	\\ \label{ass:pot:2}
	& \hbox{and fulfills} \; 
    \beta(0)\ni 0,
	\; 
	\text{$\b$ is $C^2$ in the interior of its domain $D(\b)$,}
\end{align}
whereas
\begin{align}
        \label{ass:pot:3}
	& {\pi : \erre \to \erre,\quad }
 \pi \in C^1(\erre), \quad \pi:=\Pi ' \,\, \text{ is \Lip\ continuous, {and}} 
     \\ & 	\label{ass:pot:4} 
    |\Pi(r)|\leq c_1(|r|^2+1), 
    \quad 
    |\pi(r)|\leq c_2(|r|+1),
    \quad 
    |\pi'(r)|\leq c_3, {\quad  r \in \mathbb{R},}
\end{align}
for some \last{nonnegative constants $c_1,c_2,$ and $c_3$.}
\begin{remark}
\label{rem:swtrunc}
We note that both the double well potential \eqref{Flog} and the {single-well} potential \eqref{LJ} satisfy Assumptions \eqref{ass:pot:1}--\eqref{ass:pot:3}, while Assumption \eqref{ass:pot:4} is satisfied by \eqref{Flog} only. Hence, we will treat the {single-well} potential \eqref{LJ} by adopting a proper truncation procedure {to} let its truncated {form} satisfy Assumption \eqref{ass:pot:4}.
In particular, the {single-well} potential \eqref{LJ} can be decomposed as $F_{{\rm sw}}=\Beta+\Pi$, where
{
\begin{align*}
    & \Beta({r})= \begin{cases}
  		-(1-{r}^*)\log(1-{r})
  		\, & \text{if } {r}\in[0,1),
  		\\		
  		+\infty  
  		\,
  		&\text{otherwise},
  \end{cases}
    \\
    & \Pi({r})=-\frac{{r}^3}{3}-(1-{r}^*)\frac{{r}^2}{2}-(1-{r}^*){r}+\kappa,
    \quad r \in \erre,
\end{align*}
}%
{with} $k\geq 0$ and ${r}^*\in (0,1)$. {For this latter, the growth of the corresponding perturbation $\Pi$ is of third order instead of second order. Thus, we} will consider a quadratic truncation of $\Pi$ which preserves its regularity and concavity, defined as
\begin{equation*}
\ov{\pi}(r):=
\begin{cases}
\Pi(0)+r\pi(0)+\frac{r^2}{2}\pi'(0)  &\text{for $r\leq 0$},\\
\Pi(r)  &\text{for $0<r< 1$},\\
\Pi(1)+(r-1)\pi(1)+\frac{(r-1)^2}{2}\pi'(1)  &\text{for $r\geq 1$}.
\end{cases}
\end{equation*}
Hence, we will consider a truncated form $\ov{F}_{\rm{sw}}=\Beta+\ov{\pi}$, which satisfies \eqref{ass:pot:2}--\eqref{ass:pot:4}. Since we will prove that $0\leq \ph <1$, actually ${F}_{\rm{sw}}(\ph)\equiv \ov{F}_{\rm{sw}}(\ph)$.
\end{remark}
It is well-known that $\b$ yields a maximal monotone graph in $\erre \times \erre$ with corresponding domain
$D(\b):=\{{r}\in\erre:\b({r})<+\infty\}$.}

{
\item \label{ass:h}
We suppose the interpolation function $h:\mathbb{R}\to [0,1]$ to be continuous and  such that 
\begin{align*}
    h(r) = 
    \begin{cases}
    0,
    \quad \text{$r \leq 0$},
    \\
    r ,   \quad \text{$r \in (0,1)$},
    \\
    1  ,   \quad \text{$r \geq 1$}.
    \end{cases}
\end{align*}
}

\item \label{ass:2:theta}
\abramo{
	We assume that $\Sa$ possesses a logistic growth of the form
{\begin{align} 
   \Sa =  \Sa (\ph, \pha,c)& =  
   {\theta(\ph,c)({\kappa_0}\phi_a-{\kappa_\infty}\phi_a^2),}  
   \label{def:Sasm}
\end{align}
where 
\begin{align*}
    \theta(\ph,c):=\left((c-\delta_a)_+(1-h(\phi))+\zeta\right),
    \quad \delta_a \in [0,1], \, \zeta,\kappa_0,\, \kappa_\infty > 0.
\end{align*}}%
We observe that  ${\theta=\theta(\ph,c)}$ is strictly positive and, due to {\ref{ass:h} and} the property that $0\leq c\leq 1$, it is also uniformly bounded. Moreover, we assume that
\begin{align}
	\label{def:Sc}
	\Sc & =\Sc(\ph,\pha,c)= h(\phi)(\delta_n-n)_+(1-c)-\pha c.
\end{align}%
}%
\item \label{ass:sources:3}
\abramo{%
For what concerns the source terms for the variables $\ph$ and $n$, {since} the property $0\leq n\leq 1$ is valid only in the case with a singular potential, we will need to make different assumptions {discerning the cases of smooth and  singular potentials.} In particular, the source terms \eqref{eqn:mix25} and \eqref{eqn:mix27} will be properly truncated in the case with a smooth potential.}%

{
For these reason, we postulate
\begin{align}
    \label{def:calS}
    \Sv &  
  = \Sv (\ph,n)= {\cal H}(\ph,n)-m\phi, 
  \quad 
  m >0,
\end{align}
where, for $\delta_n\in [0,1],$
\begin{align}
    \label{def:calH}
{\cal H} & ={\cal H}(\ph,n):=
 \begin{cases}
 (h(n)-\delta_n)_+h(\phi)
 & \text{when the potential is  smooth},
 \\
 (n-\delta_n)_+h(\phi)
 & \text{when the potential is   singular}.
 \end{cases}
 \end{align}
}%
{Besides, when the potential is singular, we require} the compatibility condition
\begin{align*}
	& 
	- \frac {H}  m - (\ph^0)_\Omega^-,
	\quad 
	\quad 
	 \frac {H} m + (\ph^0)_\Omega^+
	 \quad 
	\text{belong to the interior of $D(\beta),$}
\end{align*}
where $H:={\norma{\cal H}_\infty=}{\norma{(n-\delta_n)_+h(\phi)}_{\infty}}$. 
{Furthermore}, in the case {of}  a smooth potential, we assume that
\begin{align}
    \label{def:Snsm}
	\Sn & = \Sn(\ph,\pha,n)= (1-h(n))(1-h(\ph)+\pha)-\phi h(n),
\end{align}
{while, in} the case with a singular potential, we assume that
\begin{align}
     \label{def:Sn}
	\Sn & = \Sn(\ph,\pha,n)= (1-n)(1-h(\ph)+\pha)-\ph n.
\end{align}%

\item \label{ass:mobilities:4}
We assume $\mobm\in C^0(\erre^3)$ and 
{ $\mobn\in C^0(\erre^2)$}
to be globally Lipschitz continuous and
there exist positive constants $m_0$ and $M$ such that
\begin{alignat}{2} \label{mob:bou}
  & 0 < m_0 \leq \mobm(\ph,\ph_a,n),\mobn(\ph_a,c) \leq M < + \infty{,}
  &&  \quad  \ph, n{,c}\in\erre,~\ph_a \geq 0{.}
\end{alignat}

\item \label{ass:parameters:5}
For the chemotaxis sensitivities $\chi_\ph$ and $\chi_a$\andre{,} we require that 
\abramo{
\begin{align*} 
    \chi_\ph, \chi_a : 
    \begin{cases}
        \chi_\ph \geq 0 , \hspace{10.7mm} \chi_a \in( 0,1)  \quad 
        & \text{with smooth potential},
        \\
         \chi_\ph \in (0,1), \quad \chi_a  \in( 0,1)       
         & \text{with singular potential}.
    \end{cases}
\end{align*}
}%
\end{enumerate}

The first result we are going to present concerns the existence of weak solutions in both  two and three \abramo{space} dimensions.
\begin{theorem}[Existence of weak solutions, $d \in \{2,3\}$]
\label{THM:WEAK:LOC}
Suppose that \ref{ass:1:potential}--\ref{ass:parameters:5} hold. 
Moreover, let the initial data fulfill
\begin{align}\label{ass:ini:weak:1}
	& \ph^0 \in V,
	\quad 
	F(\ph^0) \in \Lx1,
	\quad
	(\ph^0)_\Omega \in D(\beta), 
	\\ \label{ass:ini:weak:2}
	& 
	\ph_a^0 \geq 0 \;\; \text{$a. e.$ in $\Omega$},
	\quad 
	\ph_a^0 {\log (\ph_a^0)} \in \Lx1,
	\\ & \label{ass:ini:weak:3}
	n^0 \in V,
	\quad 
	c^0 \in V\cap \Lx\infty,
	\quad 
	0 \leq c^0 \leq 1 \;\; \text{$a. e.$ in $\Omega$}.
\end{align}
Besides, let the threshold {Sobolev exponent}
{
\begin{align}
    \label{def:esp}
    \esp :=\;\text{arbitrary in $(1, + \infty)$ if $d=2$, and $\esp := 6$ if $d=3$.}
\end{align}
}%
Then, the multiphase Cahn--Hilliard--Keller--Segel model \Sys\ admits at least a weak solution.
Namely, there exists a quintuple $(\ph, \pha, \mu , n, c)$  such that
\begin{align}
	& \label{reg:weak:1}
	\phT \in \H1 \Vp \cap \L\infty V \cap \L4 {W} \cap \L2 {\Wx{2,\esp}},
 	\\ &  \label{reg:weak:4}
 	\pha(x,t) \geq 0 \quad \text{for $a.e. \, (x,t) \in Q$,}
 	\\ &  \label{reg:weak:5}
 	{\pha \in \C0 {{(\Wx{1,4})^*}}  \cap \L{\frac {d+2}{d+1}} {\Wx{1,\frac {d+2}{d+1}}},}
 	\\ &  \label{reg:weak:6}
 	\pha \an{\log (\pha)} \in \L\infty {\Lx1},
 	\quad 
 	\an{\ph_a^2} \an{\log (\pha)} \in \L1 {\Lx1},
 	\\ &  \label{reg:weak:8}
 	\mu \in \L2 V,
 	\\ &  \label{reg:weak:10}
 	n \in \H1 H \cap \L\infty V \cap \L2 W,
 	\\ &  \label{reg:weak:11}
  	c \in L^\infty(Q): \quad 0 \leq c(x,t) \leq 1 \quad \text{for $a.e. \, (x,t) \in Q$,}	
 	\\ &  \label{reg:weak:12}
 	c \in \H1 H \cap \L\infty V {\cap \L2 W}.
\end{align}
Besides, it fulfills the pointwise formulation
\[
\mu = -\Delta \phT + F'(\phT) \quad \text{$a.e.$ in $Q$},
\]
where, in the case of a {\it singular potential}, $F'(\ph)=\xi + \pi (\phT)$, with
\abramob{
\begin{equation}
\label{reg:weak:2}
	\xi \in \L2 {\Lx \esp},
 \end{equation} 
 and} $ \xi \in \b(\phT) $ $a.e.$ in $Q$,
along with the weak formulations, recall the definition of the source terms in  \an{\eqref{def:Sasm}--\eqref{def:Sn}},
\begin{align}
	& \label{wf:1}
	\<\dt \phv, v >_V
	+ \iO \mobm(\ph, \ph_a,n)\nabla \mu \cdot \nabla v
	- \chi_\ph \iO \mobm(\ph, \pha, n) \nabla n \cdot \nabla v 
	=
	\iO \Sv v,	
  \\
 \non
    & {
     \iO \pha(t) w
	+{\intQt} \mobn(\phi_{a},c)\nabla \phi_{a}\cdot \nabla w
	- \chi_a \intQt \ph_a \mobn(\ph_a,c) \nabla c \cdot \nabla w}
	\\ & \quad 
 \label{wf:3}
	{=
	 \iO \ph_a^0  w
	 + 
 \intQt \Sa w,
 }
	\\ \label{wf:4}
	& \iO \dt n\, v 
	+ \iO \nabla n \cdot \nabla v
	- \chi_\ph \iO \phv v
	=
	\iO \Sn v,
 \\ \label{wf:5}
	& \iO \dt c \,v 
	+ \iO \nabla c\cdot \nabla v
	- \chi_a \iO \pha v
	=
	\iO \Sc v,
\end{align}
for {almost every $t \in (0,T)$,} every $v \in V$ and 
{$w \in {\Wn},$ 
where ${\Wn}$ is defined as
\begin{align*}
	{\Wn}:= \Wx{1,d+2}.
\end{align*}
}%
Moreover, the initial conditions in \eqref{system:8} are fulfilled in the sense that
\begin{align}
	\label{wf:initial:1}
	& \ph(0) = \ph^0,
	\quad 
	n(0) = n^0,
	\quad 
	c(0) = c^0
	\quad \text{$a.e.$ in $\Omega$,}
	\\
	\label{wf:initial:2}
	& 
	\pha(0) = \ph_a^0
	\quad \text{in ${(\Wx{1,4})^*}$}. 
\end{align}

\end{theorem}

\begin{remark}\label{RMK:INDATA}
We notice that in case of the regular potential \eqref{Freg}, the second condition in \eqref{ass:ini:weak:1} is already fulfilled.
In fact, we have $F({r}) = {\cal O}({r}^4)$ as $|{r}| \to +\infty$  as well as $V \emb \Lx {\esp}$ with $\esp$ as defined in \eqref{def:esp}.

Besides, in case of singular potentials like \eqref{Flog}--\eqref{F2ob}, we have $D(\beta)=[0,1]$ and the second condition in \eqref{ass:ini:weak:1} imposes the initial datum $\ph^0$ 
to be uniformly bounded and such that $0 \leq \ph^0(x) \leq 1$ for almost every $x \in \Omega$.
Moreover, the singularity of the potential yields that the order parameter $\phT$ belongs to the  physical range, that is,
\begin{align}
	 	\non
	 	{\phT \in L^\infty(Q): \quad 0 \leq \phT(x,t) \leq 1 \quad \text{for $a.e. \, (x,t) \in Q$.}}
\end{align}
{In that case, it also holds that
\begin{align*}
  	n \in L^\infty(Q): \quad 0 \leq n(x,t) \leq 1 \quad \text{for $a.e. \, (x,t) \in Q$.}
\end{align*}}
\end{remark}

\begin{remark}\label{RMK:CHO}
The prescribed {form of $\SS$ in \eqref{def:calS}} is the typical choice one encounters in the Cahn--Hilliard--Oono equation, where the function {$\cal H$} reduces to a constant ${{\cal H}} \in (-m,m)$.
We are aware of the recent contribution \cite{GSS}, where the authors show that the last average condition in \eqref{ass:ini:weak:1} can be actually relaxed a little bit allowing pure phases to be considered as initial data. 
Due to the complexity of our system, the result does not directly apply and we left that problem open for possible future research.
\end{remark}

{
\begin{remark}
    Let us highlight some differences to the work referenced as \cite{RSchS}. The first one lies in the choice of the source $\Sa$, as observed in the Remark \ref{rem:plog}.
     In contrast, our different formulation of the chemotaxis coupling has allowed us to deduce somewhat improved regularities for $\pha$ (cf. \eqref{reg:weak:5}).
     This improvement is further manifested in a more  consistent variational framework (cf. \eqref{wf:3}): it is worthwhile to compare the two setting  of test functions
     related to the chemotactic variable.
    \end{remark}
}


Assuming 
{a more regular initial datum $\ph_a^0$ \abramob{and a smallnes condition for the chemotactic \last{sensitivity}  $\chi_a$ with respect to the magnitude of the initial datum and other parameters of the model},} we can derive some regularity results in the two space dimensions.
\begin{theorem}[Regularity result, $ d =2$]
\label{THM:REG:LOC}
Suppose that \ref{ass:1:potential}--\ref{ass:parameters:5} hold and let $d =2$. 
Moreover, in addition to \eqref{ass:ini:weak:1}, suppose {that}
\begin{align}
    \label{ini:reg:strong}
	\ph_a^0 \in H.
\end{align}
\abramob{Suppose further that the following smallness condition on the chemotactic parameter,
\begin{equation}
    \label{smallness}
    \chi_a<\left(\frac{\sqrt{1+\ov{C}}-1}{2}\right)^{\frac{1}{4}},
\end{equation}
is satisfied, where $\ov{C}$ is a positive parameter depending only on the domain $\Omega$, on the parameters $\kappa_{\infty}$ and $m_0$ and on proper norms of the initial conditions \last{(cf. \eqref{def:cbar})}.}
\newline
Then, the components $\pha$ and $c$ of the weak solution obtained from Theorem~\ref{THM:WEAK:LOC} enjoy the {additional} regularities
\begin{align}\label{reg:reg:p2}
	\pha & \in \H1 \Vp \cap \L\infty  H \cap \L2 V,
	\\ \label{reg:reg:c:p2}
	c & \in \L4 {{W}},
\end{align}
and the weak formulation \eqref{wf:3} can be {equivalently reformulated as}
\begin{align*}
	\non
	& \<\dt \pha, v>_V
	+ \iO \mobn(\ph_a,c)  \nabla \pha \cdot \nabla v
	- \chi_a \iO \pha \mobn(\ph_a,c) \nabla c \cdot \nabla v
	=  {\iO \Sa v},
 \end{align*}
{almost everywhere in $(0,T)$ and for every $v \in V$}.
Besides the initial condition $\pha(0)=\ph_a^0$ is fulfilled almost everywhere in $\Omega$.

\end{theorem}

Next, provided that the initial data is more regular, the mobility functions are constant, {and} the space dimension is {two,} we can show that there exist more regular weak solutions. Here, the first regularity result  follows.

\begin{theorem}[Regularity result on $n$ and $c$, $d \in \{2,3\}$]
\label{THM:REG:STRONG:PREL}
Suppose that \ref{ass:1:potential}--\ref{ass:parameters:5} \abramob{and \eqref{smallness}} are fulfilled,
{and assume that, besides to \ref{ass:h}, it holds that $h \in W^{1,\infty}(\erre)$.}
Moreover, 
in addition to \eqref{ass:ini:weak:1} and \eqref{ini:reg:strong}, suppose that 
\begin{align}
    \label{ini:reg:strong:prel}
    \dt n (0)& := \Delta n^0 + \chi_\ph \ph^0 + \Sn (\ph^0,\ph_a^0,n^0,c^0) \in H,
    \quad 
    c^0 \in \Hx2.
\end{align}
Then, there exist components $n$ and $c$ of a weak solution $( \phd, \pha, \xi, \mu , n, c)$ such that 
\begin{align}
   \label{reg:prel:1}
   n & \in \W{1,\infty} H \cap \H1 V \cap \L\infty {\Hx2} \cap \L2 {\Hx3},
    \\ \label{reg:prel:2}
    c & \in \H1 H \cap \L\infty {\Hx2} \cap \L2 {\Hx3}.
\end{align}
\end{theorem}

\begin{theorem}[Regularity result, $ d =2$, constant mobilities]
\label{THM:REG:STRONG}
Suppose that \ref{ass:1:potential}--\ref{ass:parameters:5} \abramob{and \eqref{smallness}} hold and let $ d =2$ and that ${\mobm\equiv \mobn} \equiv 1$.
{Besides, in addition to \ref{ass:h}, we suppose that $h \in W^{1,\infty}(\erre)$.}
Moreover, 
in addition to \eqref{ass:ini:weak:1}, \eqref{ini:reg:strong}, and \eqref{ini:reg:strong:prel}, suppose that 
\begin{align}
    \label{ini:reg:strong:bis}
    \ph^0 & \in W,
    \quad 
    \mu^0:= -\Delta \ph^0 + F'(\ph^0) \in V,
    \quad 
    \ph_a^0 \in V.
\end{align}
Then, there exists a weak solution $(\ph,\pha,\xi,\mu,n,c)$ such that 
\begin{align}\label{reg:strong:1}
	\ph & \in \W{1,\infty} \Vp \cap \H1  V \cap \L\infty {\Wx{2,\esp}},
	\\ \label{reg:strong:2}
	\xi & \in \L\infty {\Lx \esp},
    \\ \label{reg:strong:3}
	\mu & \in \L\infty {V},
    \\ \label{reg:strong:4}
	\pha & \in \H1 H \cap \C0 V \cap \L2 {\Hx2},
\end{align}
with $\esp$ being defined as in \eqref{def:esp}.

\end{theorem}

Finally, \abramob{in the case with a singular potential,} if the convex part of the double-well potential enjoys the following estimate
\begin{align}
    \label{growth}
    \exists \, c_\beta >0:
    \quad 
    |\beta ' ({r}) | \leq e^{c_\beta (|\beta ({r})| +1 )},
    \quad 
    {r} \in D(\beta),
\end{align}
we can obtain another regularity improvement.
The above condition is known to {be fulfilled} in the two dimensional setting, for instance, by the singular logarithmic potential \eqref{Flog}.
Notice that a similar version holds for the single-well potential \eqref{LJ} when its argument is close to one. Namely, {for \eqref{LJ}, condition} \eqref{growth} holds for every ${r} \in ( 1/2,1)$ instead.

\begin{theorem}[Regularity result, $ d =2$, separation property]
\label{THM:REG:SEP}
{Suppose the assumptions of the Theorem \ref{THM:REG:STRONG} are fulfilled and the double-well potential {enjoys} \eqref{growth}.}
Moreover,
in addition to \eqref{ass:ini:weak:1}, \eqref{ini:reg:strong} and \eqref{ini:reg:strong:bis}, suppose that 
\begin{align}
    \label{ini:reg:sep}
    \ph^0 \in H^4(\Omega),
    \quad 
    \mu^0 \in W.
\end{align}
Then, there exists a weak solution $(\ph,\pha,\xi,\mu,n,c)$ such that 
\begin{align}\label{reg:sep:1}
	\ph & \in \W{1,\infty} H \cap \H1  {\Hx2} \cap \L\infty {\Hx4 \cap \Wx{2,\esp}},
    \\ \label{reg:sep:2}
	\mu & \in \L\infty {\Hx2} \cap \L2 {\Hx3},
\end{align}
with $\esp$ as defined in \eqref{def:esp}.
Besides, if \eqref{growth} is fulfilled,
then there exist ${\d}_*, {\d}^* \in (0,1)${, ${\d}_*\leq {\d}^*$,}  such  that the separation property  holds:
\begin{align}
    \label{separation}
    0 < {\d}_* \leq \ph(x,t) \leq {\d}^* <1 \quad \text{{for every} $(x,t ) \in {\ov Q}.$}
\end{align}
Finally, in case the potential is single-well and fulfills \eqref{growth} for every ${r} \in ( 1/2,1)$, then there exist ${\d}^*$  such  that the separation property holds:
\begin{align}
    \label{separation:LJ}
    0 \leq \ph(x,t) \leq {\d}^* <1 \quad 
    \text{{for every} $(x,t ) \in {\ov Q}.$}
\end{align}

\end{theorem}

The last result we are going to address concerns the uniqueness of solutions. This is obtained as consequence of a suitable continuous dependence estimate that is fulfilled by regular solutions that enjoy the regularities listed in the above theorems.

\begin{theorem}[Uniqueness,  $ d =2$, constant mobilities]
\label{THM:UNIQUENESS}
{Suppose the assumptions of Theorem \ref{THM:REG:SEP} are fulfilled.}
{Besides, the source term $\Sa$ possesses the simplified form, compare with \eqref{def:Sasm},
$ \Sa=\Sa(\pha)={\kappa_0}\phi_a-{\kappa_\infty}\phi_a^2.$
}
Then there exists a unique weak solution
$(\ph,\pha, \mu,n,c)$ to the system \Sys.
Moreover, let $\{(\ph_i, \ph_a, \mu_i,n_i, c_i)\}_i$, $i=1,2$, denote a couple of weak solutions as obtained from Theorem \ref{THM:REG:SEP} associated to initial data $\{(\ph^{0}_i, \ph_{a,i}^{0},\mu^{0}_i,n^{0}_i,c^{0}_i)\}_i$
fulfilling, for $i=1,2$,
\eqref{ass:ini:weak:1}, \eqref{ini:reg:strong}, \eqref{ini:reg:strong:bis}, and \eqref{ini:reg:sep}.
Then, it holds that 
\begin{align}
    \non
    & \norma{\ph_1-\ph_2 - ((\ph_1)_\Omega-(\ph_2)_\Omega )}_{\L\infty \Vp}
    + \norma{(\ph_1)_\Omega-(\ph_2)_\Omega}_{L^\infty(0,T)}
    \\ & \qquad \non
    + \norma{\ph_{a,1}-\ph_{a,2} - ((\ph_{a,1})_\Omega-(\ph_{a,2})_\Omega)} _{\L\infty \Vp \cap \L2 H}
    + \norma{(\ph_{a,1})_\Omega-(\ph_{a,2})_\Omega}_{L^\infty(0,T)}
    \\ & \qquad  \non
    + \norma{n_1-n_2}_{\L\infty H \cap \L2 V}
    + \norma{c_1-c_2}_{\L\infty H \cap \L2 V}
     \\ & \quad  \non
   \leq K
    \big( 
    \norma{\ph_1^0-\ph_2^2 - ((\ph_1^0)_\Omega-(\ph_2^0)_\Omega)}_\Vp
    +   |(\ph_1^0)_\Omega-(\ph_2^0)_\Omega|
    \big)
      \\ & \qquad  \non
  + K \big(\norma{\ph_{a,1}^{0}-\ph_{a,2}^{0}- ((\ph_{a,1}^{0})_\Omega-(\ph_{a,2}^{0})}_\Vp
  + |(\ph_{a,1}^{0})_\Omega-(\ph_{a,2}^{0})_\Omega|
    \big)
   \\ &     \label{cd:est}
   \qquad 
    + K \big(\norma{n^{0}_1-n^{0}_2} + \norma{c^{0}_1-c^{0}_2}\big),
\end{align}
for a positive constant $K$ only depending on the data of the system.

\end{theorem}

\section{Existence of Weak Solutions}
\setcounter{equation}{0}

\abramo{
In this section {we establish the validity of} Theorem~\ref{THM:WEAK:LOC}. 
{Our approach begins with the introduction of an {approximation} for System} \eqnmix, which allows us proving the existence of a local in time solution through a Faedo--Galerkin {scheme}.
{Subsequently, we expand this local solution into a global-in-time solution using a combination of a-priori estimates, which remain uniform with respect to the discretization parameter, and temporal continuity arguments.} Finally, we pass to the {limit} letting the regularization parameter goes to zero, recovering a solution to the original system. Since the Assumptions \ref{ass:1:potential}, \ref{ass:2:theta} and \ref{ass:sources:3},  and the qualitative properties of solutions of System \eqnmix, differ between the cases of $F$ being smooth or singular, i.e., satisfying \eqref{ass:potsm:1}--\eqref{ass:potsm:2} or \eqref{ass:pot:1}--\eqref{ass:pot:3}, respectively, {we adopt two distinct regularization approaches depending on the assumed regularity of the potential $F$.}}

{In this direction, let us introduce some preliminary tools.}
\abramo{
{To begin with, given} $L,M\in \mathbb{R}$ with $L <M$, we define the truncation function $T_{L,M}$ as
\begin{equation}
    \label{Tlm}
    T_{L,M}:\mathbb{R}\rightarrow [L,M],
    \quad 
    T_{L,M}({r}):=\max \{L,\min({r},M)\}, \quad  {r}\in \mathbb{R},
\end{equation}
and notice that $T_{L,M}\in W^{1,\infty}(\mathbb{R})$.
Following standard regularizing approaches in the Keller--Segel literature (see, e.g., \cite{HH}),
we define the function $\mathcal{E}_{L,M}:\mathbb{R}\rightarrow \mathbb{R}^+$, with $\mathcal{E}_{L,M}\in C^2(\mathbb{R})$, such that
\begin{equation}
    \label{Elm}
    (\mathcal{E}_{L,M}^{\prime\prime}T_{L,M})({r})=1, \quad  {r}\in \mathbb{R}{.}
\end{equation}
{Namely, if $L<1<M$ we impose} $\mathcal{E}_{L,M}^{\prime}(1)=\mathcal{E}_{L,M}(1)=0$, {and find that}
\begin{align*}
    \non
    \mathcal{E}_{L,M}^{\prime\prime}({r})& :=
    \begin{cases}
        \frac{1}{L}, \qquad {r}\leq L,\\
        \frac{1}{{r}}, \qquad \andre{L<{r}<M},\\
        \frac{1}{M}, \quad \hspace{3mm} {r}\geq M,\\      
    \end{cases}
    \mathcal{E}_{L,M}^{\prime}({r}):=
    \begin{cases}
        \frac{{r}}{L}+\andre{\log}(L)-1, \quad\hspace{5mm} {r}\leq L,\\
        \andre{\log}({r}), \qquad\qquad\qquad \andre{L<{r}<M},\\
        \frac{{r}}{M}+\andre{\log}(M)-1, \quad \hspace{2.6mm} {r}\geq M,
    \end{cases}
    \\ 
    \mathcal{E}_{L,M}({r})& :=
    \begin{cases}
        \frac{{r}^2-L^2}{2L}+(\andre{\log}(L)-1){r}+1 ,\quad\hspace{5mm} {r}\leq L,\\
        (\andre{\log}({r})-1){r}+1, \qquad\qquad\qquad \andre{L<{r}<M},\\
        \frac{{r}^2-M^2}{2M}+(\andre{\log}(M)-1){r}+1, \quad \hspace{2.6mm}{r}\geq M.      
    \end{cases}
\end{align*}
We recall the following properties concerning the functions $T_{L,M}$ and $\mathcal{E}_{L,M}$, which are derived in \cite{HH} and which will be useful in the \andre{forthcoming} calculations:
\begin{alignat}{2}
    \label{propTE1}
    &  \mathcal{E}_{L,M}({r})\geq \tfrac{{r}^2}{2L}, \quad &&  {r}\leq 0,\,  L \in (0,{e}^{-1}),\\
    \label{propTE3}
    &  {r}\mathcal{E}_{L,M}^{\prime}({r})\leq 2\mathcal{E}_{L,M}({r})+1, \quad &&  {r}\in \mathbb{R},\,  L \in (0,{e}^{-1}).
\end{alignat}
Moreover, we give the following properties, which can be directly verified by computation and which will be useful later:
\begin{align}
    \label{sleqe1}
    &  |{r}|\leq \mathcal{E}_{L,L^{-1}}({r})+{e}-1, \quad &&  L \in (0,{e}^{-1}), \,{r\in \mathbb{R},} \\
    \label{sleqe2}
    &  ({r})_+^2\mathcal{E}_{L,L^{-1}}^{\prime}({r})+(2{e})^{-1}\geq 0, \quad &&  L \in (0,{e}^{-1}), \,{r\in \mathbb{R}}, \\
    \non
    &  ({r})_+^2\leq {\ov{C}}\left(({r})_+^2\mathcal{E}_{L,L^{-1}}^{\prime}({r})+(2{e})^{-1}\right)
    \quad && 
    \\ & \label{sleqe3}
    \qquad \quad 
    +\max\left({e}^{-\frac{2(1+{\ov{C}})}{{\ov{C}}}}\left({\ov{C}}+\tfrac{4}{27{\ov{C}}^2}\right),{e}^{\frac{2}{{\ov{C}}}}\right) ,
    && L \in \left(0,{e}^{-\frac{1+{\ov{C}}}{{\ov{C}}}}\right),\, {r\in \mathbb{R}},
\end{align}
{for any positive constant $\ov C$.}
%
}
\subsection{Approximation}
{The primary challenge in introducing an approximate scheme is {associated with the necessity of obtaining} uniform (in the approximation parameters) a-priori energy estimates for the Keller--Segel system, which are tipically derived only at a formal level in the literature (see, e.g., \cite{Szk,RSchS}). In our case the energy estimate for the coupled system \eqnmix\ is related to the formal dissipative equality \eqref{eqn:mix15}, and in order to control the chemotactic coupling term $-\chi_a\iO \pha c$ in the free energy \eqref{eqn:mix3} the boundedness of $c$ is needed. Moreover, some of the source terms in \eqref{eqn:mix15} can be controlled thanks to the boundedness and the nonnegativity of certain variables. Hence, another challenge in the design of our approximate scheme is related to the need to maintain the physical boundedness (min-max conditions) and nonnegativity characteristics of the aforementioned variables also at the approximation level. This is achieved starting from the introduction of a proper regularization of the entropy density associated to the variable $\pha$, following \cite{HH}. We observe that in \cite{HH} the chemotactic coupling between the cell density and the chemical concentration was complemented by the introduction of an artificial cell diffusion in the chemical concentration equation; this avoided the need to control the chemotactic coupling term in the a-priori estimates.
The desing of our approximation schemes relates on three main ingredients:
\begin{itemize}
    \item The introduction of a truncation of the variable $\pha$ in the chemotactic flux in equation \eqref{eqn:mix32:3};
    \item The introduction of proper truncations and positive parts of some variables involved in the source terms;
    \item The introduction of a proper regularization of the {possibly} singular potential in equation \eqref{eqn:mix32:2}. 
\end{itemize}
}%
In the case with a smooth potential $F$ satisfying \eqref{ass:potsm:1}--\eqref{ass:potsm:2}, we introduce the following regularized and truncated version of System \eqnmix, depending on the regularization parameter $\an{\eps \in(0,1)}$:
\begin{alignat}{2}
 & \non
 {\dt \phi}- \div\big(\mobm(\phi,\phi_a,n)\nabla  \mu\big)+\chi_{\phi}\div\big(\mobm(\phi,\phi_a,n)\nabla n\big)
 &&
\\ & \quad  \label{eqn:smoothreg:1}
=\left(h(n)-\delta_n\right)_+h(\ph)-m\phi
 && \quad \text{in $Q$,}\\ 
 & \label{eqn:smoothreg:2} \mu={-\Delta \phi + F'(\phi)}
  && \quad \text{in $Q$,}
 \\  
& \notag{\dt \pha} - \div\big(\mobn(\phi_a,c)\nabla \phi_a\big)+\chi_a\div\big(T_{\epsilon,\epsilon^{-1}}(\phi_a)\mobn(\phi_a,c)\nabla c\big) &&\\
& \label{eqn:smoothreg:3} \quad 
=\left((T_{-\epsilon^{-1},1+\epsilon^{-1}}(c)-\delta_a)_+(1-h(\phi))+\zeta\right)({\kappa_0}\phi_a-{\kappa_\infty}(\phi_a)_+^2)
 && \quad \text{in $Q$,}
\\ 
& \label{eqn:smoothreg:4} {\dt n} -\Delta n-\chi_{\phi} \phi= (1-h(n))(1-h(\phi)+(\phi_a)_+)-\phi h(n)
 && \quad \text{in $Q$,}
\\ 
\non
& {\dt c} - \Delta c-\chi_a (\phi_a)_+
&& 
\\ & \label{eqn:smoothreg:5}
\quad 
= h(\phi)(\delta_n-n)_+(1-T_{-\epsilon^{-1},1+\epsilon^{-1}}(c))-(\phi_a)_+T_{-\epsilon^{-1},1+\epsilon^{-1}}(c)
 && \quad \text{in $Q$.}
\end{alignat}
\Accorpa\smoothreg {eqn:smoothreg:1} {eqn:smoothreg:5}
\an{As for the boundary and initial conditions, they will be selected later on for the unified \andre{approximated} system.}
\an{T}o define the regularization of System \eqnmix\ in the case with a singular potential $F$ satisfying \eqref{ass:pot:1}--\eqref{ass:pot:4}, we introduce the Moreau--Yosida regularizations (see, e.g., \cite[pp.~28 and~39]{Brezis}) of the functional $\Beta$ and the graph $\beta$, depending on a parameter $\an{\eps \in(0,1)}$, i.e.\andre{, we set}
\begin{align*}
	&\betaeps:= \frac {I - J_{\eps}}{\eps}, 
	\quad J_{\eps}:=(I+ \eps \an{\beta})^{-1}, \\
	& 0\leq \Betaeps(\an{r}) := \min_{t\in \mathbb{R}}\left\{\frac{1}{2\eps}|t-\an{r}|^2+\Beta(t)\right\},  \quad \an{r} \in \erre,
\end{align*}
being $I$  the identity operator.
From classical theory of convex analysis and \an{A}ssumption \ref{ass:1:potential} we have that, for every $\an{\eps \in (0,1)}$,
\begin{align*}
  & \hbox{$\betaeps$ is monotone and $\frac 1 \eps$-\Lip\ continuous with $\betaeps(0)=0$,}
  \\ 
 &  |\betaeps(\an{r})| \leq |\beta^\circ(\an{r})|
  \quad \hbox{for every $\an{r}\in D(\beta)$},
\end{align*}
where $\beta^\circ(\an{r})$ indicates the element of the section $\beta(\an{r})$ having minimum modulus. Moreover, \an{it readily follows that}
\[
\Betaeps(\an{r}) =\frac{\eps}{2}|\betaeps(\an{r})|^2+\Beta(J_{\eps}(\an{r})), \qquad \an{r} \in \erre,
\]
from which, since $J_{\eps}(\an{r})$ belongs to the proper domain of $\Beta$ and using the Young inequality, we get that
\begin{equation}
\label{eqn:my:1}
|\betaeps(\an{r})|\leq \frac{C}{\eps}\left(\Betaeps(\an{r})+1\right), \qquad \an{r} \in \erre, \, \an{\eps \in (0,1).}
\end{equation}
Then, we introduce the following regularized and truncated version of System \eqnmix, depending on the regularization parameter $\an{\eps \in(0,1)}$:
\an{
\begin{alignat}{2}
 & \non
 {\dt \phi}- \div\big(\mobm(\phi,\phi_a,n)\nabla  \mu\big)+\chi_{\phi}\div\big(\mobm(\phi,\phi_a,n)\nabla n\big)
 &&
 \\ \label{eqn:reg:1} & \quad 
 =\left(T_{-\epsilon^{-1},1+\epsilon^{-1}}(n)-\delta_n\right)_+h(\ph)-m\phi
 && \quad \text{in $Q$,}
 \\ 
 & \label{eqn:reg:2} 
 \mu={-\Delta \phi + \betaeps(\ph)+\pi(\ph)}
  && \quad \text{in $Q$,}
 \\  
& \non
 \dt \pha
- \div\big(\mobn(\phi_a,c)\nabla \phi_a\big)+\chi_a\div\big(T_{\epsilon,\epsilon^{-1}}(\phi_a)\mobn(\phi_a,c)\nabla c\big)
&&
\\
& \label{eqn:reg:3} \quad 
=\left((T_{-\epsilon^{-1},1+\epsilon^{-1}}(c)-\delta_a)_+(1-h(\phi))+\zeta\right)({\kappa_0}\phi_a-{\kappa_\infty}(\phi_a)_+^2) 
&& \quad \text{in $Q$,}
\\ 
& \label{eqn:reg:4}
{\dt n} -\Delta n-\chi_{\phi} (\phi)_+
&& 
\\ & 
\quad 
= \left(1-T_{-\epsilon^{-1},1+\epsilon^{-1}}(n)\right)(1-h(\phi)+(\phi_a)_+)-(\phi)_+T_{-\epsilon^{-1},1+\epsilon^{-1}}(n)
 && \quad \text{in $Q$,}
\\ 
\non
& {\dt c} - \Delta c-\chi_a (\phi_a)_+
&&
\\ & \label{eqn:reg:5}
\quad 
= h(\phi)(\delta_n-n)_+(1-T_{-\epsilon^{-1},1+\epsilon^{-1}}(c))-(\phi_a)_+T_{-\epsilon^{-1},1+\epsilon^{-1}}(c)
 && \quad \text{in $Q$.}
\end{alignat}}%
\Accorpa\reg {eqn:reg:1} {eqn:reg:5}
{
\begin{remark}
    \label{rempos}
    We observe that in the source terms of equations \eqref{eqn:smoothreg:5} and \eqref{eqn:reg:5} (and also of equation \eqref{eqn:reg:4}), we took the positive part of $\pha$, (\andre{$\ph$, respectively}), instead of their truncations
    \andre{ 
    by two primary reasons: firstly, to maintain the integrity of the min-max principles within the equations, and secondly, to enable the reconstruction of the term} $-\chi_a\frac{d}{dt}\iO \pha c$, and the analogous term involving $\ph$ and $n$, in the a-priori estimate. The latter point could not be afforded if we considered the truncation of the variable $\pha$ (\andre{$\ph$, respectively}), instead of their positive part in the chemotactic terms in equations \eqref{eqn:smoothreg:5} and \eqref{eqn:reg:5} \andre{(equation \eqref{eqn:reg:4}, respectively)}.
\end{remark}
}%
For ease of notation, we have omitted to indicate with a subscript $\epsilon$ the dependence of the solutions of \smoothreg\ and \reg\ on the regularization parameter $\epsilon$. \\
\noindent
Given $\an{\eps \in(0,1)}$, we can prove the existence of a solution to \smoothreg\ and \reg\, at least locally in time, e.g., by means of a Faedo--Galerkin approximation strategy. \an{To unify the discussion, we} introduce the functions $q_{\an{\eps}}$,  $p$ and $F_{\eps}$, defined on $\mathbb{R}$, as
\begin{equation*}
q_{\an{\eps}}(r):=
\begin{cases}
h(r), & \text{if $F$ is smooth},\\
T_{-\epsilon^{-1},1+\epsilon^{-1}}(r), &\text{if $F$ is singular},
\end{cases}
\quad
p(r):=
\begin{cases}
r, & \text{if $F$ is smooth},\\
\andre{r_+}, &\text{if $F$ is singular},
\end{cases}
\end{equation*}
and
\begin{equation*}
F_{\eps}(r):=
\begin{cases}
F(r) & \text{if $F$ is smooth},\\
\Betaeps(\ph)+\Pi(\ph) &\text{if $F$ is singular}.
\end{cases}
\end{equation*}
We observe that $\an{q_{\eps}}$ is bounded \andre{for any $\an{\eps \in(0,1)}$} and \an{is} Lipschitz continuous. Moreover, thanks to \eqref{ass:potsm:1}, \eqref{eqn:my:1} and the polynomial growth \eqref{ass:pot:4}, we have that there exists $C>0$, eventually depending on $\eps$ in the case with a singular potential, such that
\begin{equation}
\label{fpf}
|F'_{\eps}(r)|\leq C(F_{\eps}(r)+1), 
\quad F_{\eps}(r)\geq 0, 
\quad r\in \mathbb{R}.
\end{equation}
\an{Let us now fix the boundary and initial conditions. For those, we consider no-flux  Neumann boundary conditions for all the variables and as initial conditions we will consider the $H$-projection of the original initial data as detailed below.}

Let $\{\gtest_i\}_{i\in \mathbb{N}}$ be the \andre{family of} eigenfunctions of the Laplace operator with homogeneous Neumann boundary conditions, that is, \andre{for $i=0,...,k$, $\gtest_i$ are weak solutions to}
\[
-\Delta \gtest_i=\alpha_i \gtest_i \quad \text{in} \; \Omega, \quad \dn  \gtest_i ={0} \quad \text{on} \; \an{\Gamma},
\]
with $0=\alpha_0< \alpha_1 \leq \dots \leq \an{\alpha_\gal <...}\to \infty$ the reordered sequence of eigenvalues. The sequence $\{\gtest_i\}_{i\in \mathbb{N}}$ can be chosen as an orthonormal basis in $\an{H}$ and an orthogonal basis in $\an{V}$, and, thanks to the \an{properties of} $\an{\Gamma}$, \andre{it holds that} $\{\gtest_i\}_{i\in \mathbb{N}}\subset \an{W}$. Then, we introduce the \an{$H$-}projection operator
\[
\an{P_\gal}:\an{V}\to \andre{V_\gal:=}\text{span}\{\gtest_0,\gtest_1,\dots,\gtest_\an{k}\}\an{,} \quad \an{\gal \in \enne,}
\]
\an{and notice that $\bigcup_{k=0}^{\infty}V_k$ is dense in \andre{both} $V$ and in $H$.}
Accordingly, we make the Galerkin ansatz 
\begin{align}
    \non
    \an{\phi_\gal}\an{(x,t)}& =\sum_{\an{i=0}}^{\an{\gal}} a_i^\gal(t){\gtest_i}({x}), 
    \quad 
    \an{\mu_\gal}\an{(x,t)}=\sum_{\an{i=0}}^{\an{\gal}} b_i^\gal(t)\gtest_i({x}),
    \quad
    \pham\an{(x,t)}=\sum_{\an{i=0}}^{\an{\gal}} c_i^\gal(t)\gtest_i({x}),
    \quad 
    \\ 
    \label{Gal:2}
    \an{n_\gal}\an{(x,t)} & =\sum_{\an{i=0}}^{\an{\gal}} d_i^\gal(t)\gtest_i({x}),
    \quad 
    \an{c_\gal}\an{(x,t)}=\sum_{\an{i=0}}^{\an{\gal}} e_i^\gal(t)\gtest_i({x}),
\end{align}
\an{for unknowns functions $a_i,b_i,c_i,d_i,$ and $e_i$, $i=0,...,k$}
to approximate the solutions $\phi, \mu, \phi_a, n,c$ of systems \smoothreg\ and \reg, and project  the equations onto $\an{V_\gal}$, obtaining the following Galerkin approximation of systems \smoothreg\ and and \reg:
\begin{align}
 & \non
 \iO {\dt \phi}_{\an{\gal}}\,\gtest_i+ \iO \mobm(\an{\phi_\gal},\pham,\an{n_\gal})\nabla  \left(\an{\mu_\gal}-\chi_{\phi}\an{n_\gal}\right) \cdot \nabla \gtest_i
 \\ & \quad \label{smoothregdisc:1} 
 =\iO \left((q_{\an{\eps}}(\an{n_\gal})-\delta_n)_+h(\an{\phi_\gal})-m\an{\phi_\gal}\right)\gtest_i,
 \\  \label{smoothregdisc:2}  
 &  \iO\an{\mu_\gal}\,\gtest_i=\iO\nabla \an{\phi_\gal}\cdot \nabla \gtest_i + \iO F_{\eps}'(\an{\phi_\gal})\gtest_i,
 \\ \non
&   \iO {\dt \pham}\,\gtest_i + \iO \mobn(\pham,\an{c_\gal})\nabla \pham\cdot \an{\nabla} \gtest_i
\andre{-}\chi_a\iO T_{\epsilon,\epsilon^{-1}}(\pham)\mobn(\pham,\an{c_\gal})\nabla \an{c_\gal}\cdot \an{\nabla} \gtest_i\\
&  \label{smoothregdisc:3}  \quad =   \iO \left((T_{-\epsilon^{-1},1+\epsilon^{-1}}(\an{c_\gal})-\delta_a)_+(1-h(\an{\phi_\gal}))+\zeta\right)({\kappa_0}\pham-{\kappa_\infty}(\pham)_+^2)\gtest_i,
\\\notag 
&  \notag \iO {\dt \an{n_\gal}}\, \gtest_i +\iO \nabla \an{n_\gal}\cdot \nabla\gtest_i- \chi_{\phi}\iO p(\an{\phi_\gal})\,\gtest_i\\
& \label{smoothregdisc:4}   \quad = \iO \left((1-q_{\an{\eps}}(\an{n_\gal}))(1-h(\an{\phi_\gal})+(\pham)_+)-p(\an{\phi_\gal})q_{\an{\eps}}(\an{n_\gal})\right)\gtest_i,
\\ \notag 
&  \notag \iO {\dt \an{c_\gal}}\,\gtest_i + \iO \nabla \an{c_\gal} \cdot \nabla \gtest_i-\chi_a\iO  (\pham)_+\gtest_i \\
\label{smoothregdisc:5}
& \quad =  \iO \left(h(\an{\phi_\gal})(\delta_n-\an{n_\gal})_+(1-T_{-\epsilon^{-1},1+\epsilon^{-1}}(\an{c_\gal}))-(\pham)_+T_{-\epsilon^{-1},1+\epsilon^{-1}}(\an{c_\gal})\right)\gtest_i,
\end{align}
{in $[0,t]$}, with $0<t\leq T$, for $i=0, \dots,\gal$ and with initial conditions 
\begin{equation}
    \label{discic}
    \an{\phi_\gal}(0)=\an{P_\gal}(\phi^0), \quad \pham(0)=\an{P_\gal}(\phi_{a}^0), \quad \an{n_\gal}(0)=\an{P_\gal}(n^0), \quad \an{c_\gal}(0)=\an{P_\gal}(c^0).
\end{equation}
We note that thanks to the introduction of the functions $q_{\an{\eps}}$ and $F_{\eps}$ we have written a unique \an{G}alerkin approximation \an{\eqref{smoothregdisc:1}--\eqref{discic}} which is valid for both the systems \smoothreg\ and \reg.
\begin{remark}
\label{rem:bc}
The Galerkin ansatz for the variables $\pha$ and $c$ implies that both $\pham$ and $\an{c_\gal}$ satisfy homogeneous Neumann boundary conditions. This is compliant with the boundary conditions \eqref{eqn:mix14} employed in the derivation of the model and with the constitutive assumption \eqref{mob:bou}. Indeed, given the homogeneous Neumann boundary condition for the variable $c$ in \eqref{eqn:mix14} and the assumed positivity of the mobility function $\mobn$, the second boundary condition in \eqref{eqn:mix14} implies that $\dn \pha =0$. Nevertheless, we highlight the fact that in the limit system \eqref{wf:3}, the validity of the boundary condition for the variable $\phi_a$ is lost due to its low space and time regularity.
\end{remark}
\abramo{System \an{\eqref{smoothregdisc:1}--\eqref{discic}} defines a collection of initial value problems for a system of coupled normal first-order ODEs in the variables $a_i(t),\dots,e_i(t)$, $i=0,\dots, \gal$.
Due to the Assumptions \ref{ass:1:potential}, \ref{ass:mobilities:4} on the regularity of the functions $\mobm, \mobn, F,$ \an{and} $h$, the regularity of the Moreau--Yosida approximation and the regularity in space of the eigenfunctions $\gtest_i$, the structure function of the ODEs system depends continuously on the independent variables and on the coefficients. Hence, we can apply the \an{Cauchy--}Peano existence theorem to infer that there exist a sufficiently small $t_1$ with $0<t_1\leq T$ and a \andre{corresponding} local solution $(a_i^\gal,\dots,e_i^\gal)$ of \an{\eqref{smoothregdisc:1}--\eqref{discic}}, for $i=0, \dots, \gal$.
\an{Upon combining \eqref{Gal:2} with the properties of the eigenfunctions, this \andre{readily} produces local solutions $\ph_\gal, \mu_\gal, \ph_{a,\gal}, n_\gal $ and $c_\gal$.}
We now deduce \last{a-priori} estimates, uniform in the discretization parameter $\gal$, for the solutions of system \an{\eqref{smoothregdisc:1}--\eqref{discic}}, which can be rewritten, combining the equations over $i=0, \dots, \gal$, as
\begin{align}
 & \non
 \iO {\dt \phi}_{\an{\gal}}\,v+ \iO \mobm(\an{\phi_\gal},\pham,\an{n_\gal})\nabla  \left(\an{\mu_\gal}-\chi_{\phi}\an{n_\gal}\right) \cdot \nabla v
\\
& \quad \label{smoothregdisc2:1}
=\iO \left((q_{\an{\eps}}(\an{n_\gal})-\delta_n)_+h(\an{\phi_\gal})-m\an{\phi_\gal}\right)v,
 \\ 
 & \label{smoothregdisc2:2}
 \iO\an{\mu_\gal}\,v=\iO\nabla \an{\phi_\gal}\cdot \nabla v + \iO F_{\eps}'(\an{\phi_\gal})v,
 \\ \non
& \notag \iO {\dt \pham}\,v + \iO \mobn(\pham,\an{c_\gal})\nabla \pham\cdot \nabla v
\andre{-}\chi_a\iO T_{\epsilon,\epsilon^{-1}}(\pham)\mobn(\pham,\an{c_\gal})\nabla \an{c_\gal}\cdot \nabla v\\
& \label{smoothregdisc2:3}
    \quad = \iO \left((T_{-\epsilon^{-1},1+\epsilon^{-1}}(\an{c_\gal})-\delta_a)_+(1-h(\an{\phi_\gal}))+\zeta\right)({\kappa_0}\pham-{\kappa_\infty}(\pham)_+^2)v,
\\ \non
&   \iO {\dt \an{n_\gal}}\, v +\iO \nabla \an{n_\gal}\cdot \nabla v-\chi_{\phi} \iO p(\an{\phi_\gal})\,v
\\& \quad  \label{smoothregdisc2:4}
= \iO \left((1-q_{\an{\eps}}(\an{n_\gal}))(1-h(\an{\phi_\gal})+(\pham)_+)-p(\an{\phi_\gal})q_{\an{\eps}}(\an{n_\gal})\right)v,
\\ \notag 
&  \notag \iO {\dt \an{c_\gal}}\,v + \iO \nabla \an{c_\gal} \cdot \nabla v- \chi_a\iO (\pham)_+v \\
\label{smoothregdisc2:5}
& \quad = \iO \left(h(\an{\phi_\gal})(\delta_n-\an{n_\gal})_+(1-T_{-\epsilon^{-1},1+\epsilon^{-1}}(\an{c_\gal}))-(\pham)_+T_{-\epsilon^{-1},1+\epsilon^{-1}}(\an{c_\gal})\right)v,
\end{align}
for \an{almost every} $t \in [0,t_1]$ and for all $v \in \an{V_\gal}$, with initial conditions defined in \eqref{discic}. We take $v= \an{\mu_\gal}-\chi_{\phi}\an{n_\gal}+\an{\phi_\gal}$ in \eqref{smoothregdisc2:1}, $v= -\dt \an{\phi_\gal}$ in \eqref{smoothregdisc2:2}, $v= \pham$ in \eqref{smoothregdisc2:3}, $v= \dt \an{n_\gal}+2(\chi_{\phi}^2+1)\an{n_\gal}$ in \eqref{smoothregdisc2:4} and $v= \dt \an{c_\gal}+\an{c_\gal}$ in \eqref{smoothregdisc2:5}, sum all the equations and \an{manipulate} some terms, to obtain that
\begin{align}
&\notag \frac{d}{dt} \biggl(\frac{1}{2}\norma{\an{\phi_\gal}}^2+\frac{1}{2}\norma{\nabla \an{\phi_\gal}}^2+\iO F_{\eps}(\an{\phi_\gal})-\chi_{\phi}\iO \an{\phi_\gal n_\gal}+\frac{1}{2}\norma{\pham}^2+(\chi_{\phi}^2+1)\norma{\an{n_\gal}}^2\\
&\notag
    \qquad +\frac{1}{2}\norma{\nabla \an{n_\gal}}^2
    +\frac{1}{2}\norma{\an{c_\gal}}^2+\frac{1}{2}\norma{\nabla \an{c_\gal}}^2\biggr)+m\norma{\an{\phi_\gal}}^2+\norma{\dt \an{n_\gal}}^2+2(\chi_{\phi}^2+1)\norma{\nabla \an{n_\gal}}^2\\
&\notag
    \qquad +\norma{\dt \an{c_\gal}}^2+\norma{\nabla \an{c_\gal}}^2 +\iO \mobm(\an{\phi_\gal},\pham,\an{n_\gal})\nabla \an{\mu_\gal}\cdot \nabla \an{\mu_\gal}\\
&\notag
    \qquad+\chi_{\phi}\iO \mobm(\an{\phi_\gal},\pham,\an{n_\gal})\nabla  \an{n_\gal}\cdot \nabla \an{n_\gal}+ \iO \mobn(\pham,\an{c_\gal})\nabla \pham\cdot \nabla \pham\\
&\notag \qquad +
\underbrace{\kappa_\infty \iO \left((T_{-\epsilon^{-1},1+\epsilon^{-1}}(\an{c_\gal})-\delta_a)_+(1-h(\an{\phi_\gal}))+\zeta\right)(\pham)_+^3}_{\geq 0}\\
& \notag \qquad + \underbrace{2(\chi_{\phi}^2+1)\iO \left(1-h(\an{\phi_\gal})+(\pham)_+\right)q_{\an{\eps}}(\an{n_\gal})\an{n_\gal}}_{\geq 0}\\
&\notag \qquad +\underbrace{\iO \left(h(\an{\phi_\gal})(\delta_n-\an{n_\gal})_++(\pham)_+\right)T_{-\epsilon^{-1},1+\epsilon^{-1}}(\an{c_\gal})\an{c_\gal}}_{\geq 0}\\
&\notag \quad = 2\chi_{\phi}\iO \mobm(\an{\phi_\gal},\pham,\an{n_\gal})\nabla  \an{\mu_\gal}\cdot \nabla \an{n_\gal}-\iO \mobm(\an{\phi_\gal},\pham,\an{n_\gal})\nabla \an{\mu_\gal}\cdot \nabla \an{\phi_\gal}\\
&\notag \qquad + \chi_{\phi}\iO \mobm(\an{\phi_\gal},\pham,\an{n_\gal})\nabla \an{n_\gal}\cdot \nabla \an{\phi_\gal} +\iO (q_{\an{\eps}}(\an{n_\gal})-\delta_n)_+h(\an{\phi_\gal})(\an{\mu_\gal}-(\an{\mu_\gal})_{\Omega})\\
&\notag \qquad+\iO (q_{\an{\eps}}(\an{n_\gal})-\delta_n)_+h(\an{\phi_\gal})(\an{\mu_\gal})_{\Omega} - \chi_{\phi}\iO (q_{\an{\eps}}(\an{n_\gal})-\delta_n)_+h(\an{\phi_\gal})\an{n_\gal}\\
&\notag \qquad+\iO (q_{\an{\eps}}(\an{n_\gal})-\delta_n)_+h(\an{\phi_\gal})\an{\phi_\gal}-m\iO \an{\phi_\gal}(\an{\mu_\gal}-(\an{\mu_\gal})_{\Omega})-m|\Omega|(\an{\phi_\gal})_{\Omega}(\an{\mu_\gal})_{\Omega}\\
&\notag \qquad + \chi_{\phi}(m+2(\chi_{\phi}^2+1))\iO \an{\phi_\gal n_\gal}-\chi_a\iO T_{\epsilon,\epsilon^{-1}}(\pham)\mobn(\pham,\an{c_\gal})\nabla \an{c_\gal}\cdot \nabla \pham\\
&\notag \qquad + \kappa_0\iO \left((T_{-\epsilon^{-1},1+\epsilon^{-1}}(\an{c_\gal})-\delta_a)_+(1-h(\an{\phi_\gal}))+\zeta\right)\phi^2_{\an{a,\gal}}\\
& \notag \qquad +\iO (1-q_{\an{\eps}}(\an{n_\gal}))\left(1-h(\an{\phi_\gal})+(\pham)_+\right)\dt \an{n_\gal}
\\
&\notag \qquad
+2(\chi_{\phi}^2+1)\iO \left(1-h(\an{\phi_\gal})+(\pham)_+\right)\an{n_\gal}\\
& \notag \qquad -\iO p(\an{\phi_\gal})q_{\an{\eps}}(\an{n_\gal})(\dt \an{n_\gal}+2(\chi_{\phi}^2+1)\an{n_\gal}) +\chi_{\phi}\iO \left(p(\an{\phi_\gal})-\an{\phi_\gal}\right)(\dt \an{n_\gal}+2(\chi_{\phi}^2+1)\an{n_\gal})\\
&\notag \qquad+\chi_a\iO (\pham)_+\dt \an{c_\gal}+\chi_a\iO (\pham)_+ \an{c_\gal} 
\\
&\notag \qquad
+ \iO \left(h(\an{\phi_\gal})(\delta_n-\an{n_\gal})_+(1-T_{-\epsilon^{-1},1+\epsilon^{-1}}(\an{c_\gal}))\right)\dt \an{c_\gal}\\
\label{aprioridisc1}
& \qquad+\iO h(\an{\phi_\gal})(\delta_n-\an{n_\gal})_+\an{c_\gal} - \iO (\pham)_+T_{-\epsilon^{-1},1+\epsilon^{-1}}(\an{c_\gal})\dt \an{c_\gal}.
\end{align}
\an{To} bound the \an{fifth} and \andre{eighth} terms on the right hand side of the above identity, we need to obtain estimates for $|(\an{\mu_\gal})_{\Omega}|$ and $|(\an{\phi_\gal})_{\Omega}|$. Taking $v= |\Omega|^{-1}$ in \eqref{smoothregdisc2:2}, which is \an{allowed since it belongs to $V_0$}, and using the property \eqref{fpf}, we easily obtain that
\begin{equation}
    \label{massmusm}
|(\an{\mu_\gal})_{\Omega}|\leq C\left(\iO F_{\eps}(\an{\phi_\gal})+1\right).
\end{equation}
\an{Similarly}, taking $v= |\Omega|^{-1}$ in \eqref{smoothregdisc2:1}, using Assumption \eqref{def:calS} and introducing the variable $y:=(\an{\phi_\gal})_{\Omega}$ and the constant $H:=\norma{(q_{\an{\eps}}(\an{n_\gal})-\delta_n)_+h(\an{\phi_\gal})}_{L^{\infty}(Q_\an{t_1})}$, we obtain the differential Gronwall inequality
\[
-H\leq y^{\prime}+my\leq H,
\]
which gives that
\begin{equation}
    \label{massineqsm}
    y(0) e^{-m t} +\left(1-e^{-m t}\right)\left(-\frac{H}{m}\right)\leq y(t)\leq y(0) e^{-m t} +\left(1-e^{-m t}\right)\left(\frac{H}{m}\right),
\end{equation}
for every $t\in [0,t_1]$. Hence, given the assumed regularity on $\phi^0$, \andre{and the properties of the projector operator,} we obtain that
\[
|(\an{\phi_\gal})_{\Omega}|\leq C(t).
\]
Using these facts in \eqref{aprioridisc1}, together with the Poincar\'e--Wirtinger, the Cauchy--Schwarz and the Young inequalities, integrating in time \eqref{aprioridisc1} between $0$ and $t\in [0,t_1]$ and employing Assumptions \ref{ass:2:theta}--\ref{ass:parameters:5}, the regularity properties of the Moreau--Yosida approximation and of the initial data \eqref{ass:ini:weak:1}--\eqref{ass:ini:weak:3}, we obtain that
\begin{align}
&\notag \frac{1}{2}\norma{\an{\phi_\gal} \andre{(t)}}^2+\frac{1}{2}\|\nabla \an{\phi_\gal}\andre{(t)}\|^2+\iO F_{\eps}(\an{\phi_\gal}\andre{(t)})+\frac{1}{2}\norma{\pham\andre{(t)}}^2+(\chi_{\phi}^2+1)\norma{\an{n_\gal}\andre{(t)}}^2+\frac{1}{2}\norma{\nabla \an{n_\gal}\andre{(t)}}^2\\
&\notag \qquad +\frac{1}{2}\norma{\an{c_\gal}\andre{(t)}}^2+\frac{1}{2}\norma{\nabla \an{c_\gal}\andre{(t)}}^2+\int_0^t\biggl(m\norma{\an{\phi_\gal}}^2+\norma{\dt \an{n_\gal}}^2+\norma{\nabla \an{n_\gal}}^2+\norma{\dt \an{c_\gal}}^2\biggr)\\
&\notag \qquad 
+\iot \biggl(\norma{\nabla \an{c_\gal}}^2
+m_0 \norma{\nabla \an{\mu_\gal}}^2+\chi_{\phi}m_0\norma{\nabla \an{n_\gal}}^2+m_0\norma{\nabla \pham}^2\biggr)
\\ & \non \quad 
\leq C(\phi^0,\phi_a^0,n^0,c^0)+ \frac{1}{4}\norma{\an{\phi_\gal}}^2 \\
&\notag \qquad 
+ \chi_{\phi}^2\norma{\an{n_\gal}}^2
+ \int_0^t\left(\frac{m_0}{2} \norma{\nabla \an{\mu_\gal}}^2+\frac{m_0}{2}\norma{\nabla \pham}^2+\frac{1}{2}\norma{\dt \an{n_\gal}}^2+\frac{1}{2}\norma{\dt \an{c_\gal}}^2\right)\\
& \qquad \non
+C(\epsilon)\int_0^t\biggl(\frac{1}{4}\norma{\an{\phi_\gal}}^2+ \frac{1}{2}\norma{\nabla \an{\phi_\gal}}^2+\iO F_{\eps}(\an{\phi_\gal})+\frac{1}{2}\norma{\pham}^2+\norma{\an{n_\gal}}^2\biggr)
\\ & \qquad \label{aprioridisc2}
+C(\epsilon)\int_0^t\biggl(
\frac{1}{2}\norma{\nabla \an{n_\gal}}^2+\frac{1}{2}\norma{\an{c_\gal}}^2+\frac{1}{2}\norma{\nabla \an{c_\gal}}^2\biggr)+C,
\end{align}
for any $t\in [0,t_1]$. An application of the Gronwall lemma then yields\an{, for any $t\in [0,t_1]$,} that
\begin{align}
	& \non
	\norma{\an{\phi_\gal}}_{\Lt\infty V}
	+ \norma{F_{\eps}(\an{\phi_\gal})}_{\Lt\infty {\Lx1}}
	+ \norma{\pham}_{\Lt\infty H \cap \Lt2 V}
	\\ & \quad 
	+ \norma{\nabla \an{\mu_\gal}}_{\Lt2 H}
	+ \norma{\an{n_\gal}}_{\Ht1 H \cap \Lt\infty V}
	+ \norma{\an{c_\gal}}_{\Ht1 H \cap \Lt\infty V}
	\leq C.
	\label{rig:estdisc:1}
\end{align}
The bound \eqref{rig:estdisc:1}, together with \eqref{massmusm} and the Poincar\'e--Wirtinger inequality, gives that\an{, for any $t\in [0,t_1]$,}
\begin{align}
	\label{rig:estdisc:1:bisbis}
	\norma{\an{\mu_\gal}}_{\Lt2 V}\leq C.
\end{align}
\andre{A comparison argument in \eqref{smoothregdisc2:1} and \eqref{smoothregdisc2:3} finally produces\an{, for any $t\in [0,t_1]$,}
\begin{align}
	\norma{\dt \an{\phi_\gal}}_{\Lt2 \Vp}
	\leq C,
     \quad 
     \norma{\dt \pham}_{\Lt2 \Vp}
	\leq C.
	\label{rig:estdisc:2}
\end{align}}%
\an{Next, we take} $v=-\Delta \an{n_\gal}$ in \eqref{smoothregdisc2:4}, $v=-\Delta \an{c_\gal}$ in \eqref{smoothregdisc2:5} and \an{sum} the two contributions, using \eqref{rig:estdisc:1}, integrating in time between $0$ and $t\in [0,t_1]$ and employing Assumption \ref{ass:sources:3} and the regularity properties of the initial data \eqref{ass:ini:weak:3}\an{. We} obtain that
\begin{align*}
&\notag \frac{1}{2}\norma{\nabla \an{n_\gal}\andre{(t)}}^2+\frac{1}{2}\norma{\nabla \an{c_\gal}\andre{(t)}}^2+\int_0^t\left(\norma{\Delta \an{n_\gal}}^2+\norma{\Delta \an{c_\gal}}^2\right)\leq C(n^0,c^0)\\
& \notag \qquad +C\int_0^t\left(\norma{\an{\phi_\gal}}^2+\norma{\an{n_\gal}}^2+\norma{\pham}^2\right)+\frac{1}{2}\int_0^t\left(\norma{\Delta \an{n_\gal}}^2+\norma{\Delta \an{c_\gal}}^2\right)\\
& \quad 
\leq C+\frac{1}{2}\int_0^t\left(\norma{\Delta \an{n_\gal}}^2+\norma{\Delta \an{c_\gal}}^2\right),
\end{align*}
from which\andre{, using also the assumption on the initial data $c^0$ and $n^0$, along with elliptic regularity theory,} we get the bounds\an{, for any $t\in [0,t_1]$,}
\begin{equation}
        \label{rig:estdisc:delta}
	  \norma{\an{n_\gal}}_{\Lt2 W}
	+ \norma{\an{c_\gal}}_{\Lt2 W}
	\leq C.
\end{equation}
Taking now $v=-\Delta \an{\phi_\gal}$ in \eqref{smoothregdisc2:2}, we obtain that
\begin{align*}
	\iO |\Delta \an{\phi_\gal}|^2
	+ \underbrace{\iO \b'_\eps(\an{\phi_\gal})|\nabla \an{\phi_\gal}|^2}_{\geq 0}\leq \norma{\nabla \an{\mu_\gal}}\norma{\an{\phi_\gal}}-\iO \pi'(\an{\phi_\gal})|\nabla \an{\phi_\gal}|^2.
\end{align*}
In the case with a smooth potential, thanks to Assumption \eqref{ass:potsm:2}, we have that
\[
\iO |\pi'(\an{\phi_\gal})| |\nabla \an{\phi_\gal}|^2\leq C\iO\left(1+|\an{\phi_\gal}|^q\right)|\nabla \an{\phi_\gal}|^2,
\]
with $q\in [0,4)$. Observing that $\frac{4}{q}>1$ when $q<4$, using standard Sobolev embeddings, the Young inequality and \eqref{rig:estdisc:1}, we obtain that
\begin{align*}
&\int_{\Omega}|\an{\phi_\gal}|^q|\nabla \an{\phi_\gal}|^2\leq \norma{\an{\phi_\gal}}_{\an{\infty}}^q\an{\norma{{\nabla \an{\phi_\gal}}}^2}\leq C\norma{\an{\phi_\gal}}_{\an{V}}^{\frac{4+q}{2}}\left(\norma{\an{\phi_\gal}}^{\frac{q}{2}}+\norma{\Delta \an{\phi_\gal}}^{\frac{q}{2}}\right) \\
&\quad \leq C\norma{\an{\phi_\gal}}_{\an{V}}^{2+q}+\norma{\an{\phi_\gal}}_{\an{V}}^{\frac{2(q+4)}{4-q}}+\frac{1}{2}\norma{\Delta \an{\phi_\gal}}^2\leq C+\frac{1}{2}\norma{\Delta \an{\phi_\gal}}^2.
\end{align*}
In the case with a singular potential, thanks to Assumption \eqref{ass:pot:4}, we have that
\[
\iO |\pi'(\an{\phi_\gal})| |\nabla \an{\phi_\gal}|^2\leq C\norma{\nabla \an{\phi_\gal}}^2.
\]
Due to the previous computations and using  \eqref{rig:estdisc:1}, we obtain that
\begin{align*}
	\norma{\Delta \an{\phi_\gal}}^2
	\leq  C (\norma{\nabla \an{\mu_\gal}}+1).
\end{align*}
Squaring both sides, using  again \eqref{rig:estdisc:1} along with elliptic regularity theory yield\an{, for any $t\in [0,t_1]$,}
\begin{align}
	\norma{\an{\phi_\gal}}_{\an{\Lt4 {W}}}\leq C.
	\label{rig:estdisc:3:bis}
\end{align}}

The constants in the right hand side of \eqref{rig:estdisc:1}--\eqref{rig:estdisc:2} depend on the initial data, on the domain $\Omega$, on the regularization parameter $\epsilon$ but not on the discretization parameter
$\gal$. Thanks to the \last{a-priori} estimates \eqref{rig:estdisc:1}--\eqref{rig:estdisc:2}, we may extend by continuity the local solution of system \an{\eqref{smoothregdisc:1}--\eqref{discic}} to the interval $[0, T]$ and pass to the limit in a standard way as $\gal\to \infty$ in \eqref{smoothregdisc2:1}--\eqref{smoothregdisc2:5}, obtaining the existence of a weak solution $(\ph,\mu,\pha,n,c)$ to the regularized system \smoothreg\ \an{on the whole time interval $[0,T]$}. This solution has the regularity
\begin{align}
\notag \phv & \in \H1 \Vp \cap \L\infty V \cap \L4 W,  
\\ \non \mu & \in \L2 V,\\
\notag \pha  &\in \H1 \Vp \cap \L\infty H \cap \L2 V ,\\
\notag n & \in \H1 H \cap \L\infty V \cap \L2 W,\\
\label{regm} c & \in \H1 H \cap \L\infty V \cap \L2 W,
\end{align}
and satisfies the limit system
\begin{align}
 & \label{eqn:reglim:1} \<{\dt \phi},v>_V+ \iO \mobm(\phi,\phi_{a},n)\nabla  \left(\mu-\chi_{\phi}n\right) \cdot \nabla v=\iO \left((q_{\an{\eps}}(n)-\delta_n)_+h(\phi)-m\phi\right)v,
 \\ 
 &\label{eqn:reglim:2} \iO\mu\,v=\iO\nabla \phi\cdot \nabla v + \iO F_{\eps}'(\phi)v,
 \\
 \label{eqn:reglim:3}
&   \<{\dt \phi_{a}},v>_V + \iO \mobn(\phi_{a},c)\nabla \phi_{a}\cdot \nabla v
\andre{-}\chi_a\iO T_{\epsilon,\epsilon^{-1}}(\phi_{a})\mobn(\phi_{a},c)\nabla c\cdot \nabla v\\
\notag
& \quad = 
      \iO \left((T_{-\epsilon^{-1},1+\epsilon^{-1}}(c)-\delta_a)_+(1-h(\phi))+\zeta\right)({\kappa_0}\phi_{a}-{\kappa_\infty}(\phi_{a})_+^2)v,
\\
\non
&   \iO {\dt n}\, w -\iO \Delta n\, w-\chi_{\phi}\iO  p(\phi)\,w
\\ \label{eqn:reglim:4}
& \quad 
= \iO \left((1-q_{\an{\eps}}(n))(1-h(\phi)+(\phi_{a})_+)-p(\phi)q_{\an{\eps}}(n)\right)w,
\\ \label{eqn:reglim:5}
&   \iO {\dt c}\,w - \iO \Delta c \, w-\chi_a \iO (\phi_{a})_+w \\
& \notag \quad = 
 \iO \left(h(\phi)(\delta_n-n)_+(1-T_{-\epsilon^{-1},1+\epsilon^{-1}}(c))-(\phi_{a})_+T_{-\epsilon^{-1},1+\epsilon^{-1}}(c)\right)w,
\end{align}
\Accorpa\reglim {eqn:reglim:1} {eqn:reglim:5}
for \an{almost every} $t \in [0,T]$, for all $v \in V, w \in H$, with initial conditions defined in \eqref{discic}. We now want to obtain a-priori estimates for the solutions of System \reglim\ which are uniform in $\eps$, in order to study the limit problem as $\eps\to 0$ and obtain an existence result for the original System \eqnmix. In the process of obtaining these estimates we will sometimes need to consider separately the cases with a smooth potential or with a singular potential. 
We start by obtaining a maximum and a minimum principle for equation \eqref{eqn:reglim:5}, valid both in the cases with a smooth or a singular potential, which gives that $c \in [0,1]$ almost everywhere in $\an{Q}$. This condition is expected in view of the physical interpretation of $c$ as a concentration, and allows us to prove the coercivity of the chemotaxis term $-\chi_a\iO \phi_{a}c$ in the free energy of the system uniformly in $\eps$.\\ \\
\abramo{\noindent
{\bf Minimum principle.}
To begin with, let us address the minimum principle.
Let us define $f_{-}:= \andre{c_-} = - c \chi_{\{c <0\}}$ \an{and point out that} $\{c <0\}:=\{x \in \Omega : c(x) <0\}$.
Then, we take $w= - f_-$ in \eqref{eqn:reglim:5} to find that 
\begin{align*}
	& \frac12 \frac d {dt} \norma{f_-}^2
	+ \norma{\nabla (\andre{c_-})}^2
	+ \iO  h(\phi)(\delta_n-n)_+(1-T_{-\epsilon^{-1},1+\epsilon^{-1}}(c))f_-
	\\ & \quad 
	- \iO  {(\pha)_+} {T_{-\epsilon^{-1},1+\epsilon^{-1}}(c)f_-}
	=
	0.
\end{align*}
Besides, it holds that the fourth and the fifth integrals on the \lhs\ are nonnegative as well\an{. In fact, we have that}
\begin{align*}
		&\iO  h(\phi)(\delta_n-n)_+(1-T_{-\epsilon^{-1},1+\epsilon^{-1}}(c))f_-
		\\
		&\quad = \int_{\Omega\cap \{c<0\}} \underbrace{h(\phi)(\delta_n-n)_+}_{\geq 0}\underbrace{(1-T_{-\epsilon^{-1},1+\epsilon^{-1}}(c))}_{>1} \underbrace{(-c) }_{> 0}
\geq 0, 
\end{align*}
and similarly 
\begin{align*}
		&- \iO  {(\pha)_+} {T_{-\epsilon^{-1},1+\epsilon^{-1}}(c)f_-}
		=-\int_{\Omega\cap \{c<0\}} \underbrace{(\pha)_+}_{\geq 0}\underbrace{(T_{-\epsilon^{-1},1+\epsilon^{-1}}(c))}_{\leq 0} \underbrace{(-c) }_{> 0}
\geq 0. 
\end{align*}
Going back to the first identity, this entails that 
\begin{align*}
	f_- (t) =0  \quad \text{for every $t \in [0,T]$ and  $a.e. $ in $\Omega$} 
\end{align*}
from which we conclude that
\begin{align}
	\label{min:principle:c}
	c(x,t) \geq 0 \quad \text{ for $a.e. (x,t) \in \an{Q}$}.
\end{align}
\noindent
{\bf Maximum principle.}
Next, we set  $f_{+}:= (c-1)_+ = (c-1 )\chi_{\{c >1\}}$ \an{and} take $w=  f_+$ in \eqref{eqn:reglim:5} to find that
\begin{align*}
	& \frac12 \frac d {dt} \norma{f_+}^2
	+ \int_{\Omega \cap \{ c >1\}} |{\nabla c}|^2
	+ \underbrace{\int_{\Omega \cap \{c>1\}} (\pha)_+ (T_{-\epsilon^{-1},1+\epsilon^{-1}}(c)-\chi_a ) (c-1 )}_{\geq 0}\\
	& \quad + \underbrace{\int_{\Omega \cap \{c>1\}}  {h(\phv)(\delta_n-n)_+}\left(T_{-\epsilon^{-1},1+\epsilon^{-1}}(c)-1\right)(c-1)}_{\geq 0}
	=
	0,
\end{align*}
where in the third term on the left hand side of the last equality we used the compatibility condition $\chi_a \in (0,1)$ in \ref{ass:parameters:5}.
Similarly as above this yields that 
\begin{align*}
	f_+ (t) =0  \quad \text{for every $t \in [0,T]$ and  $a.e. $ in $\Omega$},
\end{align*}
meaning that 
\begin{align}
	\label{max:principle:c}
	c(x,t) \leq 1 \quad \text{ for $a.e. (x,t) \in \an{Q}$}.
\end{align}
Upon combining \eqref{min:principle:c} and \eqref{max:principle:c}, we finally infer that 
\begin{align}
\label{minmax:principle:c}
		0 \leq c(x,t) \leq 1 \quad \text{ for $a.e. (x,t) \in \an{Q}$},
\end{align}
uniformly in $\eps$. Note that, \andre{as a consequence of} \eqref{minmax:principle:c} and the definition \eqref{Tlm}, we have that
\[
T_{-\epsilon^{-1},1+\epsilon^{-1}}(c)\equiv c,
\]
uniformly in $\eps$.
%
\begin{remark}
    \label{phia+}
    We observe that in \eqref{eqn:reglim:5} the choice of taking the positive part $(\phi_a)_+$ in the chemotactic and in the source terms implies the validity of the minimum and the maximum principles for $c$ for any value of $\epsilon$. Indeed, given $\an{\eps \in(0,1)}$, a solution $\phi_a$ for the regularized system \reglim\ is not necessarily \an{non}negative, hence to enforce the minimum and maximum principles for $c$ we need to truncate $\phi_a$ to \an{non}negative values in the chemotactic and in the source terms.
\end{remark}
In the case with a singular potential it is possible to obtain also a maximum and a minimum principle for equation \eqref{eqn:reglim:4}. This condition is also expected in view of the physical interpretation of $n$ as a concentration. \andre{Recalling} $p(\ph)\an{=}(\ph)_+$, $q_{\an{\eps}}(n)\an{=}T_{-\epsilon^{-1},1+\epsilon^{-1}}(n)$, $\chi_{\ph}\in (0,1)$\andre{,} and observing that $1-h(\ph)+(\pha)_+\geq 0$, we may apply the same calculations as the one employed to prove the maximum and minimum principles for equation \eqref{eqn:reglim:5}, obtaining that
\begin{align}
\label{minmax:principle:n}
		0 \leq n(x,t) \leq 1 \quad \text{ for $a.e. (x,t) \in \an{Q}$}.
\end{align}
We again observe that the choice $p_{\epsilon}(\ph)\an{=}(\ph)_+$ implies the validity of the minimum and the maximum principles \andre{also} \an{at the approximation level}\andre{, where the singular potentials are approximated by polynomial type potentials. In fact, given $\an{\eps \in(0,1)}$, a solution $\ph$ is not necessarily \an{non}negative, nor confined in the physical range $[0,1]$.} In the limit $\eps\to 0$ the solution $\ph$ will turn out to be nonnegative only in the case with a singular potential, so the property $n\in [0,1]$ is valid only in the latter case.
Due to \eqref{minmax:principle:n} and the definition \eqref{Tlm}, we have that
\[
T_{-\epsilon^{-1},1+\epsilon^{-1}}(n)\equiv n,
\]
uniformly in $\eps$, and accordingly
\begin{equation*}
q_{\an{\eps}}(n)\an{=}q(n)\an{=}\begin{cases}
h(n) & \text{if $F$ is smooth},\\
n &\text{if $F$ is singular},
\end{cases}
\end{equation*}
with $q(n)\in [0,1]$.
We then redefine the source terms in \reglim\ as \andre{follows}
\begin{align}
 & \label{eqn:sourcerl:1}  \Sv(\ph,n)=(q(n)-\delta_n)_+h(\phi)-m\phi\an{=:}P(\ph,n)-m\ph,
 \\ 
 \non
&   \Sa(\ph,\pha,c) = \big((c-\delta_a)_+(1-h(\phi))+\zeta\big)({\kappa_0}\phi_{a}-{\kappa_\infty}(\phi_{a})_+^2)
\\ & \label{eqn:sourcerl:2}
 \phantom{\Sa(\ph,\pha,c)}\;
\an{=} \,\theta(\ph,c)({\kappa_0}\phi_{a}-{\kappa_\infty}(\phi_{a})_+^2),
\\
\label{eqn:sourcerl:3}
&   \Sn(\ph,\pha,\andre{n})= (1-q(n))(1-h(\phi)+(\phi_{a})_+)-p(\phi)q(n),
\\ \label{eqn:sourcerl:4}
&   \Sc(\ph,\pha,n,c) = h(\phi)(\delta_n-n)_+(1-c)-(\phi_{a})_+c,
\end{align}
\Accorpa\sourcerl {eqn:sourcerl:1} {eqn:sourcerl:4}
and observe that there exist \an{constants} $C_1,C_2,C_3\geq 0$ such that
\begin{align}
	& |P(\ph,n)|\leq C_1, \quad  \ph,n\in \mathbb{R},
	\qquad
	\label{ass:Source:1}
	\\
	& \zeta\leq \an{|}\theta(\ph,c)|\leq 1+\zeta, \quad  \ph\in \mathbb{R}, c\in [0,1],
	\qquad
	\label{ass:Source:2}
	\\
	& |\Sn(\ph,\pha,n)|\leq C_2\left(|\ph|+(\pha)_++1\right), \quad  \ph,\pha,n\in \mathbb{R},
	\qquad
	\label{ass:Source:3}
	\\
	&  |\Sc(\ph,\pha,n,c)|\leq C_3\left((\pha)_++|n|+1\right), \quad  \ph,\pha,n\in \mathbb{R}, c\in [0,1].
	\qquad
	\label{ass:Source:4}
\end{align}
\Accorpa \Source {ass:Source:1} {ass:Source:4}
We now \an{move} to obtain a-priori estimates for System \reglim.
\begin{remark}
    \label{remapriorieps}
    In order to rigorously obtain a-priori estimates uniform in the regularization parameter $\epsilon$ for System \reglim, we should need to consider a time regularization of \reglim\ with time regularized functions $\phi_\an{\tau},\phi_{a\an{\tau}}$, depending on a regularization parameter $\an{\tau}$, where, given a function $u:Q\rightarrow \mathbb{R}$, we define:
    \[
    u_\an{\tau}\an{(x,t)}:=\frac{1}{\an{\tau}}\int_{t-\an{\tau}}^{t}u \an{(x, \tau)}\,d\tau,
    \]
    with $u_\an{\tau}\an{(x,t)}\an{:}=u^0(\an{x})$ for $t\leq 0$.
    In this way, if $u\in\an{\H1 \Vp \cap  \L2 V }$, we have that
    \[
    \int_0^T\<{\dt u_\an{\tau}},v>_V=\int_Q \dt u_\an{\tau} v.
    \]
    Since, with the given regularities of $\phi$ and $\phi_a$ in \eqref{regm}, it \an{readily follows that} $\phi_\an{\tau}\rightarrow \phi$ strongly in $\C0 V$, $\phi_{a\an{\tau}}\rightarrow \phi_a$ strongly in $\L2 V$ and $\dt \phi_\an{\tau} \rightarrow \dt \phi$, $\dt \phi_{a\an{\tau}} \rightarrow \dt \phi_a$ strongly in $\L2 \Vp$ as $\an{\tau}\to 0$, we should easily pass to the limit as $\an{\tau}\to 0$ in the aforementioned $\an{\tau}$-time regularized version of \reglim. Since this procedure is standard\an{,} see, e.g., \cite[Lemma 2]{ELG}\an{,} in the following we implicitly assume to have performed a time regularization of \reglim\ to obtain a-priori estimates, avoiding to report all the details for \an{simplicity}.
\end{remark}
}%
\step First estimate

{We now take $v= \mu-\chi_{\phi}n+\phi$ in \eqref{eqn:reglim:1}, $v= -\dt \phi$ in \eqref{eqn:reglim:2}, $v= \mathcal{E}_{\epsilon,\epsilon^{-1}}^{\prime}(\phi_{a})-\chi_ac$ in \eqref{eqn:reglim:3}, $w= \dt n+2(\chi_{\phi}^2+1)n$ in \eqref{eqn:reglim:4} and $w= \dt c+c$ in \eqref{eqn:reglim:5}, sum all the equations and \an{rearrange} some terms as in \eqref{aprioridisc1}. \andre{Using} \an{the identity}
\begin{equation}
\label{tgrad}
T_{\epsilon,\epsilon^{-1}}(\pha)\nabla \mathcal{E}_{\epsilon,\epsilon^{-1}}^{\prime}(\pha) = \nabla \pha,
\quad \andre{\eps \in (0,1),}
\end{equation}
we \andre{find} that
\begin{align}
&\notag \frac{d}{dt} \biggl(\frac{1}{2}\norma{\phi}^2+\frac{1}{2}\norma{\nabla \phi}^2+\iO F_{\eps}(\phi)-\chi_{\phi}\iO \phi\,n+\iO\mathcal{E}_{\epsilon,\epsilon^{-1}}(\pha)-\chi_{a}\iO \pha\,c\\
&\notag \qquad+(\chi_{\phi}^2+1)\norma{n}^2+ \frac{1}{2}\norma{\nabla n}^2+\frac{1}{2}\norma{c}^2+\frac{1}{2}\norma{\nabla c}^2\biggr)+m\norma{\phi}^2+\norma{\dt n}^2\\
&\notag \qquad +2(\chi_{\phi}^2+1)\norma{\nabla n}^2+\norma{\dt c}^2+\norma{\nabla c}^2 + \iO \mobm(\phi,\phi_{a},n)\nabla \mu\cdot \nabla \mu\\
&\notag \qquad +\chi_{\phi}\iO \mobm(\phi,\phi_{a},n)\nabla  n\cdot \nabla n+ {\kappa_\infty \iO \theta(\ph,c)(\phi_{a})_+^2\mathcal{E}^{\prime}(\pha)}\\
&\notag \qquad + \iO T_{\epsilon,\epsilon^{-1}}(\pha)\mobn(\phi_{a},c)\nabla \left(\mathcal{E}_{\epsilon,\epsilon^{-1}}^{\prime}(\pha)-\chi_ac\right)\cdot \nabla \left(\mathcal{E}_{\epsilon,\epsilon^{-1}}^{\prime}(\pha)-\chi_ac\right)\\
& \notag \qquad + \underbrace{2(\chi_{\phi}^2+1)\iO \left(1-h(\ph)+(\pha)_+\right)q(n)n}_{\geq 0}+\underbrace{\iO \left(h(\phi)(\delta_n-n)_++(\phi_{a})_+\right)c^2}_{\geq 0}\\
&\notag \quad = 2\chi_{\phi}\iO \mobm(\phi,\pha,n)\nabla  \mu\cdot \nabla n-\iO \mobm(\phi,\pha,n)\nabla \mu\cdot \nabla \phi\\
&\notag \qquad + \chi_{\phi}\iO \mobm(\phi,\pha,n)\nabla n\cdot \nabla \phi+\iO P(\phi,n)(\mu-\an{\mu_{\Omega}})+\iO P(\phi,n)\an{\mu_{\Omega}}\\
&\notag \qquad - \chi_{\phi}\iO P(\phi,n)n+\iO P(\phi,n)\phi-m\iO \phi(\mu-\an{\mu_{\Omega}})-m|\Omega|\an{\ph_{\Omega}} \, \an{\mu_{\Omega}}\\
&\notag \qquad + \chi_{\phi}(m+2(\chi_{\phi}^2+1))\iO \phi n+\chi_a\iO \left((\pha)_+-\pha\right)\dt c+  \kappa_0\iO \an{\theta(\ph,c)}\pha\mathcal{E}^{\prime}(\pha)\\
&\notag\qquad -\kappa_0\chi_a\iO \an{\theta(\ph,c)}\pha c+ \kappa_\infty \chi_a \iO \an{\theta(\ph,c)}(\phi_{a})_+^2c +\iO \Sn(\ph,\pha,n)\dt n\\
&\notag\qquad+ 2(\chi_{\phi}^2+1)\iO (1-h(\phi)+(\pha)_+)n+\chi_a\iO (\pha)_+ c-2(\chi_{\phi}^2+1)\iO p(\phi)q(n)n\\
&\non \qquad +\chi_{\phi}\iO \left(p(\phi)-\phi\right)(\dt n+2(\chi_{\phi}^2+1)n)
+ \iO \Sc(\ph,\pha,n,c)\dt c
\\ &  \label{apriorieps1}
\qquad 
+\iO h(\phi)(F\delta_n-n)_+c.
\end{align}
\an{Here, we} need uniform estimates for $\andre{{\mu}_{\Omega}}$ and $\andre{{\phi}_{\Omega}}$ in order to bound the \an{fifth} and \an{ninth} terms on the right hand side of \eqref{apriorieps1}. We obtain them separately for the case with a smooth potential and the case with a singular potential. In the former case, we proceed with similar arguments as in \eqref{massmusm} and \eqref{massineqsm}, taking $v= |\Omega|^{-1}$ in \eqref{eqn:reglim:1} and \eqref{eqn:reglim:2} and using Assumption \eqref{ass:potsm:1}. We \andre{are then lead to}
\begin{equation}
    \label{massmulim}
    |\an{\mu_{\Omega}}|\leq C\left(\iO F_{\eps}(\phi)+1\right),
\end{equation}
and
\begin{equation}
    \label{massineqlim}
    |\andre{\phi_{\Omega}(t)}|\leq C
    \quad 
    \text{for every $t\in [0,T]$,}
\end{equation}
uniformly in $\epsilon$. For what concerns the \andre{singular potential case}, the compatibility condition in \ref{ass:sources:3} plays a crucial role in constraining the mass dynamics.
Setting $y : = \an{\ph_\Omega}$, $H:=\norma{P(\ph,n)}_{L^{\infty}(\an{Q})}$ and testing  \eqref{eqn:reglim:1} by $ \an{v=}{|\Omega|}^{-1}$, we arrive at the inequalities \an{in} \eqref{massineqsm}, from which, thanks to the compatibility   in \ref{ass:sources:3}, \andre{produces}
\begin{align}\label{mass}
	\andre{\phi_{\Omega}(t)} \quad \text{belongs to the interior of $D(\beta)$ for every $t \in [0,T].$}
\end{align}
To control the mean of $\mu$, we test \eqref{eqn:reglim:2} by $v\andre{=} 1$ \an{to} find that
\begin{align}
\label{MZa}
	|\Omega| |\mu_\Omega| 
	\leq \norma{F_\eps'(\phT)}_1.
\end{align}
Let us notice from Assumptions \eqref{ass:pot:1}--\eqref{ass:pot:4} \andre{it holds} that 
\begin{align}
\label{MZb}
  \norma{F_\eps'(\phT)}_1
 	\leq 
 	\norma{\betaeps (\phT)}_1 +\norma{\pi(\phT)}_1 
 	\leq 
 	\norma{\betaeps (\phT)}_1 + C (\norma{\phT}^2 +1) .	
\end{align}
Thus, it is enough to control the term involving the regularized singular part $\betaeps$.
We test \eqref{eqn:reglim:2} by $\phT- \an{\ph_\Omega}$ to find that
\begin{align}
\label{MZ0}
	\iO \betaeps(\phT) (\phT- \an{\ph_\Omega}) 
	+ \iO \pi(\phT) (\phT- \an{\ph_\Omega}) 
	+ \norma{\nabla \phT}^2
	=  \iO \mu (\phT- \an{\ph_\Omega}).
\end{align}
On the other hand, using the mass property in \eqref{mass} and arguing as in \cite{MZ}, we find positive constants $C_F$ and $c_F$ such that
\begin{align}
    \label{MZ}
	\iO \betaeps(\phT) (\phT- \an{\ph_\Omega}) 
	\geq C_F \norma{\betaeps(\phT)}_1
	- c_F.
\end{align}
Thus, using \eqref{MZ0} in \eqref{MZ} and Assumption \eqref{ass:pot:4}, the Poincar\'e--Wirtinger and the Young inequalities we  deduce that
\begin{align}
	\label{est:MZ:1}
	C_F \norma{\betaeps(\phT)}_1 
	& \leq 
	C \norma{\nabla \mu }\norma{\nabla \phT}
	+c_F
	+ C (\norma{\phT}^2+1).
\end{align}
Hence, collecting \eqref{ass:Source:1} and \eqref{MZa}--\eqref{est:MZ:1} and using the Young inequality, we obtain that
\begin{equation}
\label{est:MZ:2}
\iO P(\phi,n)\an{\mu_{\Omega}}\leq C_1|\Omega| |\mu_\Omega| \leq C_1\norma{F_\eps'(\phT)}_1\leq \frac{m_0}{4}\norma{\nabla \mu}^2+C\norma{\nabla \ph}^2+C (\norma{\phT}^2+1).
\end{equation}
Using \andre{\eqref{minmax:principle:c}, \eqref{minmax:principle:n}, and \eqref{mass}, \eqref{est:MZ:2}} in \eqref{apriorieps1}, together with the Poincar\'e--Wirtinger, the Cauchy--Schwarz and the Young inequalities, integrating in time between $0$ and $t\in [0,T]$ and employing Assumptions \ref{ass:2:theta}--\ref{ass:parameters:5} and the regularity properties of the initial data \eqref{ass:ini:weak:1}--\eqref{ass:ini:weak:3}, we obtain that
\begin{align}
&\notag \frac{1}{2}\norma{\phi\an{(t)}}^2+\frac{1}{2}\|\nabla \phi\an{(t)}\|^2+\iO F_{\eps}(\phi\an{(t)})+\iO\mathcal{E}_{\epsilon,\epsilon^{-1}}(\pha\an{(t)})
\\ & \non
\qquad 
-\chi_{a}\iO \pha\an{(t)}\,c\an{(t)}+(\chi_{\phi}^2+1)\norma{n\an{(t)}}^2+\frac{1}{2}\norma{\nabla n\an{(t)}}^2\\
&\notag \qquad +\frac{1}{2}\norma{c\an{(t)}}^2+\frac{1}{2}\norma{\nabla c\an{(t)}}^2+\int_0^t\biggl(m\norma{\phi}^2+\norma{\dt n}^2+\norma{\nabla n}^2+\norma{\dt c}^2+\norma{\nabla c}^2\\
&\notag \qquad +m_0 \norma{\nabla \mu}^2+\chi_{\phi}m_0\norma{\nabla n}^2 \biggr)+ {\kappa_\infty \int_0^t\iO \an{\theta(\ph,c)}(\phi_{a})_+^2\mathcal{E}_{\epsilon,\epsilon^{-1}}^{\prime}(\pha)}\\
&\notag \qquad + \int_0^t\iO T_{\epsilon,\epsilon^{-1}}(\pha)\mobn(\phi_{a},c)\nabla \left(\mathcal{E}_{\epsilon,\epsilon^{-1}}^{\prime}(\pha)-\chi_ac\right)\cdot \nabla \left(\mathcal{E}_{\epsilon,\epsilon^{-1}}^{\prime}(\pha)-\chi_ac\right)\\
&\notag \quad \leq C(\phi^0,\phi_a^0,n^0,c^0)+ \frac{1}{4}\norma{\phi\andre{(t)}}^2+ \chi_{\phi}^2\norma{n\andre{(t)}}^2 + \int_0^t\left(\frac{m_0}{2} \norma{\nabla \mu}^2+\frac{1}{4}\norma{\dt n}^2+\frac{1}{4}\norma{\dt c}^2\right)\\
&\notag \qquad +C\int_0^t\biggl(\frac{1}{4}\norma{\phi}^2+ \frac{1}{2}\norma{\nabla \phi}^2+\iO F_{\eps}(\phi)+\frac{1}{2}\norma{n}^2+\frac{1}{2}\norma{\nabla n}^2\biggr)+C\\
&\notag  \qquad +\underbrace{\chi_a\int_0^t\iO \left((\pha)_+-\pha\right)\dt c}_{\an{=:}{\andre{\I}}_1}+  \underbrace{\kappa_0\int_0^t\iO \an{\theta(\ph,c)}\pha\mathcal{E}_{\epsilon,\epsilon^{-1}}^{\prime}(\pha)}_{\an{=:}{\andre{\I}}_2}-\underbrace{\kappa_0\chi_a\int_0^t\iO \an{\theta(\ph,c)}\pha c}_{\an{=:}{\andre{\I}}_3}\\
&\notag\qquad + \underbrace{\kappa_\infty \chi_a \int_0^t\iO \an{\theta(\ph,c)}(\phi_{a})_+^2c}_{\an{=:}{\andre{\I}}_4}+\underbrace{\iO \Sn(\ph,\pha,n)\dt n+ \iO \Sc(\ph,\pha,n,c)\dt c}_{\an{=:}{\andre{\I}}_5} \\
\label{apriorieps2}
& \qquad+\underbrace{\chi_a\int_0^t\iO (\pha)_+ c+ 2(\chi_{\phi}^2+1)\iO (\pha)_+n +\int_0^t\iO (\pha)_+c \,\dt c}_{\an{=:{\andre{\I}}_6}},
\end{align}
for any $t\in [0,T]$. 
The term ${\andre{\I}}_1$ can be bounded using \eqref{propTE1}, the Cauchy--Schwarz and the Young inequality\an{. Namely, it holds that}
\[
|{\andre{\I}}_1|\leq \chi_a\int_0^t\iO |(\pha)_-|\,|\dt c|\leq \frac{1}{8}\int_0^t\norma{\dt c}^2+4\chi_a^2\epsilon \int_0^t\iO\mathcal{E}_{\epsilon,\epsilon^{-1}}(\pha).
\]
Using \eqref{propTE3} and \eqref{minmax:principle:c}, we can bound ${\andre{\I}}_2$ as
\[
|{\andre{\I}}_2|\leq 2\kappa_0(\zeta+1)\int_0^t \iO\mathcal{E}_{\epsilon,\epsilon^{-1}}(\pha)+\kappa_0(\zeta+1)|\Omega|\an{T}.
\]
For what concerns the term ${\andre{\I}}_3$, we rewrite it as
\[
{\andre{\I}}_3=-\underbrace{\kappa_0\chi_a\int_0^t\iO \an{\theta(\ph,c)}(\pha)_+ c}_{\geq 0}+\kappa_0\chi_a\int_0^t\iO \an{\theta(\ph,c)}(\pha)_- c,
\]
hence, using \eqref{propTE1}, \eqref{minmax:principle:c}, the Cauchy--Schwarz and the Young inequality, we obtain that
\[
|{\andre{\I}}_3|\leq \kappa_0\chi_a(\zeta+1)\epsilon\int_0^t \iO\mathcal{E}_{\epsilon,\epsilon^{-1}}(\pha)+\frac{\kappa_0\chi_a(\zeta+1)}{2}|\Omega|\an{T}.
\]
The term ${\andre{\I}}_4$ can be bounded using \eqref{sleqe3} and \eqref{minmax:principle:c}, \an{leading to}
\[
|{\andre{\I}}_4|\leq \kappa_\infty \chi_a (1+\zeta)\int_0^t\iO (\phi_{a})_+^2\leq \frac{\kappa_\infty \zeta}{4} \int_0^t\iO \left((\phi_{a})_+^2\mathcal{E}_{\epsilon,\epsilon^{-1}}^{\prime}(\pha)+(2\an{e})^{-1}\right)+C|\Omega|\an{T}.
\]
Finally, ${\andre{\I}}_5$ can be bounded using \eqref{sleqe3}, \eqref{minmax:principle:c}, \eqref{ass:Source:3}, \eqref{ass:Source:4},  the Cauchy--Schwarz and the Young inequality, obtaining that
\begin{align*}
|{\andre{\I}}_5| \an{+|{\andre{\I}}_6|} & \leq \frac{1}{4}\int_0^t\norma{\dt n}^2+\frac{1}{8}\int_0^t\norma{\dt c}^2+C\int_0^t\left(\norma{\ph}^2+\norma{n}^2+1\right)\\
& \quad +\frac{\kappa_\infty \zeta}{4}\int_0^t\iO \left((\phi_{a})_+^2\mathcal{E}_{\epsilon,\epsilon^{-1}}^{\prime}(\pha)+(2\an{e})^{-1}\right)+C|\Omega|\an{T}.
\end{align*}
Thanks to \eqref{sleqe1}, we can treat the chemotactic term in \eqref{apriorieps2} by \an{noticing} that
\[
\chi_a\iO \pha\,c\leq \chi_a\iO|\pha|\leq \chi_a\iO \mathcal{E}_{\epsilon,\epsilon^{-1}}(\pha)+\chi_a(\an{e}-1)|\Omega|.
\]
Using the previous results in \eqref{apriorieps2}, adding to both sides the quantity 
\[\frac{\kappa_\infty}{2\an{e}} \int_0^t\iO \an{\theta(\ph,c)}\] 
and considering \eqref{sleqe2}, we obtain that
\begin{align}
&\notag \frac{1}{4}\norma{\phi\an{(t)}}^2+\frac{1}{2}\|\nabla \phi\an{(t)}\|^2+\iO F_{\eps}(\phi\an{(t)})+(1-\chi_a)\iO\mathcal{E}_{\epsilon,\epsilon^{-1}}(\pha\an{(t)})
\\ & \qquad \non
+\norma{n\an{(t)}}^2
+\frac{1}{2}\norma{\nabla n\an{(t)}}^2+\frac{1}{2}\norma{c\an{(t)}}^2_\andre{V}
\\
&\notag \qquad +\int_0^t\biggl(m\norma{\phi}^2+\frac{1}{2}\norma{\dt n}^2+\norma{\nabla n}^2+\frac{1}{2}\norma{\dt c}^2+\norma{\nabla c}^2+\frac{m_0}{2} \norma{\nabla \mu}^2+\chi_{\phi}m_0\norma{\nabla n}^2 \biggr)\\
&\notag \qquad + \frac{\kappa_\infty \zeta}{2} \underbrace{\int_0^t\iO \left((\phi_{a})_+^2\mathcal{E}_{\epsilon,\epsilon^{-1}}^{\prime}(\pha)+(\andre{2e})^{-1}\right)}_{\geq 0}+\underbrace{\kappa_0\chi_a\int_0^t\iO \an{\theta(\ph,c)}(\pha)_+ c}_{\geq 0}\\
&\notag \qquad + \int_0^t\iO T_{\epsilon,\epsilon^{-1}}(\pha)\mobn(\phi_{a},c)\nabla \left(\mathcal{E}_{\epsilon,\epsilon^{-1}}^{\prime}(\pha)-\chi_ac\right)\cdot \nabla \left(\mathcal{E}_{\epsilon,\epsilon^{-1}}^{\prime}(\pha)-\chi_ac\right)\\
&\notag \quad \leq C(\phi^0,\phi_a^0,n^0,c^0\an{)}+C+C\int_0^t\biggl(\frac{1}{4}\norma{\phi}^2+ \frac{1}{2}\norma{\nabla \phi}^2+\iO F_{\eps}(\phi) \an{\biggr)}\\
\label{apriorieps3}
& \qquad 
\an{ +C \iot \biggl( }(1-\chi_a)\iO\mathcal{E}_{\epsilon,\epsilon^{-1}}(\pha)
+\norma{n}^2+\frac{1}{2}\norma{\nabla n}^2+\frac{1}{2}\norma{c}^2+\frac{1}{2}\norma{\nabla c}^2\biggr),
\end{align}
for any $t\in [0,T]$. 
An application of the Gronwall lemma and of the properties \eqref{sleqe1}, \eqref{sleqe3}, then yields that
\begin{align}
	& \non
	\norma{\phi}_{\L\infty V}
	+ \norma{F_{\eps}(\phi)}_{\L\infty {\Lx1}}
	+ \norma{\pha}_{\L\infty {\Lx1}}
        + \norma{(\pha)_+}_{\L2 H}
	\\ & \quad 
	+ \norma{\nabla \mu}_{\L2 H}
	+ \norma{n}_{\H1 H \cap \L\infty V}
	+ \norma{c}_{\H1 H \cap \L\infty V}
	\leq C.
	\label{rig:estreg:1}
\end{align}
Also, as a consequence of property \eqref{propTE1} and of \eqref{apriorieps3}, we have that
\begin{equation}
    \label{rig:estreg:1bis}
    \norma{(\pha)_-}_{\L\infty H}\leq 2\epsilon \iO\mathcal{E}_{\epsilon,\epsilon^{-1}}(\pha)\leq C\epsilon.
\end{equation}
\an{Combining} \eqref{rig:estreg:1} \an{with} \eqref{rig:estreg:1bis}\an{,} we get that
\begin{equation}
    \label{rig:estreg:1tris}
    \norma{\pha}_{\L2 H}\leq C.
\end{equation}
In the case with a smooth potential, the bound \eqref{rig:estreg:1}, together with \eqref{massmulim}, gives that $\norma{\mu_{\Omega}}_{L^{\infty}(0,T)}\leq C$\an{, whereas, in} the case with a singular potential, the bounds \eqref{MZa}, \eqref{MZb}, \eqref{est:MZ:1} and \eqref{rig:estreg:1}  \an{just give} that $\norma{\mu_{\Omega}}_{L^{2}(0,T)}\leq C$. Hence, the bound \eqref{rig:estreg:1}, together with the Poincar\'e--Wirtinger inequality, gives\an{, in both cases,} that
\begin{align}
	\label{rig:estreg:1:bisbis}
	\norma{\mu}_{\an{\L2 V}}\leq C.
\end{align}

\step
Second estimate

\an{Next, taking} $w=-\Delta n$ in \eqref{eqn:reglim:4}, $w=-\Delta c$ in \eqref{eqn:reglim:5}, which are
feasible test functions \andre{due} \an{to} \eqref{regm}, and summing the two contributions, using \eqref{rig:estreg:1} and \eqref{rig:estreg:1tris}, \an{integrating} in time between $0$ and $t\in [0,T]$ and employing \eqref{ass:Source:3}, \eqref{ass:Source:4} and the regularity properties of the initial data \eqref{ass:ini:weak:3}, we obtain that
\begin{align*}
\notag 
& \frac{1}{2}\norma{\nabla n\andre{(t)}}^2+\frac{1}{2}\norma{\nabla c\andre{(t)}}^2+\int_0^t\left(\norma{\Delta n}^2+\norma{\Delta c}^2\right)\leq C(n^0,c^0)
\\ & \notag \qquad +C\int_0^t\left(\norma{\phi}^2+\norma{n}^2+\norma{\pha}^2\right)+\frac{1}{2}\int_0^t\left(\norma{\Delta n}^2+\norma{\Delta c}^2\right)\\
& \quad 
\leq C+\frac{1}{2}\int_0^t\left(\norma{\Delta n}^2+\norma{\Delta c}^2\right),
\end{align*}
from which\an{, using also elliptic regularity,} we get the bounds
\begin{equation}
        \label{rig:estreg:delta}
	  \norma{n}_{\an{\L2 W}}
	+ \norma{c}_{\an{\L2 W}}
	\leq C.
\end{equation}
\step
Third estimate

A comparison argument in \eqref{eqn:reglim:1}
\andre{then} produces
\begin{equation}
        \label{rig:estreg:2}
	\norma{\dt \phi}_{\an{\L2 \Vp}}
	\leq C.
\end{equation}
}
\step
Fourth estimate

We observe that the estimate \eqref{rig:estdisc:3:bis} is uniform in $\eps$, so we get that
\begin{align}
	\norma{\phT}_{\L4 {W}}\leq C.
	\label{rig:est:3:bis}
\end{align}
In the case with a smooth potential, thanks to \eqref{rig:est:3:bis}, given the assumed regularity of $F$ and the bound \eqref{rig:estreg:1:bisbis}, we have that $\norma{\mu-F'(\ph)}_{\L2 {\Lx \esp}}\leq C$, with $\esp$ as defined in \eqref{def:esp}. Hence, elliptic regularity theory applied to \eqref{eqn:reglim:2} gives that 
\begin{align}
	 \norma{\phT}_{\L2 {\Wx{2,\esp}}}
	\leq C.
	\label{rig:est:3:sm}
\end{align}
In the case with a singular potential, we can consider \eqref{eqn:reglim:2}
\an{as a family of time-dependent elliptic problems} with maximal monotone perturbations as follows:
\begin{align*}
	\begin{cases}
	-\Delta \phT 
	+ \betaeps(\phT) 
	 =  f_\ph:= 
	 \mu - \pi( \phT)
	 \quad 
	 & \text{in $\an{\Omega},$}
	 \\
	 \dn \phT =0 \quad & \text{on $\an{\Gamma}.$}
	 \end{cases}
\end{align*}
Since Assumption \eqref{ass:pot:4} and the above computations ensure that the forcing term $f_\ph \in \L2 V$, standard arguments allow us to infer that
\begin{align}
	 \norma{\phT}_{\L2 {\Wx{2,\esp}}}
	+ \norma{\betaeps(\phT)}_{\L2 {\Lx \esp}}
	\leq C\an{,}
	\label{rig:est:3}
\end{align}
where we also use elliptic regularity theory and the continuous embedding $V \emb \Lx \esp$ with $\esp$ as defined in \eqref{def:esp}.

\step 
Fi\an{f}th estimate

To conclude, we test \eqref{eqn:reglim:3} by $\mathcal{E}_{\epsilon,\epsilon^{-1}}^{\prime}(\pha)$ to obtain that
\begin{align*}
	 & \frac d {dt} \iO \mathcal{E}_{\epsilon,\epsilon^{-1}}(\pha)
	 + m_0 \iO \mathcal{E}_{\epsilon,\epsilon^{-1}}''(\pha)|\nabla \phi_a|^2
	 + \zeta \kappa_\infty \iO (\ph_a)_+^2 \mathcal{E}_{\epsilon,\epsilon^{-1}}^{\prime}(\pha)
	 \\ & \quad 
	\leq 
	- \chi_a \iO   \pha \, \Delta c 
	+ \kappa_0(1+\zeta) \iO  \pha \mathcal{E}_{\epsilon,\epsilon^{-1}}^{\prime}(\pha)
	\\ & \quad 
	\leq 
	C( \norma{\Delta c}^2
	+  \norma{\pha}^2 )
	+ \kappa_0 (1+\zeta) \iO |\pha| \mathcal{E}_{\epsilon,\epsilon^{-1}}^{\prime}(\pha)
	\\ & \quad 
	 \leq 
	 2\kappa_0 (1+\zeta) \iO \mathcal{E}_{\epsilon,\epsilon^{-1}}(\pha)
	 + C\leq C.
\end{align*}
This allows us to deduce the additional bound
\an{
\begin{align}
	\label{rig:est:4}
	\norma{({\mathcal{E}_{\epsilon,\epsilon^{-1}}''(\pha)})^{1/2} \, \nabla \phi_a}_{\L2 H} \leq C.
\end{align}
}

\step 
Sixth estimate

Let us now obtain some information on the time derivative of $\pha$. In this direction, let us notice that \eqref{apriorieps3} yields, recalling \eqref{mob:bou}, in particular that
\begin{align}
	& \non
          \norma{T_{\epsilon,\epsilon^{-1}}(\pha) \mobn(\pha,c)\nabla (\mathcal{E}_{\epsilon,\epsilon^{-1}}^{\prime}(\pha) -\chi_a c)}_{L^{\frac{4}{3}}(Q)}\\
	&\quad \label{rig:est:5} 
           \leq 
	M\norma{\andre{\left(T_{\epsilon,\epsilon^{-1}}(\pha)\right) ^{1/2}}}_{L^{4}(Q)}\norma{\left(T_{\epsilon,\epsilon^{-1}}(\pha)\right) ^{1/2}\nabla (\mathcal{E}_{\epsilon,\epsilon^{-1}}^{\prime}(\pha) -\chi_a c)}_{L^{2}(Q)}
	\leq C.
\end{align}
On the other hand, \andre{we obtain that, for $z \in \Wx{1,{4}}$,} it holds that
\begin{align*}
	&\iO  T_{\epsilon,\epsilon^{-1}}(\pha) \mobn(\ph_a,c) \nabla (\mathcal{E}_{\epsilon,\epsilon^{-1}}^{\prime}(\pha) - \chi_a c)\cdot \nabla z
	\\ & \quad \leq 
	\norma{T_{\epsilon,\epsilon^{-1}}(\pha) \mobn(\pha,c)\nabla (\mathcal{E}_{\epsilon,\epsilon^{-1}}^{\prime}(\pha) -\chi_a c)}_{\frac {4}{3}} \norma{\nabla z}_{4}
	\\ & \quad \leq 
	M\norma{\andre{T_{\epsilon,\epsilon^{-1}}(\pha)}}^{1/2}\norma{(T_{\epsilon,\epsilon^{-1}}(\pha) ) ^{1/2} \, \nabla (\mathcal{E}_{\epsilon,\epsilon^{-1}}^{\prime}(\pha) -\chi_a c)} \norma{z}_{\andre{\Wx{1,{4}}}}
	\leq C.
\end{align*}
Besides, \andre{we owe to the continuous embedding} $\andre{\Wx{1,{4}}} \emb \Lx{\infty}$, \andre{to derive that}
\begin{align*}
	\iO \theta(\phT,c) (\kappa_0 \pha - \kappa_\infty \an{(\ph_a)_+^2})z
	\leq 
	C ( 1+ \norma{\ph_a}^2_2)\norma{z}_\infty 
	\leq 
	C ( 1+ \norma{\pha}^2_2)\norma{z}_{\andre{\Wx{1,{4}}}}.
\end{align*}
Combining the above estimate it is then a standard matter to derive from \abramob{\eqref{eqn:reglim:3}} that
\begin{align}
	\label{rig:est:6}
		\norma{\dt \pha}_{\L1 {(\andre{\Wx{1,{4}}})^*}} \andre{\leq C}.
\end{align}

\subsection{Passing to the Limit} 

In this section, we  aim at detailing the passage to the limit $\eps\to 0$. Hence, we now employ  a rigorous notation $(\ph^\eps, \ph_a^\eps, \mu^\eps , n^\eps, c^\eps)$ to indicate the approximate solutions.
Given that the limit passage as $\eps \to 0$ is standard  for the majority of terms, our emphasis will be directed towards the novelties that necessitated {\it ad hoc} treatment. Consequently, our primary attention will be focused on the equation involving the chemotactic variable $\pha$.

First, let us recall that 
$\ph^\eps, \ph_a^\eps, \mu^\eps , n^\eps,$ and $c^\eps$ satisfy the estimates established in the previous section with a positive constant $C$ that it is independent of $\eps$. From those, 
Banach--Alaoglu theorem entails the existence of limit functions $\ph, \ph_a, \mu , n,$ \last{and} $c$ such that, up to a not relabelled subsequence, as $\eps \to 0$,
\begin{alignat*}{2}
	\ph^\eps & \to \phT \quad && \text{weakly-star in $\L\infty V $},
    \\
    & {} && \quad \text{and weakly in $\H1 \Vp \cap\L4 {W} \cap \L2 {\Wx{2,\esp}}$}, 
	   \\
   {\ph_a^\eps} & {\to \pha \quad} && {\text{weakly in $\L 2 {H}$,}}
   \\
	\mu^\eps & \to \mu \quad && \text{weakly in $\L2 V$,}
	\\
	n^\eps & \to n \quad && \text{weakly-star in $\H1 H \cap \L\infty V \cap \L2 W$,}
    \\
    c^\eps & \to c \quad && \text{weakly-star in $\H1 H \cap \L\infty V {\cap \L2 W} \cap L^\infty(Q)$,}
\end{alignat*}
\abramob{and, in the case with a singular potential, the existence of a limit function $\xi$ such that, up to a not relabelled subsequence, as $\eps \to 0$,
\begin{alignat*}{2}
    \betaeps(\ph^\eps) & \to \xi \quad && \text{weakly in $ \L2 {\Lx \esp}$,}
\end{alignat*}}%
with exponent $\esp$ be defined as in \eqref{def:esp}.
Besides, the min-max property in \eqref{minmax:principle:c} is valid for $c$, whereas the min-max property \eqref{minmax:principle:n} is valid for $n$ in case of singular potentials.
Then, standard compactness arguments imply that, as $\eps \to 0$,
\begin{alignat*}{2}
	\ph^\eps,c^\eps,n^\eps  & \to \phT, c,n \quad &&\text{strongly in $\C0 {\Hx {1-\eta}} \cap \L2 V$} \quad  \text{for every $\eta >0$},
\end{alignat*}
and almost everywhere in $Q$.
\last{From the strong convergence of $\ph^\eps$ and the pointwise convergences at disposal, it is a standard matter to recover in the limit the inclusion $\xi \in \beta(\ph)$ $a.e.$ in $Q$}.
Now, to pass to the limit in the nonlinear terms involving $\ph_a^\eps$ in \eqref{eqn:reglim:3}, also strong convergence of $\ph_a^\eps$ to $\pha$ has to be shown. This can be achieved upon combining \eqref{rig:est:6} with  some information on the gradient  $\nabla \ph_a^\eps$ and the Aubin--Lions theorem.
Thus, from the aforementioned bounds and the interpolation $\L\infty H \cap \L2 V \emb L^{\frac {2(d+2)}d}(Q)$, we infer that
\begin{align*}
		\norma{\nabla c^\eps}_{L^{\frac {2(d+2)}d}(Q)} \leq C.
\end{align*}
Besides, from the above estimates, we infer that $\norma{T_{\eps,\eps^{-1}}(\ph_a^\eps)}_{L^2(Q)} \leq C$ which yields, using the above bound, that 
\begin{align*}
	\norma{T_{\eps,\eps^{-1}}(\ph_a^\eps) \nabla c^\eps}_{L^{\frac {d+2}{d+1}}(Q)}\leq C. 
\end{align*}
It is worth noticing that, for $d\in \{2,3\}$,  $1<\frac 54 \leq \frac {d+2}{d+1} \leq \frac 43$.
Combining this latter with the bound \eqref{rig:est:5} and the identity \eqref{tgrad}, we obtain that
\begin{align}
    \label{nabla:pha}
	\norma{\nabla \ph_a^\eps}_{L^{\frac {d+2}{d+1}}(Q)} \leq C.
\end{align}
Therefore,  \eqref{rig:est:6} and \eqref{nabla:pha}, along with the generalized Aubin--Lions theorem in the form \cite[Cor. 4, Sec. 8]{simon}, produce, as $\eps \to 0$, 
\begin{align*}
	& \ph_a^\eps \to \ph_a \quad \text{strongly in ${\L{\frac {d+2}{d+1}} {\Lx q} }$},
	\quad 
    q\in \big[1,\tfrac{d(d+2)}{d^2-2}\big).
\end{align*}
The range of exponents mentioned above, for which $2< \tfrac {15}7\leq\tfrac{d(d+2)}{d^2-2}\leq 4$, is actually  not so crucial as the above strong convergence allows us to infer that $\ph_a^\eps \to \ph_a$ also almost everywhere in $Q$, in particular.
Thus, in view of the previous bounds along with Vitali's theorem, as $\eps \to 0$,
\begin{align*}
	\ph_a^\eps \to \ph_a \quad \text{strongly in ${L^p(Q)}$,  $p<2$}.
\end{align*}

Upon combining the properties in \ref{ass:mobilities:4} with the above strong and almost everywhere convergences, we find that, as $\eps \to 0$,
\begin{align*}
    \mobm(\ph^\eps,\ph_a^\eps,n^\eps) \to \mobm(\ph,\ph_a,n),
    \quad
    \mobn(\ph_a^\eps,c^\eps) \to \mobn(\ph_a,c) 
    \quad \text{strongly in $L^q(Q)$,  $q \in [1,\infty)$},
\end{align*}
and $a. e.$ in $Q$.
Next, we  consider $w \in \Wn$, multiply \eqref{eqn:reglim:3} by a function $\varrho \in C_c^{\infty}([0,T))$ \abramob{and integrate in time between $0$ and $T$, obtaining}, after integration by parts, that
%
\begin{align}
    & \non -\int_{Q} \andre{\varrho'} \ph_a^{\eps} w
    +\andre{\intQ} \varrho\, \mobn(\phi_{a}^{\eps},c^{\eps})\nabla \phi_{a}^{\eps}\cdot \nabla w
    \andre{-}\chi_a\andre{\intQ} \varrho\, T_{\epsilon,\epsilon^{-1}}(\phi_{a}^{\eps})\mobn(\phi_a^{{\eps}},c^{\eps})\nabla c^{\eps}\cdot \nabla w\\
    &\label{rhophia} \quad =
    \varrho(0)\int_{\Omega}\ph_a^{\eps}(0)w
    + \andre{\intQ} \varrho \,\Sa\last{(\ph^\eps,\ph_a^\eps,c^\eps)} w.
\end{align}
Then, let us show how to pass to the limit in all the delicate terms. 
The first \abramob{terms} on the left\last{-hand} \abramob{and} \last{right-hand} \abramob{sides} readily pass to the limit by using the above \abramob{strong} convergence \abramob{for $\phi_{a}^{\eps}$} and that \abramob{$\Wn \emb \Lx \infty$}, \abramob{considering also that $\varrho \in C_c^{\infty}([0,T))$}. 
As the second integral is concerned
{we first notice that, combining the pointwise convergences above with \ref{ass:mobilities:4} and the Lebesgue dominated convergence theorem, $\mobn(\phi_{a}^{\eps},c^{\eps})\nabla w\to \mobn(\phi_{a},c)\nabla w$ strongly in $L^q(0,T;\Lx{d+2})$, $q \in [1,\infty)$.
Then, it readily follows that,} as $\eps \to 0$,
\begin{align*}
    & {\intQ} {\varrho}\, \mobn(\phi_{a}^{\eps},c^{\eps})\nabla \phi_{a}^{\eps}\cdot \nabla w
    \to 
    {\intQ} {\varrho}\, \mobn(\phi_{a},c)\nabla \phi_{a}\cdot \nabla w,
    \quad 
   w \in {\Wn}.
\end{align*}
We now move to the third integral and notice that, as $\eps \to 0$,
\begin{align*}
    & \chi_a \intQ  \varrho \, T_{\eps,\eps^{-1}}(\ph^\eps_a) \mobn(\ph_a^\eps,c^\eps) \nabla c^\eps \cdot \nabla w
    \to 
    \chi_a \intQ  \varrho \, \pha \mobn(\pha,c) \nabla c \cdot \nabla w,
    \quad 
   w \in {\Wn}.
\end{align*}
Actually, to pass to the limit, using the weak-strong convergence principle, it suffices that $ w \in W^{1, 3+ \gamma}(\Omega)$ for some  $\gamma>0$. 
In fact, we can combine the following properties: $T_{\eps,\eps^{-1}}(\ph^\eps_a)$ converges strongly in $L^{p}(Q)$ for any $p<2$, $\mobn(\ph_a^\eps,c^\eps)$ converges pointwise and strongly in any $L^p(Q)$, whereas $\nabla c^\eps$, due to \eqref{rig:estreg:1} and \eqref{rig:estreg:delta}, converges weakly in $\L{q_1} {\Lx{q_2}}$ for  some $q_1>2, q_2<6$, so that it is enough that $ \nabla w\in \Lx{3+ \gamma}$ for some $\gamma>0$. 
The other integrals can be treated arguing straightforwardly.
In particular, the last one, which is quadratic in $\ph_a^\eps$ due to \ref{def:Sasm} can be dealt with by using the pointwise convergence of $\ph_a^\eps$ and the dominated convergence theorem.
Thus, letting $\eps \to 0$ in \eqref{rhophia} leads to 
\begin{align*}
    & -\int_{Q} \andre{\varrho'} \ph_a w
    +\andre{\intQ} \varrho\, \mobn(\phi_{a},c)\nabla \phi_{a}\cdot \nabla w
    \andre{-}\chi_a\andre{\intQ} \varrho\, \phi_{a}\mobn(\pha,c)\nabla c\cdot \nabla w\\
    &\quad =
    \varrho(0)\iO \ph_a(0) w 
    + \andre{\intQ} \varrho \,\Sa\last{(\ph,\ph_a,c)}  w,
    \quad w \in \Wn, \, \varrho \in C_c^{\infty}([0,T)),
\end{align*}
\abramob{which, thanks to the application of Lemma \ref{LEM:dist}, gives \eqref{wf:3} for almost every $t\in (0,T)$.}
Thus, we can pass to the limit as $\eps \to 0$ in \eqref{eqn:reglim:1}--\eqref{eqn:reglim:5} to obtain \eqref{wf:1}--\eqref{wf:5}, \abramob{fulfilling the initial conditions \eqref{wf:initial:1}}.

Let us then show how to use \eqref{rig:est:6} to infer some continuity \abramob{in time} property of the variable $\ph_a^\eps$ \abramob{and recover the initial condition in the sense of \eqref{wf:initial:2}}. The strategy is largely inspired by \cite{RSchS} so that we just briefly repeat the main argument \abramob{adapting the \last{technique} to our framework}.
\abramob{Starting from \eqref{eqn:reglim:3}} using similar computations as above and that $\Wn \emb \Lx\infty$, we have
\begin{align*}
    \norma{\dt \ph_a^\eps}_{(\Wx{1,4})^*}
   &  \leq 
   C\norma{(T_{\epsilon,\epsilon^{-1}}(\pha)) ^{1/2}}_{{4}}
   \norma{(T_{\epsilon,\epsilon^{-1}}(\pha)) ^{1/2}\nabla (\mathcal{E}_{\eps,\eps^{-1}}^{'}(\pha) -\chi_a c)}
   + C
   \\ & \quad
   + C\norma{\abramob{(\ph_a^\eps)_+}}^2
	\leq f_1^\eps+C+f_2^\eps\last{,}
\end{align*}
\abramob{where we also used \eqref{rig:estreg:1bis},} \last{and it holds that} $\norma{f_1^\eps}_{L^{\frac 43}(0,T)} \leq C $ for a positive constant independent of $\eps$.
Next, {\abramob{for $r \geq 0$}}, we set $\Phi(r) = r {{\cal E}'}_{\eps,\eps^{-1}}(e + \abramob{\sqrt{r}}) $, observing that $\Phi$ is convex \abramob{and} increasing.
Applying $\Phi$ to the estimate above and integrating in time leads us to
\begin{align} \non
  \int_0^T \Phi\big( \| \dt \ph_a^\eps\|_{(\Wx{1,4})^*)} \big)
   & \leq \int_0^T \Phi(f_{1}^\eps + C + f_{2}^\eps)\\
 \non
   & \leq C \int_0^T \Phi(f_{1}^\eps) + C + C \int_0^T \Phi(f_{2}^\eps)\\
 \label{const:expl}
   & \leq C + C \int_0^T \Phi\big(C \| \ph_a^\eps\|^2 \big) \leq C,
\end{align}
\abramob{where in the last step we used \eqref{apriorieps3},} and the latter $C>0$ represents a \last{computable} constant that relies solely on the known data associated with the problem.
Consider now $0\leq \tau < t \leq T$ and notice that
\begin{equation*} 
 \frac{\| \ph_a^\eps(t) - \ph_a^\eps(\tau) \|_{(\Wx{1,4})^*}}{|t-\tau|}
   \leq \int_\tau^t \frac1{|t-\tau|} \| \dt \ph_a^\eps(r) \|_{(\Wx{1,4})^*} \, {\rm dr}.
\end{equation*}
Using that $\Phi$ is nondecreasing and convex, and applying
Jensen's inequality, we obtain
\begin{align*} \non
  \Phi\bigg( \frac{\| \ph_a^\eps(t) - \ph_a^\eps(\tau) \|_{ (\Wx{1,4})^*}}{|t-\tau|} \bigg)
   & \leq \Phi\bigg( \int_\tau^t \frac1{|t-\tau|} \| \dt \ph_a^\eps(r) \|_{(\Wx{1,4})^*} \, {\rm dr} \bigg)\\
 \non  
  & \leq  \int_\tau^t \frac1{|t-\tau|} \Phi\big( \| \dt \ph_a^\eps(r) \|_{(\Wx{1,4})^*} \big) \, {\rm dr}\\
  & \leq \frac1{|t-\tau|} \int_0^T \Phi\big( \| \dt \ph_a^\eps(r) \|_{(\Wx{1,4})^*} \big) \, {\rm dr}
    \leq \frac C{|t-\tau|},
\end{align*}
with the same $C>0$ as in \eqref{const:expl}.
Subsequently, using the strict monotonicity of $\Phi$ once more, we infer that
\begin{equation*} 
  \frac{\| \ph_a^\eps(t) - \ph_a^\eps(\tau) \|_{(\Wx{1,4})^*}}{|t-\tau|} 
   \leq \Phi^{-1} \Big(\frac{C}{|t-\tau|} \Big),
\end{equation*}
whence
\begin{equation*} 
  \| \ph_a^\eps(t) - \ph_a^\eps(\tau) \|_{(\Wx{1,4})^*}
   \leq | t - \tau |\Phi^{-1}\Big(\frac C{|t-\tau|}\Big).
\end{equation*}
Upon recognizing that $\Phi^{-1}$ is strictly sublinear at infinity, 
verified through a direct check, we obtain an equicontinuity property: 
for any $\badeps > 0$, there exists \last{$\delta > 0$} such that for every $\eps \in (0,1)$ and every $0 \leq \tau < t \leq T$ 
\begin{equation*} 
  |t - \tau| < \delta 
  \quad \Rightarrow \quad 
  \| \ph_a^\eps(t) - \ph_a^\eps(\tau) \|_{(\Wx{1,4})^*} < \badeps.
\end{equation*}
Due to  \eqref{rig:estreg:1}, we also have
\begin{equation*} 
  \| \ph_a^\eps \|_{L^\infty(0,T;L^1(\Omega))} \leq C.
\end{equation*}
Hence, observing that $L^1(\Omega) \last{\subset\subset} {(\Wx{1,4})^*}$ with compact embedding,
if we take as $\calZ$ a generic (reflexive) Banach space  such that 
\begin{equation*}
  L^1(\Omega) \subset \subset \calZ \subset {(\Wx{1,4})^*)},
\end{equation*}
using some interpolation one checks that Ascoli's theorem
can be applied to the sequence $\{\ph_a^\eps\}$ in the space 
$\C0 {\cal Z}$ so to obtain
\begin{align*} 
  \ph_a^\eps \to \pha \quad \text{strongly in }\,C^0([0,T];{\calZ})
\end{align*}
and, in particular, in $\C0 {(\Wx{1,4})^*}$.

%
%
%

We finalize the proof by demonstrating how the remaining regularity properties mentioned in the theorem can be achieved.
From \eqref{rig:estreg:1bis}, letting $\eps \to 0$, we infer by semicontinuity of norms that 
\begin{align*}
    \text{$(\pha)_- =0 \quad a. e.$ in $Q$, \, so that \, $\pha\geq 0 \quad a.e.$ in $Q$,}
\end{align*}
which proves \eqref{reg:weak:4}.
Using this, along with  the bound at disposal, allow us to employ the dominate convergence theorem to infer that 
\begin{align*}
    \kappa_\infty \iO \ph_a^2 \log (\ph_a)
    =
    \kappa_\infty \iO (\ph_a)_+^2 \log (\ph_a)
    =
   \kappa_\infty 
   \lim_{\eps \to 0}
    \iO (\ph_a^\eps)_+^2 \mathcal{E}_{\epsilon,\epsilon^{-1}}^{\prime} (\ph_a^\eps)\leq C
\end{align*}
so that \eqref{reg:weak:6} follows.
Let us incidentally notice that this entails that $\pha\in\L2H$, in particular.
%
This concludes the proof.

\color{black}

\section{Regularity Results}
\setcounter{equation}{0}

\label{SEC:REG}
\an{In this section, our focus turns to exploring the regularity characteristics of weak solutions.}
The proofs we are going to derive are related to  the regularity of weak solutions and
rely on suitable  higher-order a-priori estimates. \an{To derive these estimations, we follow a formal approach, directly handling the original system \Sys, thus avoiding any unnecessary additional technicalities. Nevertheless, a rigorous approach would entail incorporating the previously introduced approximation.}
In this direction, let us also remark that similar results for a comparable system have been obtained in \cite{RSchS}.
Finally, let us recall that \an{the} results to follow are restricted to the two dimensional setting \andre{$d=2$}, where better Sobolev embedding estimates hold.

\begin{proof}[Proof of Theorem~\ref{THM:REG:LOC}]

\andre{To begin with, we}
test equation \eqref{system:4} by $\pha $ to infer that
\begin{align}\label{reg:proof:1}
	\frac 12 \frac d {dt} \norma{\pha}^2
	+ \andre{m_0} \norma{\nabla \pha}^2
	+ \kappa_\infty \norma{\pha}^3_3
	\leq C \norma{\pha}^2
	+ \chi_a \iO \pha \nabla c \cdot \nabla \pha.
\end{align}
To handle the last term on the \rhs\ we use the Young and \last{Lady\v zhenskaya} inequalities (for $d=2$) to obtain that
\abramob{
\begin{align}
	& \non \chi_a \iO \pha \nabla c \cdot \nabla \pha
	 \\ & \quad \non
  \leq \chi_a \norma{\pha}_4 \norma{\nabla c}_4 \norma{\nabla \pha}
	 \leq
	 \chi_a  C_{1,\Omega}^2\norma{\pha}^{1/2}\norma{\pha}^{1/2}_V \norma{c}_V^{1/2}\norma{c}_{\Hx2}^{1/2} \norma{\nabla \pha}
	 \\ & \non \quad \leq
	 \chi_aC_{1,\Omega}^2 C_{2,\Omega}^{\frac{1}{2}}\norma{c}_V^{1/2}\norma{\pha}^{1/2} (\norma{\pha}^{1/2}+ \norma{\nabla \pha}^{1/2}) (\norma{\Delta c}_{}^{1/2}+1) \norma{\nabla \pha}
	  \\ & \non \quad \leq
	 \chi_aC_{\Omega} \norma{c}_V^{1/2}\norma{\pha}\norma{\nabla \pha}+\chi_aC_{\Omega} \norma{c}_V^{1/2}\norma{\pha}\norma{\Delta c}_{}^{1/2}\norma{\nabla \pha}
     \\ & \label{reg:proof:2} \qquad + \chi_aC_{\Omega} \norma{c}_V^{1/2}\norma{\pha}^{\frac{1}{2}}\norma{\nabla \pha}^{\frac{3}{2}}+\chi_aC_{\Omega} \norma{c}_V^{1/2}\norma{\pha}^{\frac{1}{2}}\norma{\Delta c}_{}^{1/2}\norma{\nabla \pha}^{\frac{3}{2}},
\end{align}
where $C_{1,\Omega}$ is the positive constant in the \last{Lady\v zhenskaya} inequality and $C_{2,\Omega}$ is the positive constant in the elliptic regularity estimates, both depending \last{only} on the geometry of the domain, and $C_{\Omega}:=C_{1,\Omega}^2 C_{2,\Omega}^{\frac{1}{2}}$.
We now observe that some a-priori estimates in \eqref{rig:estreg:1} remain valid also in the limit, as $\eps \to 0$, by \last{the} weak convergence and weak lower semicontinuity of norms. In particular, we have that
\begin{equation}
    \label{apriorilimit}
        \norma{\pha}_{\L\infty {\Lx1}}+\norma{n}_{\L\infty V}
	+ \norma{c}_{\H1 H \cap \L\infty V}
	\leq C_0,
\end{equation}
where $C_0:=C(\phi^0,\phi_a^0,n^0,c^0\an{)}+C$ \last{arose} from \eqref{apriorieps3} \last{and} depends only on proper norms of the initial conditions and on $|\Omega|$.}
Next,  we test \eqref{system:6} by $-\Delta c$ to infer 
\begin{align*}
	\norma{\Delta c}^2
	= 
	- \iO (\abramob{\dt c}+ \chi_a \pha +\Sc) \Delta c.	 
\end{align*}
Recalling Theorem~\ref{THM:WEAK:LOC}, \abramob{ \eqref{def:Sc} and using the Cauchy--Schwarz and Young inequalities lead us to
\begin{align*}
	\norma{\Delta c}^2 \leq (\chi_a^2+1)\norma{\pha}^2+\norma{\dt c}^2+\norma{n}^2+\frac{3}{4}\norma{\Delta c}^2,
\end{align*}
which gives, employing \eqref{apriorilimit}, that
\begin{align}
        \label{reg:proof:3}
	\norma{\Delta c}^2 \leq 4(\chi_a^2+1)\norma{\pha}^2+4C_0^2.
\end{align}
Using \eqref{apriorilimit} and \eqref{reg:proof:3} in \eqref{reg:proof:2}, \last{the} Young inequality and keeping track of the exact values of some constants when needed, we obtain that, upon introducing a positive constant $\iota <1$, \last{that}
\begin{align*}
	 \chi_a \iO \pha \nabla c \cdot \nabla \pha
	 &\leq
	 m_0\iota \norma{\nabla \pha}^2+C\norma{\pha}^2+\frac{\chi_a^2C_{\Omega}^2C_0}{m_0\iota}\norma{\pha}^2\norma{\Delta c}
    \\ & \quad 
    +\frac{27C_0^2C_{\Omega}^4\chi_a^4}{4m_0^3\iota^3}\norma{\pha}^2\norma{\Delta c}^2
	  \\ &  \leq
	 m_0\iota \norma{\nabla \pha}^2+C\norma{\pha}^2+C+\left(\frac{27C_0^2C_{\Omega}^4\chi_a^4(\chi_a^2+1)}{m_0^3\iota^3}+\frac{{\last{\eps_0}}}{C_0}\right)\norma{\pha}^4,
\end{align*}
where ${\last{\eps_0}}$ is an arbitrarily small positive constant. We then observe that the following interpolation \last{estimate} is valid
\begin{equation}
    \label{reg:proof:int}
    \norma{\pha}^4=\left(\iO \phi_a^{\frac{3}{2}}\phi_a^{\frac{1}{2}}\right)^2\leq \norma{\pha}^3_3\norma{\pha}_1\leq C_0\norma{\pha}^3_3\,
\end{equation}
where in the last step we used \eqref{apriorilimit}. Collecting the previous results in \eqref{reg:proof:1}, we end up with
\begin{align}
       & \non
	\frac 12 \frac d {dt} \norma{\pha}^2
	+ (m_0-\iota) \norma{\nabla \pha}^2
	+ \left(\kappa_\infty-{\last{\eps_0}}-\frac{27C_0^3C_{\Omega}^4\chi_a^4(\chi_a^2+1)}{m_0^3\iota^3} \right)\norma{\pha}^3_3
	 \\ & \quad  \label{reg:proof:4}
 \leq C \norma{\pha}^2
	+ C.
\end{align}
Defining the constant
\begin{align}
    \label{def:cbar}
    \ov{C}:=\frac{(\kappa_\infty-{\last{\eps_0}})m_0^3\iota^3}{27C_0^3C_{\Omega}^4},
\end{align}
which depends only on the parameters $\kappa_\infty,m_0$, on the initial conditions through $C_0$ and on the domain through $C_{\Omega}$, in the smallness hypothesis \eqref{smallness} we may integrate \eqref{reg:proof:4} over time, use condition \eqref{ini:reg:strong} on the initial data and the
Gronwall's lemma to get 
\begin{align*}
	\norma{\pha}_{\L\infty H \cap \L2 V \last{\cap L^3(Q)}}
	\leq C.
\end{align*}
Moreover, squaring the inequality \eqref{reg:proof:3}, using \eqref{reg:proof:int} and \eqref{reg:proof:4}, employing also 
elliptic regularity theory we obtain that 
\begin{align*}
	\norma{c }_{\L4 {\Hx2}}
	\leq C.
\end{align*}
}%
\an{From this \andre{it is} a standard matter to derive from a comparison argument in \eqref{system:6} that}
\begin{align*}
	\norma{\dt \pha}_{\L2 \Vp}
	\leq C
\end{align*}
\andre{concluding} the proof.

%

\end{proof}

\begin{remark}\label{RMK:CHKS}
Unfortunately, even \andre{when the mobility $\mobn$ is constant, e.g., $\mobn\equiv 1$,} \an{we are unable to establish the aforementioned theorem in three dimensions.}
Of course the crucial term is the last on the \rhs\ of \eqref{reg:proof:1}.
The same strategy fails to work due to the different Sobolev's embeddings in dimension three.
\an{Once can also notice that}\andre{, using integration by parts,}
\begin{align*}	\chi_a \iO \pha \nabla c \cdot \nabla \pha
	= \frac {\chi_a}2 \iO \nabla c \cdot \nabla (\ph_a^2)
	= -\frac {\chi_a}2 \iO \Delta  c \, \ph_a^2,
\end{align*}
but this does not help as for $\Delta c$ we just have the $L^2$-bound given in Theorem~\ref{THM:WEAK:LOC}\andre{.}
\an{This contrasts with the situation described in } \cite[Thm.~2.2]{RSchS}, where the role of the variable $c$ in the cross-diffusion term is played by an order parameter $\ph$ which solves a (singular) Cahn--Hilliard equation instead of a parabolic one.
\an{Given that the Cahn--Hilliard equation is fourth-order in space, it offers additional regularity for $\Delta \phT$, enabling the utilization of the \andre{aforementioned} argument.}
\end{remark}

\begin{proof}[Proof of Theorem \ref{THM:REG:STRONG:PREL}]
To begin with, we  start with proving more regularity for the chemotactic and nutrient variables. Consider equation \eqref{system:6} and observe that, due to the above results, it holds that $f_c:= \Sc + \chi_a \pha \in \L2 V$. Testing then \eqref{system:6} by $\Delta^2 c=- \Delta(-\Delta \an{c)}$, integrating by parts, \an{using \ref{ass:sources:3}, and} Young's inequality produces
\begin{align*}
    \frac 12 \frac d {dt} \norma{\Delta c}^2
    + \norma{\nabla \Delta c}^2
    =
    \iO \nabla f_c \cdot \nabla \Delta c
    \leq \frac 12 \norma{\nabla \Delta c}^2
    + \frac 12 \norma{\nabla f_c}^2.
\end{align*}
Thus, it readily follows after integration over time\an{,} using the second assumption on the initial condition 
in \eqref{ini:reg:strong:prel}\an{,} and elliptic regularity theory  that $c \in \L\infty {\Hx2} \cap \L2 {\Hx3}$.

Next, we move to the chemotactic variable and differentiate \eqref{system:5} with respect to time, integrate over time for an arbitrary $t \in [0,T]$,  and test the resulting equation by $\dt n$ to obtain
\begin{align*}
   \frac 12 \iO |\dt n(t)|^2
    + \int_{Q_t} |\nabla \dt n|^2
    =
   \frac 12 \iO |\dt n(0)|^2
    + \int_{Q_t} \dt (\Sn) \dt n
    + \chi_\ph \int_{Q_t} \dt \ph \, \dt n.
\end{align*}
The first term on the \rhs\ can be bounded owing to the first condition in \eqref{ini:reg:strong:prel}, whereas by Young's inequality we bound the second one as
\begin{align*}
    \int_{Q_t} \dt (\Sn) \dt n
    = \iot \< \dt n, \dt (\Sn) >
    \leq 
     \frac 14 \int_{Q_t} \norma{\dt n}^2_V
     + C \iot \norma{\dt (\Sn)}_*^2.
\end{align*}
\an{The second term on the \rhs\ can be bounded, recalling \ref{ass:sources:3} \andre{and using the assumption on $h$}, by}
\begin{align*}
    \iot \norma{\dt (\Sn)}_*^2
    \leq C (\norma{\dt \ph}_*^2 + \norma{\dt \pha}_*^2 + \andre{\norma{\dt n}^2}).
\end{align*}
Finally, for the last term\andre{,} we observe that 
\begin{align*}
    \chi_\ph \int_{Q_t} \dt \ph \, \dt n 
    = 
   \chi_\ph \iot  \<   \dt n ,\dt \ph>
   \leq
   \frac 14 {\andre{\iot}} \norma{\dt n}^2_V
   + C \iot \norma{\dt \ph}_*^2.
\end{align*}
Therefore, Gronwall's lemma yields that 
$n  \in \W{1,\infty} H \cap \H1 V$. Moreover, arguing exactly as above, we also infer that $n \in \L\infty {\Hx2} \cap \L2 {\Hx3}$ concluding the proof.
\end{proof}

\begin{proof}[Proof of Theorem \ref{THM:REG:STRONG}]
First, we differentiate equation \eqref{system:3}  with respect to time to derive 
\begin{align}   \label{eq:mu:dt}
    \dt \mu = - \Delta \dt \ph + \beta'(\ph)\dt \ph -\lambda \dt \ph
    \quad \text{in $Q$}.
\end{align}
Then, we test \eqref{system:1} by $\dt \mu$, the above \eqref{eq:mu:dt} by $\dt \ph$,
\eqref{system:4} by $\dt \pha$, and add the resulting identities to infer that
\begin{align}
   & \non  \frac 12 \frac d{dt} \norma{\nabla \mu}^2
    + \norma{\dt \ph}^2_V
    + \iO \beta'(\ph)|\dt \ph|^2
    + \norma{\dt \pha}^2
    + \frac 12 \frac d{dt} \norma{\nabla \pha}^2
     \\ & \non \quad 
   = \iO \andre{\SS(\ph,n)} \dt \mu 
   - \chi_\ph \iO \nabla  n \cdot \nabla (\dt \mu)
    + (1+ \lambda)\norma{\dt \ph}^2
    \\ & \qquad 
    - \chi_a \iO (\nabla \pha \cdot \nabla c + \pha \Delta c)\dt \pha
   + \iO \theta \andre{(\ph,c)}(\kappa_0 \pha - \kappa_\infty \ph_a^2) \dt \pha
   = \sum_{i=1}^5 \I_i.
   \label{reg:stima}
\end{align}
\an{We point \andre{out} that above, for} convenience, we also add to both sides the term $\norma{\dt \ph}^2$.
\an{Then, we estimate the integrals on the \rhs. Using integration by parts we obtain that}
\begin{align}\label{dtS:est}    \I_1 = \iO \andre{\SS(\ph,n)} \dt \mu 
    = 
    \frac d {dt}\iO \andre{\SS(\ph,n)} \mu
    - \iO \dt (\andre{\SS(\ph,n)}) \mu.
\end{align}
Now, the first term on the \rhs\ can be moved on the \lhs\ of the above identity, whereas
the other can be bounded by the Young inequality as
\begin{align*}   
    - \iO \dt (\andre{\SS(\ph,n)}) \mu
    \leq 
    \delta (\norma{\dt \ph}^2 + \norma{\andre{\dt n}}^2)
    + C_\delta \norma{\mu}^2,
\end{align*}
for every $\delta >0$, \andre{due to the properties of $\SS$ required in \ref{ass:sources:3} which entails that $\cal H$ is uniformly bounded}.
Besides, using the Young and \Holder\ inequalities, we find
\begin{align*}
    \I_3 
    & 
    \leq \delta \norma{\dt \ph}^2_V
    +C_\delta \norma{\dt \ph}^2_*,
    \\
    \I_4 & 
   \leq 
   C \norma{\nabla \pha}\norma{\nabla c}_\infty\norma{\dt \pha}
    + C \norma{\pha}_4\norma{\Delta c}_4\norma{\dt \pha}
    \\  &
    \leq \delta \norma{\dt \pha}^2
    + C_\delta \norma{\pha}_V^2\norma{c}_{\Wx{2,4}}^2 ,
    \\
   \I_5 & \leq 
   \delta \norma{\dt \pha}^2
    + C_\delta (1+ \norma{\pha}^4_V),
\end{align*}
where we also \andre{employ} the embedding $\Hx3 \emb\Wx{2,4} \emb \Lx \infty$  that \andre{entails} $t \mapsto \norma{c(t)}_{\Wx{2,4}}^2 \in L^1(0,T)$ due to Theorem \ref{THM:REG:STRONG:PREL}.
For $\I_2$ we need to integrate by parts. Thus, we integrate \eqref{reg:stima} over time and  notice that 
\begin{align*}
      \iot \I_2 
     &  = \chi_\ph \int_{Q_t} \nabla  (\dt n) \cdot \nabla  \mu
      - \chi_\ph \iO \nabla n(t) \cdot \nabla  \mu (t)
      +\chi_\ph \iO \nabla  n^0 \cdot \nabla  \mu^0
      \\ & \leq 
      \an{C \int_{Q_t} (|\nabla \dt n|^2+|\nabla \mu|^2)
       +}\frac 14 \iO |\nabla \mu(t)|^2
     + C \iO |\nabla n(t)|^2
    + C ( \norma{\nabla n^0}^2+\norma{\nabla \mu^0}^2)
\end{align*}
for an arbitrary $t \in [0,T]$. Then,  for the first term on the \rhs\ of \eqref{dtS:est} it holds that 
\begin{align*}
    - \iO \andre{\SS(\ph,n)} \mu
    & = \iO \andre{\SS(\ph,n)} (\mu - \mu_\Omega)
    - \iO \andre{\SS(\ph,n)}  \mu_\Omega
    \geq 
    - \frac 18 \norma{\nabla \mu}^2
    - |\mu_\Omega|
    -C 
   \\ &  \geq 
    - \frac 18 \norma{\nabla \mu}^2
    - c_1(\norma{\beta(\ph)}_1 +\norma{\ph}_1 )
       \\ &  \geq 
    - \frac 18 \norma{\nabla \mu}^2
    - c_1 \norma{\beta(\ph)}_1
    - c_2
\end{align*}
for computable positive constants $c_1,c_2$.
On the other hand, arguing as above by using \eqref{MZ} and testing \eqref{system:3} by $\ph-\ph_\Omega$, we infer that
\begin{align*}
    C_F \norma{\beta(\ph)}_1
    \leq C (1+ \norma{\nabla \mu})
\end{align*}
with the same constant $C_F$.

To recover the full $V$\andre{-}norm of $\pha$ we also test \eqref{system:4} by $\pha$ to get, after similar manipulations,
\begin{align*}
    \frac 12 \frac d{dt} \norma{\pha}^2
    \leq 
    C \norma{\pha}^2
    + \chi_a \norma{\pha}_6\norma{\nabla \pha}\norma{\nabla c}_3
    \leq 
    C (1+ \norma{c}_{\Hx2}^2) \norma{\pha}^2_V.
\end{align*}
Upon adding the above estimates\an{, and integrating over time,} we have
\begin{align*}
   &  \frac 12 \Big( 
    \frac 12 \norma{\nabla \mu(t)}^2
    - 2\iO \SS\big(\ph(t),\andre{n}(t)\big) \mu(t)
    + \norma{\pha(t)}^2_V
    \Big)
    \\  & \qquad 
    + (1-\delta) \iot \norma{\abramob{\dt \ph}}_V^2
    + (1-3\delta) \int_{Q_t} |\dt \pha|^2
    \\  & \quad 
    \leq 
   C \big( 
    \norma{\mu^0}^2_V
    + \norma{\ph^0}^2
    + \norma{\ph_a^0}^2_V
    + \norma{n^0}^2_V
    \big)
    \\  & \qquad 
    +  C_\delta \iot \norma{\dt \ph}^2_*
    +  C_\delta \iot (1+\norma{c}_{\Wx{2,4}}^2 
    +\norma{c}_{\Hx2}^2+ \norma{\pha}^2_V) \norma{\pha}^2_V.
\end{align*}
Furthermore, due to the above observation we have 
\begin{align*}     
    & \frac 12 \norma{\nabla \mu(t)}^2
    - 2\iO \SS\big(\ph(t),\andre{n}(t)\big) \mu(t)
    + \norma{\pha(t)}^2_V
     \\ & \quad 
     \geq 
     \frac 38 \norma{\nabla \mu (t)}^2
     - C (\norma{\nabla \mu (t)} +1)
     + \frac 12 \norma{\pha(t)}^2_V
     \\ & \quad 
     \geq
      \frac 14 \norma{\nabla \mu (t)}^2
     + \frac 12 \norma{\pha(t)}^2_V
     -C_*
\end{align*}
for a computable positive constant $C_*$. We then add to both sides the constant $C_*$, adjust $\delta\in (0,1)$ small enough, and invoke Gronwall's lemma to infer that 
\begin{align*}
    \norma{ \nabla \mu}_{\L\infty H }
    + \norma{ \pha}_{\H1 H \cap \L\infty V }
    + \norma{ \ph}_{\H1  V }
    \leq C.
\end{align*}

Next, comparison argument in \eqref{system:3} readily shows that $\mu_\Omega$ is bounded in $L^\infty(0,T)$ so that, using \Poincare's inequality \andre{we find}
\begin{align*}
        \norma{\mu}_{\L\infty V }
    \leq C.
\end{align*}

Once this is at disposal, we can read \eqref{system:3} as an elliptic equation with forcing term bounded in $\L\infty V$ and \andre{thus} obtain  \andre{that}
\begin{align*}
     \norma{\ph}_{\L\infty {\Wx {2,\esp}} }
    +  \norma{\beta(\ph)}_{\L\infty {\Lx {\esp}} }
    \leq C
\end{align*}
with $\esp$ as in the statement.
Finally, comparison in \eqref{system:1} readily entails that 
\begin{align*}
     \norma{\dt \ph}_{\L\infty \Vp }
    \leq C,
\end{align*}
and elliptic regularity in \eqref{system:3} that 
\begin{align*}
     \norma{\pha}_{\L2 {\Hx2} }
    \leq C\andre{,}
\end{align*}
completing the proof.
\end{proof}

\begin{proof}[Proof of Theorem \ref{THM:REG:SEP}]
Here, \an{we pursue a strategy akin to the one employed in the proof of} \cite[Thm. 2.4]{RSchS}. Indeed, despite our model is of multiphase nature, the current result is mainly  focused on the Cahn--Hilliard structure so that the main ideas \an{can be extended to the current scenario}.

First, using \eqref{growth}, we readily obtain that
\begin{align*}
     \norma{\beta'(\ph)}_{\L\infty {\Lx \esp} }
    \leq C
    \quad 
    \text{for any $\esp \in [1,\infty)$.}
\end{align*}
Then, we differentiate \eqref{system:1}\andre{, using \eqref{def:calS},} with respect to time to obtain
\andre{
\begin{align*}
    \partial_{tt} \ph 
    -\Delta \dt \mu 
    =  \dt (\SS(\ph,n))
     =
     - m \dt \ph 
     + \partial_\ph {\cal H}(\ph,n) \dt \ph 
     + \partial_n  {\cal H}(\ph,n) \dt  n
     \quad 
     \an{\text{in $Q$.}}
\end{align*}
Now, as concerns partial derivatives of $\cal H$, we have
\begin{align*}
    \partial_\ph {\cal H}(\ph,n) & = \begin{cases}
    (h(n)-\d_n)_+ h'(\ph) \qquad \text{when the potential is  regular,}
    \\
     (n-\d_n)_+ h'(\ph) \quad \hspace{1cm}\text{when the potential is singular,}
    \end{cases}
    \\ 
    \partial_n {\cal H}(\ph,n) & = \begin{cases}
    h'(n)  h(\ph) \, \chi_{\{h(n)>\d_n\}} \quad
    \hspace{1mm} \text{when the potential is regular,}
    \\
       h(\ph)\,\chi_{\{n>\d_n\}} \quad \hspace{1.45cm}\text{when the potential is singular,}
    \end{cases}
\end{align*}
and notice that, due to the assumption on $h$, $\norma{ \partial_\ph {\cal H}(\ph,n) +  \partial_n {\cal H}(\ph,n)}_{L^\infty(Q)} \leq C$ for some positive constant $C$.
}
Testing it by $\dt \ph$ leads us to
\begin{align*}
    \frac 12 \frac d {dt} \norma{\dt \ph }^2
    + \iO \nabla \dt \mu \cdot \nabla \dt \ph
    \leq 
    C (\norma{\dt \ph }^2+\andre{\norma{\dt n }^2}).
\end{align*}
For the second term on the \lhs, we notice that
\begin{align*}
    \iO \nabla \dt \mu \cdot \nabla \dt \ph
    = \norma{\Delta \dt \ph}^2
    - \iO \beta'(\ph) \dt \ph \Delta \dt \ph
     + \lambda \norma{\nabla  \dt \ph}^2.
\end{align*}
Combining the above lines, we infer that
\begin{align*}
    \frac 12 \frac d {dt} \norma{\dt \ph }^2
    + \norma{\Delta \dt \ph}^2
    \leq 
    C (\norma{\dt \ph }_V^2+\norma{\dt \pha }^2)
    +\iO \beta'(\ph) \dt \ph \Delta \dt \ph.
\end{align*}
On the other hand \an{it holds that}
\begin{align*}
    \iO \beta'(\ph) \dt \ph \Delta \dt \ph
    \leq C \norma{\beta'(\ph)}_4  \norma{ \dt \ph }_4 \norma{\Delta \dt \ph}
    \leq 
    C \norma{\dt \ph}^2_V
    +\frac 12 \norma{\Delta \dt \ph}^2.
\end{align*}
Then, we integrate over time, \an{and} use that $\dt\ph(0)= \Delta \mu^0 - \chi_\ph \Delta n^0 + \SS(\ph^0,\andre{n}^0) \in H$ to find that
\begin{align*}
    \norma{\dt \ph}_{\L\infty H \cap \L2 {\Hx2}}
    \leq C.
\end{align*}

Thus, we now consider \eqref{system:1} as an elliptic equation in term of $\mu$ and observe that the forcing term is bounded in $\L\infty H \cap \L2 V$ due to Theorems \ref{THM:REG:STRONG:PREL} and \ref{THM:REG:STRONG}. It \an{then} follows that 
\begin{align*}
    \norma{\mu}_{\L\infty {\Hx2} \cap \L2 {\Hx3}}
    \leq C
\end{align*}
which also entails that, by Sobolev's embeddings,
\begin{align*}
    \norma{\mu}_{L^\infty(Q)}
    \leq C
\end{align*}
Once this regularity is proved it is a standard matter to derive the separation principle by arguing as done, e.g., in \cite{CSS, S}. Of course, if the potential is of single-well \an{type}, we can derive the separation just where the convex part of the  potential explodes
and infer \eqref{separation:LJ}.
\an{In this direction, it worth noticing that $\H1 W$  is continuously embedded in $C^0(\ov Q)$ so that the separation property holds for every point of the parabolic cylinder and not just almost everywhere.}

\end{proof}

\begin{proof}[Proof of Theorem \ref{THM:UNIQUENESS}]
For the uniqueness result, we proceed following the same lines of argument employed in the proof of  \cite[Thm. 2.8]{RSchS} having care to handle the additional equation involving $n$ and $c$. First, we set the notation
\begin{align*}
    \ph& := \ph_1-\ph_2,
    \quad 
    \pha:= \ph_{a,1}-\ph_{a,2},
    \quad 
    \mu:= \mu_1-\mu_2,
    \quad 
    n:= n_1-n_2,
    \quad 
    c:= c_1-c_2,
    \\
    {\cal S}^i& :=\Sv(\ph_i, \andre{n}_{i}),
    \quad   
     {\cal S}_{n}^i :=\Sn(\ph_i, \ph_{a,i}, n_i),
     \quad 
     {\cal S}_{c}^i :=\Sc(\ph_i, \ph_{a,i}, \andre{ c_i})
    \quad 
    \text{for $i=1,2.$}
\end{align*}
Recall that, due to assumption \ref{ass:sources:3}, there exists a positive constant $C$ such that 
\begin{align*}
   |{\cal S}^1-{\cal S}^2|
   & \leq C (|\ph| + |\andre{n}|),
   \quad 
   |{\cal S}_{n}^1-{\cal S}_{n}^2|
   \leq C (|\ph| + |\pha| + |n|),
   \\
   |{\cal S}_{c}^1-{\cal S}_{c}^2|
   & \leq C (|\ph| + |\pha| \andre{+ |c|}).
\end{align*}
Using this notation, we write \Sys\ for the differences to realize that
\begin{alignat}{2}
	\label{system:cd:1}
	& \dt \ph - \Delta \mu +   \chi_\ph \Delta n = {\cal S}^1- {\cal S}^2
	\quad && \text{in $Q$,}
	\\
	\label{system:cd:2}
	& \mu = -  \Delta \ph + F'(\ph_1)-F'(\ph_2) 
	\quad && \text{in $Q$,}
	\\
	\label{system:cd:3}
	& \dt \pha - \Delta \pha + \chi_a \div(\pha \nabla c_1 + \ph_{a,2} \nabla c) = \kappa_0 \pha
    - \kappa_\infty \pha (\ph_{a,1} + \ph_{a,2})
	\quad && \text{in $Q$,}
	\\
	\label{system:cd:4}
	& \dt n - \Delta n - \chi_\ph \phv = {\cal S}_n^1- {\cal S}_n^2 
	\quad && \text{in $Q$,}
	\\
	\label{system:cd:5}
	& \dt c - \Delta c  - \chi_a \pha = {\cal S}_c^1- {\cal S}_c^2
	\quad && \text{in $Q$,}	
    \\
    \label{system:cd:6}
	& \dn \ph = \dn \mu= \dn \pha = \dn n =  \dn c =0 	\qquad && \text{on $\Sigma$,}
	\\
    \non
    &  
	\ph(0) = \ph^{0}_1-\ph^{0}_2,
	\quad 
	\pha(0) = \ph_{a,1}^{0}-\ph_{a,2}^{0},
	\quad && 
	\\  & \quad
    \label{system:cd:7}
    n(0) = n^{0}_1-n^{0}_2,
	\quad 
	c(0) = c^{0}_1-c^{0}_2
	\qquad && \text{in $\Omega$.}
\end{alignat}
We observe that testing \eqref{system:cd:1} by \an{the constant} $\an{|\Omega|^{-1}}$ produces the identity
\begin{align}
    \label{cd:mean}
    \ph_\Omega'
    = \frac d{dt} \ph_\Omega
    = \frac 1 {|\Omega|}
    \iO ({\cal S}^1- {\cal S}^2) =  ({\cal S}^1- {\cal S}^2)_\Omega.
\end{align}
Multiplying the above by $\ph_\Omega$, one gets
\begin{align}\label{cd:mean:est:1}
    \frac 12 \frac d{dt} |\ph_\Omega|^2
    \leq 
    |\ph_\Omega|^2
    + C (\norma{\ph}^2+\norma{\andre{n}}^2).
\end{align}
Then, we subtract \eqref{cd:mean} from  \eqref{system:cd:1} and test the resulting identity by $\NN (\ph-\ph_\Omega)$ leading to
\begin{align}
    \non
    &\frac 12 \frac d {dt}\norma{\ph-\ph_\Omega}_*^2
    + \iO (\mu-\mu_\Omega)(\ph-\ph_\Omega)
     = 
    \iO \big({\cal S}^1- {\cal S}^2 -({\cal S}^1- {\cal S}^2)_\Omega \big) \NN(\ph-\ph_\Omega)
    \\ & \quad 
    \leq 
    C (\norma{\ph}^2+\norma{\andre{n}}^2+\norma{\ph-\ph_\Omega}^2_*),
    \label{cd:ineq:1}
\end{align}
where we also used that $ \iO \mu_\Omega (\ph-\ph_\Omega) =\mu_\Omega\iO  (\ph-\ph_\Omega) =0$.  
On the other hand, for every $\delta >0$, it holds that 
\begin{align}
     \nonumber
    C \norma{\ph}^2
    &\leq 
    C (\norma{\ph-\ph_\Omega}^2 + |\ph_\Omega|^2)
    \leq 
    C( \norma{\ph-\ph_\Omega}_V \norma{\ph-\ph_\Omega}_*
    + |\ph_\Omega|^2)
    \\ &  
    \leq \delta \norma{\nabla \ph}^2
    + C_\delta \norma{\ph-\ph_\Omega}_*^2
    + C |\ph_\Omega|^2.
    \label{cd:ph:interp}
\end{align}
Moreover, the second term on the \lhs\ of \eqref{cd:ineq:1} can be bo\an{u}nded as
\begin{align*}
    \iO (\mu-\mu_\Omega)(\ph-\ph_\Omega)
    = \norma{\nabla \ph}^2
    + \iO(F'(\ph_1)-F'(\ph_2) ) (\ph-\ph_\Omega)
\end{align*}
so that, combining with \eqref{cd:mean:est:1}, we obtain
\begin{align}
    \non
    & \frac 12 \frac d {dt} \big ( \norma{\ph-\ph_\Omega}_*^2 +|\ph_\Omega|^2 \big)
    + (1-2\delta) \norma{\nabla \ph}^2
   \\ & \quad  
   \leq 
      C_\delta \norma{\ph-\ph_\Omega}_*^2
    + C (\norma{\andre{n}}^2+|\ph_\Omega|^2)
    + \iO|\beta(\ph_1)-\beta(\ph_2) | |\ph_\Omega|,
    \non
\end{align}
where we also estimate the nonconvex contribution of the potential as 
\begin{align*}
    \iO |\pi(\ph_1)-\pi(\ph_2) | |\ph-\ph_\Omega|
    \leq \delta \norma{\nabla \ph}^2
    +  C_\delta \norma{\ph-\ph_\Omega}_*^2
    + C |\ph_\Omega|^2.
\end{align*}
Next, we test \eqref{system:cd:4} by $ n$ and \eqref{system:cd:5} by $ c$.
Adding the resulting equalities leads us to obtain 
\begin{align*}
   & \frac 1 2 \frac d {dt} \big(\norma{n}^2+\norma{c}^2\Big)
   +  \norma{\nabla n}^2
   + \norma{\nabla c}^2
   \\ & \quad 
   =  \iO ({\cal S}_n^1-{\cal S}_n^2) n
   +  \chi_\ph \iO \ph \, n
   +   \iO ({\cal S}_c^1-{\cal S}_c^2) c
   +  \chi_a \iO \pha c
   \\ & \quad \leq
   \delta \norma{\nabla \phi}^2
   + C_{\delta } \norma{\ph-\ph_\Omega}^2_*
   + C|\ph_\Omega|^2
   + C(\norma{n}^2+\norma{c}^2)
    \\ & \qquad
    +  C (\norma{\pha-(\pha)_\Omega}^2 + |(\pha)_\Omega|^2).
\end{align*}
where we used the \Holder\ and \Poincare\ inequalities, \ref{ass:sources:3}, as well as \eqref{cd:ph:interp}. Adding this to the above inequality produces
\begin{align}
    \non
    & \frac 12 \frac d {dt} \big ( \norma{\ph-\ph_\Omega}_*^2 +|\ph_\Omega|^2 + \norma{n}^2+\norma{c}^2 \big)
    + (1-3\delta) \norma{\nabla \ph}^2
     +  \norma{\nabla n}^2
   + \norma{\nabla c}^2
   \\ & \quad   \non
   \leq 
      C_\delta \norma{\ph-\ph_\Omega}_*^2
     +C( |\ph_\Omega|^2 + |(\pha)_\Omega|^2)
     +  C^* \norma{\pha-(\pha)_\Omega}^2 
    \\ & \qquad 
    + C(\norma{n}^2+\norma{c}^2)
    + \iO|\beta(\ph_1)-\beta(\ph_2) | |\ph_\Omega|,
    \label{cd:phi:est}
\end{align}
for an explicit and computable constant that we term $C^*$.

We then repeat similar  arguments as above to handle the mean value of  $\pha$. Namely, we multiply \eqref{system:cd:3} by $\an{|\Omega|^{-1}}$ obtaining
\begin{align}
    \label{cd:mean:pha}
    (\pha)_\Omega'
    = \frac d {dt} (\pha)_\Omega
    = \kappa_0 (\pha)_\Omega
    - \kappa_\infty (\ph_{a,1}^2-\ph_{a,2}^2)_\Omega.
\end{align}
Then, we test the above by $(\ph_a)_\Omega$ to infer that \begin{align*}
     \frac 12 \frac d{dt} |(\pha)_\Omega|^2
    \leq 
    \eta \norma{\pha - (\pha)_\Omega}^2
    + C_\eta (\norma{\ph_{a,1}}^2
    + \norma{\ph_{a,2}}^2+1)|(\pha)_\Omega|^2,
\end{align*}
for a positive constant $\eta$, yet to be selected. Then, we subtract  \eqref{cd:mean:pha} to \eqref{system:cd:3} and test the difference by $\NN (\pha-(\pha)_\Omega)$ to infer that
\begin{align*}
    \non
    &\frac 12 \frac d {dt}\norma{\pha-(\pha)_\Omega}_*^2
    +\norma{\pha-(\pha)_\Omega}^2
  \\ & \quad 
  \leq
    \chi_a \iO (\pha \nabla c_1 + \ph_{a,2} \nabla c)\cdot \nabla \NN(\pha-(\pha)_\Omega)
    + \kappa_0 \norma{\pha-(\pha)_\Omega}^2_*
    \\ & \qquad 
    - \kappa_\infty \iO \big(\ph_{a,1}^2 - \ph_{a,2}^2 - (\ph_{a,1}^2)_\Omega 
    + (\ph_{a,2}^2)_\Omega \big)\NN(\pha-(\pha)_\Omega).
\end{align*}
In the order, using the same computations as in the proof of \cite[Thm. 2.8]{RSchS}, we have
\begin{align*}
     & \chi_a \iO (\pha \nabla c_1 + \ph_{a,2} \nabla c)\cdot \nabla \NN(\pha-(\pha)_\Omega)
    \\ & \quad 
    \leq 
     \eta (\norma{\nabla c}^2
     +\norma{\pha-(\pha)_\Omega}^2)
     + C\eta |(\pha)_\Omega|^2
     + C_\eta (\norma{c_1}^2_{W^{2,6}(\Omega)} + \norma{c_2}^4_{6}) \norma{\pha-(\pha)_\Omega}^2_*,
\end{align*}
and 
\begin{align*}
     &
     - \kappa_\infty \iO \big(\ph_{a,1}^2 - \ph_{a,2}^2 - (\ph_{a,1}^2)_\Omega 
    + (\ph_{a,2}^2)_\Omega \big)\NN(\pha-(\pha)_\Omega)
   \\ & \quad 
   \leq 
    \eta \norma{\pha-(\pha)_\Omega}^2
    + C_\eta (\norma{c_1}^{4} + \norma{c_2}^{4} + 1) \norma{\pha-(\pha)_\Omega}^2_*
    + C |(\pha)_\Omega|^2,
\end{align*}
so that 
\begin{align}
    \non
    &\frac 12 \frac d {dt} \big( \norma{\pha-(\pha)_\Omega}_*^2
    + |(\pha)_\Omega|^2\big)
    - \eta \norma{\nabla c}^2
     + (1-2\eta  ) \norma{\pha-(\pha)_\Omega}^2
     \\ & \quad \non
    \leq
      C_\eta (\norma{c_1}^2 + \norma{c_2}^2+ 1)|(\pha)_\Omega|^2
      \\ & \qquad 
      + C_\eta (\norma{c_1}^{4} + \norma{c_2}^4_6+ \norma{c_1}^2_{W^{2,6}(\Omega)} + 1) \norma{\pha-(\pha)_\Omega}^2_*.
      \label{cd:pha:est}
\end{align}
We then choose $\delta =  1/6$ in \eqref{cd:phi:est}.
Then we add the resulting inequality tested by a positive constant ${\andre{\omega}}$ yet to be selected to find that 
\begin{align*}
        \non
    &\frac 12 \frac d {dt} \Big( 
     {\andre{\omega}} \norma{\ph-\ph_\Omega}_*^2 + {\andre{\omega}}|\ph_\Omega|^2
    + \norma{\pha-(\pha)_\Omega}_*^2
    + |(\pha)_\Omega|^2
    +{\andre{\omega}}\norma{n}^2
    +{\andre{\omega}}\norma{c}^2
    \Big)
    \\ & \qquad 
    + \frac {\andre{\omega}}2 \norma{\nabla \ph}^2
    + {\andre{\omega}} \norma{\nabla n}^2
    + ({\andre{\omega}} - \eta ) \norma{\nabla c}^2
      \\ & \qquad \non 
      + (1-2\eta  ) \norma{\pha-(\pha)_\Omega}^2
     \\ & \quad \non
    \leq
    {\andre{\omega}} C  \norma{\ph-\ph_\Omega}_*^2
      + {\andre{\omega}} C(|\ph_\Omega|^2 + |(\pha)_\Omega|^2  )
      + {\andre{\omega}} C^*\norma{\pha- (\pha)_\Omega}^2
      \\ & \qquad 
      + {\andre{\omega}} C ( \norma{n}^2+\norma{c}^2)
    + {\andre{\omega}} \iO |\beta(\ph_1)-\beta(\ph_2) | |\ph_\Omega|
      \\ & \qquad 
     + C_\eta (\norma{c_1}^2 + \norma{c_2}^2+ 1)|(\pha)_\Omega|^2
      \\ & \qquad 
      + C_\eta (\norma{c_1}^4+ \norma{c_2}^4_6+ \norma{c_1}^2_{W^{2,6}(\Omega)} + 1) \norma{\pha-(\pha)_\Omega}^2_*.
\end{align*}
Finally, we select ${\andre{\omega}}= {\andre{\omega}}^* :=  1/2  \min\{1, 1/C^*\}$ so to absorb the term involving $C^*$ on the \lhs, and highlight that all the above constants $C$ are now independent of ${\andre{\omega}}^*$ as it is fixed. 
It is clear that the only term that need to be handled is the last one on the \rhs. 
The simplest case occurs when \andre{${\cal H}={{\cal H}}(\ph,n)$} is a constant function and this has been analyzed in \cite{GGM} to deal with the Cahn--Hilliard--Oono equation, see the forthcoming Remark \ref{RMK:Oono} below.
However, the scenario for more general $h$ is more delicate and forced us to assume  \an{\eqref{separation} or \eqref{separation:LJ}, respectively. That entails} that
\begin{align*}
    \beta (\ph_1)- \beta (\ph_2)=  \ell\ph,
    \quad 
    \text{with} 
    \quad 
     \ell:=\int_0^1 \beta'(s \ph_1+ (1-s) \ph_2 ) {\rm ds}.
\end{align*}
Thus, by the \Holder\ and Young inequalities we find that
\begin{align*}
    &\iO|\beta(\ph_1)-\beta(\ph_2) | |\ph_\Omega|
     \leq \norma{\ell}\norma{\ph}|\ph_\Omega|
    \leq \norma{\ell}(\norma{\ph-\ph_\Omega} - |\ph_\Omega|) |\ph_\Omega|
    \\ & \quad 
    \leq 
    \norma{\ell}(\norma{\ph-\ph_\Omega} + |\ph_\Omega|)|\ph_\Omega|
    \leq 
   \norma{\ell}(\norma{\nabla \ph} + |\ph_\Omega|) |\ph_\Omega|
     \\ & \quad 
   \leq 
    \frac {{\andre{\omega}}^*}4  \norma{\nabla \ph}^2
    + C|\ph_\Omega|^2( \norma{\ell}^2+1)  
       \leq 
    \frac {{\andre{\omega}}^*}4  \norma{\nabla \ph}^2
    + C|\ph_\Omega|^2( \norma{\beta'(\ph_1)}^2+\norma{\beta'(\ph_2)}^2+1)  .
\end{align*}
We can now adjust $\eta$  small enough, for instance $\eta = 1/2 \min \{{\andre{\omega}}^*/2, 1/4\}$, and apply Gronwall's lemma to conclude.



\end{proof}

\begin{remark}  \label{RMK:Oono}
When $\andre{{\cal H}=\andre{{\cal H}(\ph, n)}}$ is a constant function, \eqref{cd:mean} reduces to
\begin{align*}
    \ph_\Omega' + m\ph_\Omega= 0.
\end{align*}
This is the scenario one encounters, e.g., in the Cahn--Hilliard--Oono equation (see \cite{GGM}).
Testing by $\sign \ph_\Omega$ produces 
\begin{align*}
    \frac 12 \frac d{dt} |\ph_\Omega|^2
    +|\ph_\Omega| =0.
\end{align*}
Thus,  observing that 
\begin{align*}
    \iO |\beta(\ph_1)-\beta(\ph_2) | |\ph_\Omega|
    \leq C (\norma{\beta(\ph_1)}_1+\norma{\beta(\ph_2)}_1) |\ph_\Omega|,
\end{align*}
one realizes that the above integral can be controlled by the Gronwall lemma provided to make use of  the above identity and no additional \andre{properties are} is required for $\beta$. 
\end{remark}

\vskip 6mm
\noindent{\bf Acknowledgements}
\noindent
The authors acknowledge some support 
from the MIUR-PRIN Grant 2020F3NCPX, {\it Mathematics for industry 4.0} (Math4I4) and
their affiliation to the GNAMPA (Gruppo Nazionale per l'Analisi Matematica, 
la Probabilit\`a e le loro Applicazioni) of INdAM (Isti\-tuto 
Nazionale di Alta Matematica). 
\an{Finally, AS acknowledges support by MUR, grant Dipartimento di Eccellenza 2023-2027.}

\End{document}

\bye